%% file: Dynamic_Alternative_Routing.tex
\documentclass[]{article}

\PassOptionsToPackage{%
	style	= alphabetic%
}{biblatex}

\input{DAR_preamble_biblio.tex}

\input{DAR_preamble_main.tex}
\title{\sffamily Metastability in Loss Networks\\with Dynamic Alternative Routing}
\author{\Large \sffamily Sam Olesker-Taylor}
\date{}

\begin{document}

\maketitle

\renewcommand{\abstractname}{\sffamily Abstract\smallskip}

\begin{abstract}
	Consider $N$ stations interconnected with links, each of capacity $K$, forming a complete graph.
	Calls arrive to each link at rate $\lambda$ and depart at rate $1$.
	If a call arrives to a link $x y$, connecting stations $x$ and $y$, which is at capacity, then a third station $z$ is chosen uniformly at random and the call is attempted to be routed via $z$:
		if both links $x z$ and $z y$ have spare capacity, then the call is held simultaneously on these two;
		otherwise the call is lost.
	
	We analyse an approximation of this model.
	We show rigorously that there are three phases according to the traffic intensity $\alpha \cq \lambda/K$:
		for $\alpha \in (0,\alpha_c) \cup (1,\infty)$, the system has mixing time logarithmic in the number of links $n \cq \binom N2$;
		for $\alpha \in (\alpha_c,1)$ the system has mixing time exponential in $n$, the number of links.
		Here $\alpha_c \cq \tfrac13 (5 \sqrt{10} - 13) \approx 0.937$ is an explicit critical threshold with a simple interpretation.
	We also consider allowing multiple rerouting attempts. This has little effect on the overall behaviour; it does not remove the metastability phase.
	
	Finally, we add \emph{trunk reservation}: in this, some number $\sigma$ of circuits are reserved; a rerouting attempt is only accepted if at least $\sigma+1$ circuits are available.
	We show that if $\sigma$ is chosen sufficiently large, depending only on $\alpha$, not $K$ or $n$, then the metastability phase is removed.
\end{abstract}

%
%

\small
\begin{quote}
	\begin{description}
		\item [Keywords:]
		loss network, dynamic alternative routing, metastability, mixing times
		
		\item [MSC 2020 subject classifications:]
		60K20, 60K25, 60K30; 90B15, 90B18, 90B22
	\end{description}
\end{quote}
\normalsize

\blfootnote{%
	Email: \href{mailto:sam.ot@posteo.co.uk}{sam.ot@posteo.co.uk}%
\hfill%
	Statistical Laboratory, University of Cambridge, UK%
\\%
	Website: \href{https://mathematicalsam.wordpress.com}{mathematicalsam.wordpress.com}%
\hfill
	Supported by EPSRC Doctoral Training Grant \#1885554%
}

\vfill \setcounter{tocdepth}{1}
\sffamily \boldmath
\tableofcontents
\normalfont \unboldmath
\normalsize

\vspace*{\bigskipamount}

\romannumbering

\newpage
\section{Introduction to Model, Main Results and Outline}
\label{sec:intro}

\subsection{Introduction to Model}
\label{sec:intro:model}

We analyse a popular \textit{stochastic loss network with dynamic alternative routing}.

\begin{quote}
	Suppose that $N$ nodes are linked to form a complete graph.
	Between any pair of nodes, call requests arrive at rate $\lambda$ and there is a link of capacity $K$.
	If there is a spare circuit on the link joining the end points of a call, then the call is accepted and carried by that circuit.
	Otherwise the call chooses at random a two-link path joining its end-points:
		the call is accepted on that path if both links have a spare circuit;
		otherwise it is lost.
	
	Calls release links on which they are held at rate $1$:
		those on a single-link route release this one link, while those on a two-link route release both links simultaneously.
\end{quote}
This has been the subject of substantial attention over the years; see, eg, \cite{GHK:bistability,GK:dyn-routing,K:stoc-net,K:dyn-routing} and further references at \cite{K:dyn-routing-web}. See \cite[\S3.7]{KY:stoc-net-book} for a more modern and particularly readable overview.
Dynamic alternative routing was implemented in BT's UK telecom network in 1996.

We analyse a slight approximation to this model, which was suggested to us by \textcite{K:dar:private}.
\begin{quote}
	All units of capacity, across the entire network, are released independently and at rate~$1$.
	
	Rerouted calls take one unit of capacity on a single link, chosen uniformly, as opposed to choosing a two-link path and taking one unit of capacity on each.
	The rate at which reroutings happen is doubled, to preserve the overall rate at which capacity is taken.
\end{quote}
See \S\ref{sec:motivation-comparison} for justification as to why this model well-approximates the original.
We now state precisely our model, which we denote $\DAR_n(\alpha,K)$, standing for \textit{dynamic alternative routing};
write $\lambda \cq \alpha K$.
\begin{quote}
There are $n$ links, labelled $\bra{ 1, ..., n }$, each of capacity $K$.
Given that a call is held on some link, it departs after an exponential-$1$ time.
Calls arrive to each link according to independent Poisson streams, with state-dependent rates:
	if the state of the system is $x = (x_1, ..., x_n)$,
	then the arrival rate, which is the same for each (non-full) link, is
	\[
		\lambda \rbb{ 1 + 2f(1-f) }
	\Quad{where}
		f \cq \tfrac1n \sumt[n]{j=1} \one{x_j = K},
	\]
	ie $f$ is the proportion of links which are full in $x$;
	calls do not arrive to full links.
\end{quote}
It is not immediately clear that this is the correct variable-rate arrival process; see \S\ref{sec:motivation-comparison} for justification.
We call $\lambda$ the \textit{arrival rate}, $K$ the \textit{capacity}, $\alpha = \lambda/K$ the \textit{traffic intensity} and $n$ the \textit{number of links}.
We denote the unique invariant distribution of these dynamics by $\Pi$.

We can write down the generator of this continuous-time Markov chain explicitly.
The state space is $\Omega \cq \bra{ 0, 1, ..., K }^n$.
The interpretation of $x \in \Omega$ is that link $i$ has $x_i$ calls on it, for each $i \in [n]$.
For each $j \in [n]$,
define the maps $m_{j,\pm} : \Omega \to \Omega$ by adding/subtracting $1$ from link $j$:
\[
	\rbb{ m_{j, +}(x) }_i
\cq
\begin{cases}
	(x_j + 1) \wedge K
		&\text{if}\quad
	i = j,
\\
	x_i
		&\text{if}\quad
	i \ne j;
\end{cases}
\qquad
	\rbb{ m_{j,-}(x) }_i
\cq
\begin{cases}
	(x_j - 1) \vee 0
		&\text{if}\quad
	i = j,
\\
	x_i
		&\text{if}\quad
	i \ne j.
\end{cases}
\]
The generator $\mcl$ of $\DAR_n(\alpha, K)$ is then defined by the following action on functions $g : \Omega \to \mbr$:
\begin{gather*}
	\mcl g(x)
\cq
	\sumt{j \in [n]}
	\nu(x) \rbb{ g\rbb{ m_{j,+}(x) } - g(x) }
+	\sumt{j \in [n]}
	x_j \rbb{ g\rbb{ m_{j,-}(x) } - g(x) }
\Quad{for}
	x \in \Omega,
\\
	\text{where}
\quad
	\nu(x)
\cq
	\alpha K \rbb{ 1 + 2 f(x) \rbr{ 1 - f(x) } }
\Quad{and}
	f(x)
\cq
	\tfrac1n \sumt{j \in [n]} \one{ x_j = K }.
\end{gather*}
Each non-full link has a call added at rate $\nu(x)$ and a call is removed from link $j \in [n]$ at rate $x_j$.

\medskip

We also consider two ways in which the original model can be extended.
The first adds \textit{retries}, with parameter $\rho$:
	here the idea is that instead of trying a pair of links and losing the call if either of these is full, the system tries to reroute using a pair of links $\rho$ times, stopping if a try is successful and losing the call if all $\rho$ fail.
Finally, a \textit{trunk reservation}, with parameter $\sigma$, is added:
	here instead of accepting a rerouting request if both of the links have at least one spare circuit, there must be at least $\sigma+1$ spare circuits.
This reduces to the original model if $\rho = 1$ and $\sigma = 0$.
Again, we analyse an approximation to these models, using the same adjustments as described above.

\subsection{Statement of Results}

It is well-known that the original model exhibits two phase transitions, depending on the ratio $\alpha \cq \lambda/K$, in the limit \ninf with $\alpha$ and $K$ fixed.
This is because an ODE representing the proportion of links which are full has two fixed points when $\alpha \in (\alpha_c, 1)$, for a specific $\alpha_c \in (0,1)$.
This causes the system to have metastability.
A quantitative version of metastability has never been pursued rigorously, to the best of our knowledge. Only heuristics, non-rigorous approximations and simulations have been employed.
References to past work are deferred to~\S\ref{sec:intro:prev-mot-comp:prev}.

We rigorously establish the same phase transition for our model, which is a slight simplification of the original.
%
We derive an appropriate fixed point equation in \S\ref{sec:prelim:def-alphac}, giving
\(
	\alpha_c = \tfrac13(5\sqrt{10} - 13).
\)

For a Markov chain $X = (X^t)_{t\ge0}$ with invariant distribution $\Pi$,
for $\eps \in (0,1)$,
define
\[
	\tmix(\eps)
\cq
	\inf\brb{ t \ge 0 \midb \maxt{x} \tvb{ \prb[x]{ X^t \in \cdot } - \Pi } \le \eps }.
\]
When considering the mixing time of a $\DAR_n(\alpha,K)$ system,
we write $\tmix( \, \cdot \,; \alpha, K, n)$.

\begin{mainthm}[Mixing Time for $\DAR_n(\alpha,K)$ System]
\label{res:intro:main}
	Let $\alpha_c \cq \tfrac13(5\sqrt{10} - 13)$.
	\begin{itemize}[itemsep = 0pt, topsep = \smallskipamount, label = \bcdot]
		\item 
		\emph{Fast Mixing.}
		Suppose $\alpha < \alpha_c$ or $\alpha > 1$.
		There exists a constant $C$ so that
		for
			all $\eps \in (0,1)$
		and
			all $K$ and $n$ sufficiently large,
		we have
		\(
			\tmix(\eps; \alpha, K, n) \le C \log n.
		\)
		
		\item 
		\emph{Slow Mixing.}
		Suppose $\alpha_c < \alpha < 1$.
		There exists a positive constant $c$ so that
		for
			all $\eps \in (0,\tfrac12)$
		and
			all $K$ and $n$ sufficiently large,
		we have
		\(
			\tmix(\eps; \alpha, K, n) \ge e^{cn}.
		\)
	\end{itemize}
\end{mainthm}

%
%

\begin{Proof}[Proof References]
\renewcommand{\qedsymbol}{\ensuremath{\triangle}}
The fast mixing part is proved in \cref{res:low,res:high}, for $\alpha < \alpha_c$ and $\alpha > 1$, respectively.
The slow mixing part is proved in \cref{res:slow}.
\end{Proof}

\medskip

From a network engineering point of view, metastability is a highly undesirable property.
As such, one wishes to adapt the model so as to remove this metastability.
We consider two extensions alluded to earlier:
	the first allows \textit{retries}
while
	the second adds \textit{trunk reservation}.

\begin{itemize}[itemsep = 0pt, topsep = \smallskipamount, label = \ensuremath{\bcdot}]
	\item 
	For \cref{res:intro:retry}, recall that $\rho$ is the number of tries, ie $\rho-1$ retries.
	We define $\alpha_c(\rho)$ analogously to $\alpha_c$, but for $\rho$ tries, so $\alpha_c = \alpha_c(1)$; see \S\ref{sec:prelim:def-alphac}
		and \S\ref{sec:retries}.
%
	Metastability still exists;
	in fact $\rho \mapsto \alpha_c(\rho) : \mbn \to (0,1)$ is a decreasing map, so in some sense it gets worse.
	
	\item 
	For \cref{res:intro:tr}, recall that $\sigma$ is the number of circuits reserved; we assume that $\sigma$ is sufficiently large in terms only of $\alpha$, independent of both $K$ and $n$.
	Metastability is removed by reserving this `small' number of circuits: there is fast mixing for all $\alpha < 1$.
\end{itemize}

\noindent
These are studied in \S\ref{sec:retries} and \S\ref{sec:tr}, respectively, where the following statements are proved.
We add an extra parameter $\rho$ or $\sigma$ to $\tmix$ to indicate the number of retries or circuits reserved, respectively.

\begin{mainthm}[Mixing Time with Retries]
\label{res:intro:retry}
	For all $\rho \in \mbn$, there exists a unique constant $\alpha_c(\rho)$
	with the following properties.
	\begin{itemize}[itemsep = 0pt, topsep = \smallskipamount, label = \ensuremath{\bcdot}]
		\item 
		\emph{Fast Mixing.}
		Suppose $\alpha < \alpha_c(\rho)$ or $\alpha > 1$.
		There exists a constant $C$ so that,
		for
			all $\eps \in (0,1)$
		and
			all $K$ and $n$ sufficiently large,
		we have
		\(
			\tmix(\eps; \alpha, K, n; \rho) \le C \log n.
		\)
		
		\item 
		\emph{Slow Mixing.}
		Suppose $\alpha_c(\rho) < \alpha < 1$.
		There exists a positive constant $c$ so that,
		for
			all $\eps \in (0,\tfrac12)$
		and
			all $K$ and $n$ sufficiently large,
		we have
		\(
			\tmix(\eps; \alpha, K, n; \rho) \ge e^{cn}.
		\)
	\end{itemize}
	Further, the map $\rho \mapsto \alpha_c(\rho) : \mbn \to \mbr$ is strictly decreasing and satisfies $\alpha_c(\rho) \in (0,1)$ for all $\rho \in \mbn$.
\end{mainthm}

The similarity, almost equivalence, between \cref{res:intro:main,res:intro:retry} is not unexpected once one realises the underlying reason behind the metastability.
This is described in \S\ref{sec:outline-criticalthresh:crit}.


\begin{mainthm}[Mixing Time with Trunk Reservation]
\label{res:intro:tr}
	For all $\alpha \in (0,\infty) \setminus \bra{1}$,
	there exists a $\sigma_*(\alpha)$ so that if $\sigma \ge \sigma_*(\alpha)$ then
	for
		all $\eps \in (0,1)$
	and
		all $K$ and $n$ sufficiently large,
	we have
	$\tmix(\eps) \le 60 \log n$.
	
	Furthermore, with these parameters, when $\alpha < 1$, the proportion of links which are full, in equilibrium, may be made as small as desired by taking $K$ sufficiently large (independently of $n$).
\end{mainthm}

An important observation is that the trunk reservation parameter $\sigma_*$ depends only on $\alpha$, not on $K$ or $n$. This means that when the system is scaled up (ie $K$ and $n$ increase), the number of reserved links does not need to increase.
It is somewhat remarkable that reserving this small number of links (not growing as $\Kinf$) removes metastability.

The underlying structure of the proof of \cref{res:intro:tr} follows that of \cref{res:intro:main}, but the details are significantly different.
The main point is to show how the additional reservation of circuits can be used to extend the fast mixing with slow arrivals regime all the way to the entirety of $\alpha \in (0,1)$.
This requires significant additional analysis.

\subsection{Motivation for Critical $\alpha_c$ and Outline of Proof}
\label{sec:outline-criticalthresh}

We next explain the underlying reasons \emph{why} such a critical $\alpha_c$ appears.
The intuition is based on the original model, which has a nice interpretation in terms of rerouting, not our approximation via a variable-rate Poisson arrivals with rate-$1$ departures.
We use this intuition to give a brief outline of the proof.
The actual proof uses the approximation for technical reasons only.

\subsubsection{Critical Threshold}
\label{sec:outline-criticalthresh:crit}

The maximal service rate, ie rate at which capacity is released, of a given link is $K$.
Ignoring for the moment reroutings, it is clear that there should be significantly different behaviour for $\alpha > 1$ compared with $\alpha < 1$.
Indeed, calls always arrive faster than they depart if $\alpha > 1$, giving rise to a bias towards adding calls.
Contrastingly, 
Further, if $\alpha < 1$ and the number of calls on a given link is larger than $\alpha K$, then there is bias towards removing calls. Thus the equilibrium probability hat the link is full tends to $0$ as $\Kinf$.
If the system starts empty, then large deviations results should imply that the proportion of full links does not become significant high for a long time.

The existence of $\alpha_c$, however, is less obvious.
Key to understanding the underlying reason for its existence is realising that \emph{each} rerouted call holds \emph{two} units of capacity.
Thus, while
	\emph{calls} arrive to the system at rate $\lambda n$, ie $\lambda$ for each link,
	\emph{capacity} is requested at rate
	\[
		n \cdot \rbb{
			\lambda \cdot \rbb{ 1 - f(x) }
		+	2 \lambda \cdot f(x) \cdot \pr{\text{accept reroute}}
		}
	\Quad{and}
		\pr{\text{accept reroute}} \approx \rbb{ 1 - f(x) }^2,
	\]
	where $f(x) \cq \tfrac1n \abs{ \bra{ j \mid x_j = K } }$ is the proportion of full links.
The capacity requests are uniformly distributed over the non-full links.
Thus the \textit{effective arrival rate}, ie rate at which capacity is requested from a given link, is approximately
\[
	\nu(x)
\cq
	\lambda \rbb{ 1 + 2 f(x) \rbr{ 1 - f(x) } }.
\]
	
We have $\nu(x) \approx \lambda$ when the proportion $f(x)$ of full links is small.
However, there may exist $x \in \Omega$ such that $\nu(x) > K$, even if $\lambda < K$.
Capacity is being requested at a rate faster than it is released in this case.
There is thus drift towards adding calls, like in the description of $\alpha > 1$ above.
If $\alpha$ is sufficiently small, then $\nu(x) < K$ for all $x \in \Omega$.
The parameter $\alpha_c$ is exactly the critical threshold above which we can choose $x \in \Omega$ so that $\nu(x) \ge K$.

\smallskip

The acceptance probability can only increase when there is more than one rerouting attempt, which is the case of \cref{res:intro:retry}.
Exactly the same heuristics thus show that for each $\rho \in \mbn$ there is some critical $\alpha_c(\rho) \in (0,1)$ and further that $\rho \mapsto \alpha_c(\rho) : \mbn \to (0,1)$ is a decreasing map.

\subsubsection{Outline of Proof}

We now have a fairly good idea of the qualitative behaviour of the system in each regime.
We use this to give an outline of the proof.
We use a path coupling argument for the fast mixing regimes.
The coupling between two systems that we use is natural:
	match up the call arrivals and reroutes;
	pair up calls in progress and match their departures where possible;
	let the remaining (`extra') calls depart independently.
This is explained rigorously in \S\ref{sec:fast:var-coup:coup}.

We require a `burn-in' period before attempting to couple two systems. This period is long enough so that, while not necessarily mixed in total variation, the systems have certain typical properties---namely we want the proportion of full links to be roughly correct.
We derive a variant on the \textit{variable length path coupling} technique introduced by \textcite{HV:var-path-coup}. The variant allows for the requirement that both systems exhibit some `typical behaviour' throughout the time interval of interest.
Roughly, we use the departures of the `extra' calls to couple the two systems when $\alpha < \alpha_c$, while we additionally use failed rerouting attempts when $\alpha > 1$.

\smallskip

We use a hitting time approach for the slow mixing regime, ie $\alpha_c < \alpha < 1$.
The motivation for the critical threshold above tells us that there is some $f_* \in (0,1)$ such that
	$\nu(x) < K$ if $f(x) \in [f_* - \eps, f_*]$
and
	$\nu(x) > K$ if $f(x) \in [f_*, f_* + \eps]$,
for some small $\eps > 0$.
Start one system $X$ from $x$ with $f(x) = f_* + \eps$ and another system $Y$ from $y$ with $f(y) = f_* - \eps$.

Rerouting chooses an arbitrary two-link path.
It is thus reasonable to believe that the events that certain links are full are approximately independent.
See \S\ref{sec:intro:prev-mot-comp:prev} for more details on this.
Thus $f(X)$ and $f(Y)$ concentrate around their expectations.
Standard large deviations results then imply that $f(X^t) > f_* > f(Y^t)$ for all $t \le T$ for some large time $T$.
This gives to slow mixing.

\subsubsection{Metastability Here and in Other Models}

For $\alpha_c < \alpha < 1$,
our heuristics imply the following description,
paraphrased from \cite[\S3.7]{KY:stoc-net-book}.

\begin{quote}
	Fix a time period $[0,T]$, arrival rate $\lambda$ and capacity $K$;
	let the number of links $\ninf$.
	The system will freeze in one of two modes:
		either \textit{low-blocking}, where there are relatively few blocked links (and this number decreases as $K$ grows),
		or \textit{high-blocking}, where a proportion bounded away from 0 (independent of $K$) of the links are blocked.
\end{quote}

This \textit{metastability}, or \textit{bistability}, appears throughout mixing literature.
One particularly pertinent example is the Chayes--Machta dynamics in the mean field random cluster model, with percolation parameter $p = \lambda/n \asymp 1/n$.
\textcite{BS:rc-meta} establish the existence of two critical parameters $\lambda_\pm$:
	there is fast mixing in both the low- and high-density regimes, corresponding to $\lambda \in [0, \lambda_-)$ and $\lambda \in (\lambda_+, 1]$ respectively;
	the low- and high-density phases are metastable when $\lambda \in (\lambda_-, \lambda_+)$, giving rise to slow mixing.
See \cite{BS:rc-meta,GLP:rc-meta:jour} for further details.

Another example of metastability is the Ising spin model, eg on a torus $\mbz_n^d$, with $d$ fixed:
	there is fast mixing at high temperature,
	but	slow mixing at low temperature.
This is of a slightly different flavour as metastability is a consequence of symmetry of the underlying spin system.

%
%


\section{Previous Work, Motivation and Comparison}

\subsection{Previous Work}
\label{sec:intro:prev-mot-comp:prev}

Metastability has been known for the original model for a long time.
The precise model described above is not usually used, but rather an approximation to it:
	when a call is using two links, instead of releasing the two links simultaneously after an exponential-1 time, each link is given an independent exponential-1 timer and is released upon the ringing of this timer.

This substantially simplifies the technical details.
	Enumerate the links as $1, ..., n$, so $n = \binom N2$.
	Let $X_j^t$ be the number of calls on link $j \in \bra{1,...,n}$ at time $t \ge 0$.
	Then $X \cq (X_j^t \mid j \in \bra{1,...,n}, \, t \ge 0)$ is a Markov chain.
	Moreover whether three links formed a triangle or not was important in determining their behaviour in the original model;
	the links in the new model are \textit{exchangeable}, ie can be permuted arbitrarily without affecting the equilibrium behaviour.
This \textit{exchangeable model} is the one that has been most studied.
Details below justify why it well-approximates the original.

\medskip

One approach to determine if a system has metastability is to look at a differential equation approximation for an appropriate statistic of the system.
To this end, write $f_k(t)$ for the proportion of links with precisely $k \in [K]$ units of capacity in use at time $t$; so $f_K(t)$ is the proportion blocked.
The following claim is shown in \cite[\S4.3]{K:stoc-net} or \cite[\S3.7]{KY:stoc-net-book}.
For the exchangeable model, $(f_k(\cdot))_{k \in [K]}$ converges weakly to the solution of a multi-dimensional ODE as \ninf.
The proportion of links blocked at a fixed point of this ODE is given by a solution $B$ to the equation
\[
	B = E\rbb{ \alpha\rbr{ 1 + 2B(1-B) }, K }
\Quad{where}
	E(\beta,K) \cq \rbb{ (\beta K)^K / K!} \big/ \rbb{ \sumt[\infty]{k=0} (\beta K)^k / K! }.
\label{eq:history:erlang-fp}
\tag{$*$}
\]
This has a natural interpretation:
	it is the equilibrium probability that an Erlang link with arrival rate $\beta K$ and capacity $K$,
	by which we mean the Markov chain on $\bra{0,...,K}$ with transitions $k \to k+1$ at rate $\beta K$ and $k \to k-1$ at rate $k$.
For $\beta \in (0,\infty)$ fixed, independent of $K$,
we have
\[
\begin{aligned}
	E(\beta,K)
&
=
	\rbb{ 1 + K/(\beta K) + K(K-1)/(\beta K)^2 + \cdots + K!/(\beta K)^K }^{-1}
\\&
\to
	\max\bra{1 - 1/\beta,0}
\quad
	\asKinf;
\end{aligned}
\]
moreover, $E(\beta,K) \ge \max\bra{1-1/\beta,0}$ for all $\beta$ and $K$.
One can derive approximate solutions to the fixed point equation \eqref{eq:history:erlang-fp} from this.
One sees that there exists an $\alpha_c \in (0,1)$ so that if $\alpha \in (\alpha_c,1)$ then \eqref{eq:history:erlang-fp} has two distinct solutions,
say $0 < B_1 < B_2 < 1$.
This implies~metastability:
\begin{itemize}[noitemsep, topsep = \smallskipamount]
	\item 
	if we start the system from the full state, then the proportion of full links converges to~$B_2$;
	
	\item 
	if we start the system from the empty state, then the proportion of full links converges to~$B_1$.
\end{itemize}


These results were proved first by \textcite{GHK:bistability} for the exchangeable~model.
\textcite{CH:ar-original} then combined the techniques from \cite{GHK:bistability} with those of \textcite{H:triangles} to prove the same result for the original model, ie with the graph structure.
Moreover, \citeauthor{CH:ar-original} showed that the limiting ODE is the same for the original model as for the exchangeable model.
This is significant justification that the exchangeable model well-approximates the original.

The exchangeable model can be thought of as adding the assumption ``rerouted pairs release capacity independently'' to the original. Our model adds ``rerouted pairs take capacity independently''.
It is not difficult to write down the corresponding ODE for our model and to see that it has the same limit, as \ninf, as both other models.
In particular, the fixed points are the same.

\smallskip

\textcite{M:indep-block} showed that starting from the assumption that links block independently, one can derive the same ODE.
Informally then, any `sufficiently diverse' rerouting scheme should give rise to an appropriate independence structure and a related ODE.

\smallskip

In related work, \textcite{MR:equilibrium-erlang} study the equilibrium states of large networks of Erlang queues.
Roughly, when a call arrives at a full node, they consider one of two options:
	either additional processing time or extra capacity is required;
	in the latter case they use precisely the \emph{dynamic alternative routing} algorithm which we are considering in this paper.
They study various properties, including stability of the underloaded regime, corresponding to $\alpha < 1$.

\medskip

The extensions of \emph{retries} and \emph{trunk reservation} are considered in \cite[\S5]{GHK:bistability}.
Similar ODE convergence and fixed point analysis to that outlined above is given there.


\subsection{Motivation for Our Model and Comparison with Other Models}
\label{sec:motivation-comparison}

Next, we determine the rate at which calls arrive indirectly, ie via rerouting, to a specific route in the exchangeable model.
It only depends on the state of the system via the proportion $f$ of links full; we denote it $\lambda \cdot \tilde r(f)$.
The total arrival rate, ie direct and indirect, to a specific route is then $\lambda \rbr{ 1 + \tilde r(f) }$.
We show that $\tilde r(f) \to r(f) \cq 2 f(1 - f)$ as \ninf, for $f$ independent of $n$.
The rates, in the limit \ninf, are thus the same in our model as in the exchangeable model.

Suppose a proportion $f$ of the links are full.
The rate at which reroutings are attempted is $\lambda f n$.
Suppose that link $k$ is not full.
Simple counting shows that the probability that the randomly chosen pair contains $k$ and is accepted, ie contains $k$ and another non-full link, is $((1-f)n-1)/\binom n2$.~%
So
\[
	\tilde r(f)
=
	fn\rbb{(1-f)n-1}/\binomt n2
=
	2 f(1-f) \cdot \rbb{1 + \Oh{1/n}}.
\]

The arrival rate was $\lambda (1 + 2f(1-f))$ in our model. We interpret this as rate $\lambda$ directly and $2 \lambda f(1-f)$ indirectly.
Hence the rate in the two models are the same as \ninf.

Thus, while our model may appear rather artificial, we can obtain it by simplifying the original rerouting model, with all its complications:
	first let rerouted calls release their resources independently, which gives the previously-studied exchangeable model;
	then let them take their resources independently;
	finally simplify terms by taking the (\ninf)-limit.

\smallskip

We strongly believe that the fundamental behaviour of our approximate model is the same as that of the original.
One key to this belief is the straightforward to derive fact that the same ODE, and hence with the same fixed point as previously, is satisfied by our model.
This approximate model was suggested to use by \textcite{K:dar:private}.
His justification for its validity was~natural~(paraphrased):
\begin{quote} \slshape
	if the first approximation (paired calls depart independently), which is widely studied and used, is legitimate, then the second (paired calls arrive independently) should~be~too.
\end{quote}

\section{Preliminaries}
\label{sec:prelim}

We now explicitly define some notation and terminology, then the critical threshold $\alpha_c$.
Next, we explain a stochastic domination procedure, which will be crucial for our analysis.
Lastly, we state a result on the mixing time of a single Erlang link and describe a discretisation of the model.

\subsection{Notation and Terminology}
\label{sec:prelim:notation-terminology}

Here we collect some terminology and notation, some of which will be repetition from earlier.

\paragraph{Terminology}
\begin{itemize}[itemsep = 0pt, topsep = \smallskipamount, label = {\ensuremath{*}}]
\item 
We say a link \textit{reroutes} if it is full and a call arrives to it, and hence requests a rerouting.

\item 
We call our model a \textit{dynamic alternative routing system}, with parameters $\alpha$, $K$ and $n$, and denote it by $\DAR_n(\alpha,K)$; we also write $\Pi = \Pi_{\alpha,K,n}$ for its invariant distribution.

\item 
A \textit{single Erlang link} with parameters $\alpha$ and $K$, denoted $\Er(\alpha,K)$, is a single link of capacity $K$ to which calls arrive at rate $\alpha K$ and each depart at rate 1.

\item 
A \textit{product Erlang system} with parameters $\alpha$, $K$ and $n$, denoted $\Er(\alpha,K)^n$, is a system of $n$ independent $\Er(\alpha,K)$ links.

\item 
We call $\alpha$ the \textit{traffic intensity}, $K$ the capacity and $\lambda = \alpha K$ the \textit{arrival rate}.
For $f \in [0,1]$, we call $\beta(f) \cq \alpha( 1 + 2 f(1-f) )$ the \textit{effective traffic intensity} for blocking level $f$.

\item 
An Erlang link or system is \textit{subcritical} if its traffic intensity is strictly less than 1 and \textit{supercritical} otherwise, ie at least 1.
Further, it is \textit{very supercritical} if its traffic intensity is strictly greater than $\sqrt2$.

\item 
We abbreviate \textit{uniformly at random} by \textit{uar}.
\end{itemize}

\paragraph{Notation}
\begin{itemize}[itemsep = 0pt, topsep = \smallskipamount, label = {\ensuremath{*}}]
\item 
Write $\Omega \cq \bra{0,...,K}^n$ for the state space.
We use the standard partial order:\par
{\centering
	$e \le e'$ if and only if $e_i \le e'_i$ for all $i = 1, ..., n$.%
\par}

\item 
For $e \in \Omega$, write
\(
	\varphi(e) \cq \tfrac1n \abs{ \bra{ j \in [n] \mid e_j = K } }.
\)
For $e \in \Omega$, call $e$ a \textit{$\varphi(e)$-blocking} state.

\item 
The invariant distribution $\pi$ of a single $\Er(\beta, K)$ link is,
writing $\nu = \beta K$,
given by
\[
	\pi_\ell \cq \frac{\nu^\ell}{\ell!} \rbbb{ \sum_{k=0}^K \frac{\nu^k}{k!} }^{-1},
\quad
	\ell = 0,1,...,K,
\]
and moreover, in the limit $\Kinf$ with $\beta$ fixed, we have
\[
	\pi_K \cq E(\nu,K) \cq \frac{\nu^K}{K!} \rbbb{ \sum_{k=0}^K \frac{\nu^k}{k!} }^{-1}
\to
	\begin{cases}
		1 - 1/\beta	& \text{ if $\beta \ge 1$}, \\
		0			& \text{ if $\beta \le 1$}.
	\end{cases}
\]
The invariant distribution of a product $\Er(\beta, K)^n$ system is $\pi^n$, by independence of its links.

\item 
For $\alpha \in (0,\infty)$ and $f \in [0,1]$ we define
\[
	r(f) \cq 2f(1-f)
\Quad{and}
	\beta(f) \cq \alpha \rbb{ 1 + r(f) };
\]
note that $\beta(0) = \alpha$ and $\beta(f) \le \beta(\tfrac12) = \tfrac32 \alpha$ for all $f \in [0,1]$.
Also write
\[
	p(f) \cq 1 - 1/\beta(f);
\Quad{note that}
	p(0) = 1 - 1/\alpha.
\]
For $\beta(f) \ge 1$, the equilibrium probability that $\Er(\beta(f),K)$ is full is $p(f)$ in the limit~\Kinf.

\item 
For any process $Z = (Z^t)_{t\ge0}$ taking values in $\Omega$ indexed by time,
write $\varphi^Z_t \cq \varphi(Z^t)$ for $t \ge 0$.
\end{itemize}

\subsection{Definition of Critical $\alpha_c$}
\label{sec:prelim:def-alphac}

We now define the critical threshold $\alpha_c$.
The diversity of the routing in the $\DAR_n(\alpha, K)$ suggests that links should block approximately independently.
Some simple algebra shows that
\[
	p(f) > f
\Quad{if and only if}
	h(f) \cq 2 f^3 - 4 f^2 + f + 1 - 1/\alpha > 0.
\]
We thus define
\[
	\alpha_c
\cq
	\sup
	&\brb{ \alpha \in (0,\infty) \midb p(f) < f \: \forall \, f \in [0,1] }
\\
=
	\hphantom{\sup}
	\mathllap{\inf}
	&\brb{ \alpha \in (0,\infty) \midb \exists \, f \in [0,1] \ST h(f) > 0 }.
\]
Direct calculation, noting that only the constant term in the polynomial $h$ depends on $\alpha$, shows~that
\[
	\alpha_c
=
	\tfrac13 \rbb{ 5 \sqrt{10} - 13 }
\approx
	0.937129.
\]

\subsection{Stochastic Domination}
\label{sec:prelim:stoc-dom}

We describe how to stochastically dominate our system in given sets. This will be key to our analysis.
We say that one system, $Y$, stochastically dominates another, $X$, from above, and write $X \lesssim Y$, if there exists a coupling of the two systems so that when we start $X^0 = Y^0$ the systems satisfy $X^t \le Y^t$ for all times $t$.
Recall that for vectors $x,y \in \Omega = \bra{0,...,K}^n$, we say $x \le y$ if and only if $x_i \le y_i$ for all $i = 1,...,n$; it is a \textit{partial} ordering on $\Omega$.
We say that $Y$ stochastically dominates $X$ from below, and write $X \gtrsim Y$, if the reverse inequality holds.

For $e \in \Omega = \bra{0, ..., K}^n$,
recall that $\varphi(e) = \tfrac1n \abs{ \bra{ j \in [n] \mid e_j = K } }$.

\begin{lem}
\label{res:prelim:stoch-dom:gen}
	Let $X \sim \DAR_n(\alpha,K)$.
	Write
	\(
		\tau_A
	\cq
		\inf\bra{ t \ge 0 \mid X^t \notin A }
	\)
	for the \emph{exit time} of a set $A \subseteq \Omega$.
	There exists a coupling with the following properties.
	\begin{itemize}[itemsep = 0pt, topsep = \smallskipamount, label = \bcdot]
		\item 
		Let
			$f \in [0,1]$,
			$A \cq \bra{ e \mid \varphi(e) \in [f, 1 - f] }$
		and
			$Y \sim \Er(\beta(f), K)^n$.
		Then $X^s \ge Y^s$~for~all~$s \le \tau_A$.
		
		\item 
		Let
			$f \in [0,\tfrac12]$,
			$A \cq \bra{ e \mid \varphi(e) \le f }$
		and
			$Y \sim \Er(\beta(f), K)^n$.
		Then $X^s \le Y^s$ for all $s \le \tau_A$.
		
		\item 
		Let
			$f \in [\tfrac12,1]$,
			$A \cq \bra{ e \mid \varphi(e) \ge f }$
		and
			$Y \sim \Er(\beta(f), K)^n$.
		Then $X^s \ge Y^s$ for all $s \le \tau_A$.
	\end{itemize}
\end{lem}

\begin{rmkt*}
	We use the terminology
		``$X$ stochastically dominates $Y$ from above (or below) while in $A$''
	and write
		``$X \gtrsim Y$ (or $X \lesssim Y$) while $X \in A$''
	to refer to the events of \cref{res:prelim:stoch-dom:gen}.
\end{rmkt*}

\begin{Proof}
We describe the coupling explicitly.
The three properties then follow immediately.
First, we couple departures.
This is simple: we pair up calls where possible link-by-link so that they depart together; the `extra' calls on each link depart independently.

Next, we couple arrivals.
Calls arrive at different links independently in both systems.
Thus it suffices to couple arrivals a single link in each, say $X_1$ and $Y_1$.
Suppose $X$ is in state $x$; the state of $Y$ is irrelevant.
Calls arrive at rate $\beta_x \cq \beta(\varphi(x))$ to $X_1$ and at rate $\beta$ to $Y_1$.
To couple, let a call arrive at rate $\max\bra{\beta_x, \beta}$; it is sent to both $X_1$ and $X_2$ for independent approval. Upon arrival, the call is accepted to $X$ with probability $\min\bra{\beta_x/\beta, 1}\one{X_1 \ne K}$ and to $Y$ with probability $\min\bra{\beta/\beta_x, 1}\one{Y_1 \ne K}$. This indicator forces the call to be declined if the link is full.

The map $f \mapsto \beta(f)$ is increasing on $[0,\tfrac12]$, decreasing on $[\tfrac12, 1]$ and satisfies $\beta(f) = \beta(1-f)$; hence the need for $\varphi(e) \in [f,1-f]$, $f \in [0,\tfrac12]$ and $f \in [\tfrac12,1]$ in the three cases, respectively.
\end{Proof}

\begin{cor}
\label{res:prelim:stoch-dom:cor}
	Let $X \sim \DAR_n(\alpha,K)$.
	Then $\Er(\alpha,K)^n \lesssim X \lesssim \Er(\tfrac32\alpha,K)^n$.
\end{cor}

\begin{Proof}
	This follows immediately from \cref{res:prelim:stoch-dom:gen} and $\alpha \le \beta(f) \le \beta(\tfrac12) = \tfrac32 \alpha$ for all $f \in [0,1]$.
\end{Proof}

Finally, since the product Erlang systems are independent queues, there is monotonicity.
Recall that two distributions $\Theta$ and $\tilde \Theta$ satisfy $\Theta \lesssim \tilde \Theta$ if they can be coupled so that two realisations, $\theta$ and $\tilde \theta$, respectively, satisfy $\theta \le \tilde \theta$.

\begin{lem}
\label{res:prelim:stoch-dom:mon}
	Let $\beta, \tilde \beta \in (0,\infty)$.
	Let $\Theta$ and $\tilde \Theta$ be two distributions on $\Omega$.
	The following hold.
	\begin{itemize}[itemsep = 0pt, topsep = \smallskipamount]
		\item 
		Suppose that $\beta \le \tilde \beta$.
		Let $E \sim \Er(\beta, K)^n$ and $\tilde E \sim \Er(\tilde \beta, K)^n$ with $E^0 \le \tilde E^0$.
		Then $E \lesssim \tilde E$.
		
		\item 
		Suppose that $\Theta \lesssim \tilde \Theta$.
		Let $E, \tilde E \sim \Er(\beta, K)^n$ with $E^0 \sim \Theta$ and $\tilde E^0 \sim \tilde \Theta$.
		Then $E \lesssim \tilde E$.
	\end{itemize}
\end{lem}

\begin{Proof}
	These claims are immediate consequences of properties of standard Poisson processes.
\end{Proof}

%

Observe that if we did not make our final approximation, namely uncoupling the release of two-link calls, then we would not be able to do these dominations so easily. In particular, two calls may arrive at the same time, so we can never dominate from above by an Erlang system in which calls only arriving one at a time.
This is the only place the approximation is really needed.


\subsection{Mixing Time for Erlang Systems}
\label{sec:prelim:er-mix}

We use, repeatedly, the mixing time for Erlang systems in the fast mixing proofs.

\begin{lem}
\label{res:prelim:product-mixing}
	For
		all $\beta \in (0,\infty)$,
		all $K \ge 4$
	and
		all $n$ sufficiently large,
	writing $\tmix^{\beta,K;n}(\cdot)$ for the mixing time of an $\Er(\beta,K)^n$ system,
	we have
	\[
		\tmix^{\beta,K;n}(1/n) \le 10 \log n.
	\]
\end{lem}

\begin{Proof}
See \cref{res:er-mix:multi-link:mixing} in \cref{app:er-mix}.
Take $\alpha \cq \beta$ and $\sigma \cq 0$ in the notation there.
\end{Proof}

\subsection{Discretisation and Set of States Visited}
\label{sec:prelim:discrete}

We next describe a discretisation of our system.
This allows us to control the set of states visited by the system in a given time.

\begin{defn}
\label{def:prelim:discrete}
Define a discrete-time process $Z \cq (Z^m)_{m \in \mbn_0}$ by the following step distribution.
\begin{quote}
	Draw $B \sim \Bern\rbr{ 1/(2 \alpha + 1) }$.
	
	\begin{itemize}
		\item 
		Suppose $B = 1$.
		
		Select a slot uniformly at random amongst all $K n$ and set this slot to be empty.
		
		\item 
		Suppose $B = 0$.
		Write $f$ for the current proportion of full links.
		
		Sample $B' \sim \Bern\rbb{ \tfrac12 \rbr{ 1 + 2 f(1 - f) } }$.
			If $B' = 1$, then
			choose a link uniformly at random and add a call to this link if it is not already full.
			Do nothing if $B' = 0$.
	\end{itemize}
\end{quote}
Let $S \cq (S_m)_{m \in \mbn}$ be the jump times of a rate-$(2 \alpha + 1) K n$ Poisson process with $0 \cq S_0 < S_1 < S_2 < \cdots$.
Define the continuous-time process $X \cq (X^t)_{t\ge0}$ by $X^t \cq Z^m$ for $t \in [S_m, S_{m+1})$.
\end{defn}

The following lemma is straightforward to prove.

\begin{lem}
\label{res:prelim:discrete}
	We have $X \sim \DAR_n(\alpha, K)$.
	Let $T \in [0, \infty)$ and $M \in \mbn_0$.
	Then
	\begin{gather*}
		\brb{ X^t \midb t \in [0,T] }
	\subseteq
		\brb{ Z^m \midb m \in [0,M] \cap \mbz }
	\Quad{on the event}
		\bra{ S_M \ge T }.
	\end{gather*}
\end{lem}

\section{Slow Mixing in Interim Regime: $\alpha_c < \alpha < 1$}
\label{sec:slow}

This section is devoted to the establishment of slow mixing in the interim regime $\alpha \in (\alpha_c,1)$.
Our aim is to show slow mixing via two analogous results:
	first,
		if the system starts from a `generic' high-blocking state,
		then it remains in such a state for exponentially long;
	second,
		if it starts in a low-blocking state, eg the empty configuration,
		then it remains in such a state for exponentially long.
The first result will hold due to $\alpha > \alpha_c$, and the second due to $\alpha < 1$.

If $\alpha > 1$, then it is clear that the system will `fill up' quickly, and so it would not be the case that the system remains in a low-blocking state for exponentially long. Of course, this does not imply any result on mixing times, but merely says that \emph{this particular method} will not be helpful.

The precise statement that we prove in this section is the following.

\begin{thm}
\label{res:slow}
	For all $\alpha \in (\alpha_c,1)$,
	there exists a constant $c$ so that,
	for all $K$ and $n$ sufficiently large,
	for all $t \le e^{cn}$,
	we have
	\[
		\MAX{x \in \Omega} \,
		\tvb{ \prb[x]{ X^t \in \cdot } - \Pi }
	\ge
		\tfrac12 - e^{-cn}.
	\]
	Thus
	for all $K$ sufficiently large,
	all $\eps \in (0,\tfrac12)$
	and all $n$ sufficiently large,
	we have
	$\tmix(\eps) \ge e^{cn}$.
\end{thm}


We now state the first of the two results that we wish to prove.
It is the formalisation of ``moving from high-blocking to low-blocking takes exponentially long''.

\begin{prop} \label{res:slow:high-to-low}
	For all $\alpha > \alpha_c$,
	there exist constants $c > 0$ and $\delta \in (0,1)$ and an $x \in \Omega$ so that,
	for all $K$ and $n$ sufficiently large,
	\[
		\tau
	\cq
		\inf\brb{ t \ge 0 \midb \tfrac1n \abs{\bra{j \mid X_j^t = K }} \le \delta }
	\Quad{satisfies}
		\prb[x]{ \tau \le e^{c n} }
	\le
		e^{-c n};
	\]
	in words, at least a proportion $\delta$ of links are full for time $e^{cn}$ with probability at least $1 - e^{-cn}$.
\end{prop}

\begin{Proof}
Recall the following notation:
for $e \in \bra{0,...,K}^n$ and $f \in [0,1]$,
write
\[
	\varphi(e) \cq \tfrac1n \absb{ \brb{ j \midb e_j = K } },
\quad
	\beta(f) \cq \alpha\rbb{1 + 2f(1-f)}
\Quad{and}
	p(f) \cq 1 - 1/\beta(f);
\]
in words,
	$\beta(f)$ is the effective traffic intensity when the blocking proportion is $f$
and
	$p(f)$ is the equilibrium probability that a supercritical $\Er(\beta(f),K)$ link is full in the limit $\Kinf$.

Assume that $\alpha \le 1$.
Since $\alpha > \alpha_c$, by definition there exists $\delta \in [0,1]$ so that $p(\delta) > \delta$.
We may assume that $\delta \in [0, \tfrac12]$ since $p(f) = p(1-f)$ for all $f$.
We may further assume that $\delta \in (0, \tfrac12)$ since $p(0) = 0 < \alpha$ and $p(\tfrac12) = 1 - (\tfrac32 \alpha)^{-1} < \tfrac12$.
Choose $\eta > 0$ sufficiently small so that $p(\delta) > \delta + 2 \eta$.

By \cref{res:prelim:stoch-dom:cor}, we have $X \lesssim U \sim \Er(\tfrac32\alpha, K)^n$.
By \cref{res:prelim:stoch-dom:gen}, we have
\[
	X \gtrsim L \sim \Er(\beta(\delta), K)^n
\QUAD{while}
	\varphi(X) \in [\delta,1-\delta].
\]
Write $\Pi_L$ and $\Pi_U$ for the invariant distributions of $L$ and of $U$, respectively.
Also write
\[
	\mfu \cq \brb{ e \in \Omega \midb \varphi(e) \le \tfrac12 },
\quad
	\mfl \cq \brb{ e \in \Omega \midb \varphi(e) \ge \delta }
\QUAD{and}
	\mfs \cq \mfu \cap \mfl.
\]
Recall that $\delta < \tfrac12$.
Write $\tau_\mfu$, $\tau_\mfl$ and $\tau_\mfs$ for the exit times of $\mfu$, of $\mfl$ and of $\mfs$, respectively, by $X$.
Clearly, $\tau \ge \tau_\mfs = \tau_\mfu \wedge \tau_\mfl$.
By \cref{res:prelim:stoch-dom:gen},
we can couple so that
\[
	L^t \le X^t \le U^t
\Quad{for all}
	t \in [0, \tau_\mfs].
\]
In particular if $X$ has exited $\tau_\mfs$, then either $L$ has exited $\mfl$ or $U$ has exited $\mfu$.
Write $\sigma_\mfl$ and $\sigma_\mfu$ for these two times, respectively.
Hence
\(
	\tau_\mfs
\ge
	\sigma_\mfu \wedge \sigma_\mfl.
\)
Thus, for any $x \in \Omega$ and $t \ge 0$, we have
\[
	\prb[x]{ \tau \le t }
\le
	\prb[x]{ \tau_\mfs \le t }
\le
	\prb[x]{ \sigma_\mfu \wedge \sigma_\mfl \le t }
\le
	\prb[x]{ \sigma_\mfu \le t }
+	\prb[x]{ \sigma_\mfl \le t }.
\label{eq:slow:tau-sigma}
\nt
\]

Write $\Pi_L$ and $\Pi_U$ for the invariant distribution of $L$ and of $U$, respectively.
We have
\(
	\Pi_L \lesssim \Pi_U
\)
since $L \lesssim U$ by \cref{res:prelim:stoch-dom:gen}.
Note that $\mfu$ is an `up-set': if $e \in \mfu$ and $e' \le e$, then $e' \in \mfu$ too.
By \cref{res:prelim:stoch-dom:mon},
we then have
\[
	\prb[\Pi_L]{ \sigma_\mfu \le t }
\le
	\prb[\Pi_U]{ \sigma_\mfu \le t }.
\]
Combined with \eqref{eq:slow:tau-sigma},
we deduce that
\[
	\prb[\Pi_L]{ \tau \le t }
\le
	\prb[\Pi_L]{ \sigma_\mfu \le t }
+	\prb[\Pi_L]{ \sigma_\mfl \le t }
\le
	\prb[\Pi_U]{ \sigma_\mfu \le t }
+	\prb[\Pi_L]{ \sigma_\mfl \le t }.
\label{eq:slow:exit-time-dom}
\nt
\]

\smallskip

We first analyse $\sigma_\mfl$, ie how long it takes $L$ to leave $\mfl$, when started from its invariant distribution $\Pi_L$.
The independence of the coordinates and the definition of $p$ implies that
\[
	\exb[\Pi_L]{ \varphi(L^0) }
\to
	p(\delta)
\ge
	\delta + 2\eta
\quad
	\asKinf.
\]
The the product nature of the system thus implies that
for sufficiently large $K$,
we have
\[
	n \varphi(\mcl) \gtrsim \Bin(n, \delta + \eta)
\Quad{when}
	\mcl \sim \Pi_L.
\]
Concentration of the Binomial, eg Hoeffding's inequality, then gives
\[
	\Pi_L\rbb{\mfl^c}
=
	\prb[\Pi_L]{ \varphi(L^0) < \delta }
\le
	\prb[\Pi_L]{ \Bin(n,\delta+\eta) < \delta n }
\le
	\expb{ - 2 \eta^2 n }.
\label{eq:slow:inv-exp-decay}
\nt
\]

The above holds with $\mcl$ replaced by $L^t$ when $L^0 \sim \Pi_L$ since $\Pi_L$ is invariant for $L$.
Consider the discrete-time process $Z$ coupled to $X \sim \DAR_n(\alpha, K)$ as in \cref{def:prelim:discrete}:
	$X^t = Z^m$ for $t$ with $S_m \le t < S_{m+1}$.
If we run $X$ for a time $e^{a' n}$, for some constant $a' > 0$, then $Z$ takes at most $e^{2 a' n}$ steps with probability at least $1 - e^{- b' n}$, for some constant $b' > 0$, by concentration of the Poisson distribution.
Applying \cref{res:prelim:discrete} and \eqref{eq:slow:inv-exp-decay} along with a union bound over the steps of the discrete-time chain $Z$ shows that
there exist constants $a_1, b_1 > 0$ so that
\[
	\prb[\Pi_L]{ \sigma_\mfl \le e^{a_1 n} }
\le
	e^{-b_1 n}.
\label{eq:slow:lower}
\nt
\]

\smallskip

The analysis of $\sigma_\mfu$ is similar.
We can take $\eta \cq \tfrac1{20}$ here, for example, since $\max_{f \in [0,1]} p(f) = p(\tfrac12) \le \tfrac13 < \tfrac12$.
We obtain constants $a_2, b_2 > 0$ so that
\[
	\prb[\Pi_U]{ \sigma_\mfu \le e^{a_2 n} }
\le
	e^{-b_2 n}.
\label{eq:slow:upper}
\nt
\]

\smallskip

Now set $a \cq \min\bra{a_1,a_2}$ and $b \cq \min\bra{b_1,b_2}$.
Plugging (\ref{eq:slow:lower}, \ref{eq:slow:upper}) into \eqref{eq:slow:exit-time-dom} gives
\[
	\prb[\Pi_L]{ \tau \le e^{a n} }
\le
	\prb[\Pi_L]{ \sigma_\mfu \wedge \sigma_\mfl \le e^{a n} }
\le
	2 e^{-b n}.
\]
Since $\mbp_{\Pi_L}$ is an averaging measure, this implies that there exists some $x \in \Omega$ so that
\[
	\prb[x]{ \tau \le e^{a n} }
\le
	2 e^{-b n}.
\]

\smallskip

If $\alpha > 1$, then simply stochastically lower bound $X \gtrsim \Er(\alpha,K)^n$, using \cref{res:prelim:stoch-dom:cor}. The above argument, but requiring only $\sigma_\mfl$, applies, as here $p(0) = 1 - 1/\beta(0) = 1 - 1/\alpha > 0$.
	%
\end{Proof}

We have just shown that if we start in high-blocking, then it takes exponentially long to get to low-blocking. We now show that the converse is also true.
Write $0 \cq (0,...,0) \in \bra{0,...,K}^n$ for the state where every link is empty.

\begin{prop} \label{res:slow:low-to-high}
	For all $\alpha < 1$ and all $f \in (0,1)$,
	there exists a constant $c > 0$ so that,
	for all $K$ and $n$ sufficiently large,
	\[
		\tau \cq \brb{ t \ge 0 \midb \tfrac1n \abs{\bra{ j \mid X_j^t = K }} \ge f }
	\Quad{satisfies}
		\prb[0]{ \tau \le e^{cn} } \le e^{-cn};
	\]
	in words, at most a proportion $f$ are full for time $e^{cn}$ with probability at least $1 - e^{-cn}$.
\end{prop}

\begin{Proof}
Recall the following notation:
for $e \in \bra{0,...,K}^n$ and $f \in [0,1]$, write
\[
	\varphi(e) \cq \tfrac1n \absb{ \brb{ j \midb e_j = K } }
\Quad{and}
	\beta(f) \cq \alpha\rbb{1 + 2f(1-f)};
\]
in words, $\beta(f)$ is the effective traffic intensity, when the blocking proportion is $f$.

Since $\alpha < 1$ there exists an $\delta \in (0,\tfrac12)$ so that $\beta(\delta) < 1$.
By monotonicity (in $f$) of the condition
\(
	\varphi(X^t) \ge f
\)
in the definition of $\tau$,
we may assume that $f < \delta$.

By \cref{res:prelim:stoch-dom:gen}, we have
\[
	X \lesssim U \sim \Er(\beta(\delta),K)^n
\Quad{while}
	\varphi(X^t) \le \delta,
\]
and so in particular up until time $\tau$.
Further, by \cref{res:prelim:stoch-dom:mon}, we may assume that $U^0 \sim \Pi_U$, the invariant distribution of $U$.
Hence it suffices to prove the statement with $\tau$ replaced by
\[
	\sigma_\mfu \cq \inf\brb{t \ge 0 \midb \varphi(U^t) \ge f}.
\]

For a single $\Er(\beta(\delta),K)$ link, the probability of being full tends to 0 as $\Kinf$, since $\beta(\delta) < 1$, and so certainly becomes at most $\tfrac12 \delta$ if $K$ is sufficiently large.
The proof is completed analogously to the previous one, using concentration of the proportion of full links in a product Erlang system.
\end{Proof}

From these hitting time results we obtain our mixing lower bound..

\begin{Proof}[Proof of \cref{res:slow}]
Fix $\alpha \in (\alpha_c,1)$.
Choose $\delta$, $c_1$ and $x$ as guaranteed by \cref{res:slow:high-to-low}.
Then choose $c_2$ as guaranteed by \cref{res:slow:low-to-high} with $f \cq \delta$.
Set $c = \min\bra{c_1,c_2}$.

Let $X,Y \sim \DAR_n(\alpha,K)$ with $X^0 = x$ and $Y^0 = y$.
Let
\[
	\tau_X \cq \inf\brb{ t \ge 0 \midb \varphi(X^t) \le f }
\Quad{and}
	\tau_Y \cq \inf\brb{ t \ge 0 \midb \varphi(Y^t) \ge f },
\]
recalling that $\varphi(e) \cq \varphie$ for $e \in \bra{0,...,K}^n$.
Since $\varphi(x) > f$ and $\varphi(0) = 0 < f$, if $X^t = Y^t$ then $\tau_X \wedge \tau_Y \le t$.
However, \cref{res:slow:high-to-low,res:slow:low-to-high} tell us, respectively, that
\[
	\prb[x]{ \tau_X \le e^{c n} }
\le
	e^{-c n}
\Quad{and}
	\prb[0]{ \tau_Y \le e^{c n} }
\le
	e^{-c n}.
\]
We thus deduce that
\[
	\tvb{ \pr[x]{ X^t \in \cdot } - \pr[0]{ Y^t \in \cdot } }
\ge
	1 - 2 e^{-c n}
\Quad{when}
	t \le e^{c n},
\]
using the union bound.
The claim now follows by the triangle inequality.
\end{Proof}


\counterwithin{thm}{subsection}
\section{Fast Mixing in Edge Regimes: $\alpha < \alpha_c$ or $\alpha > 1$}

We establish fast mixing for the regimes $\alpha < \alpha_c$ and $\alpha > 1$ in this section.
As mentioned before, we use a variable length path coupling argument, introduced by \textcite{HV:var-path-coup}. We use a minor variant of their result, given in \cref{res:var-length-coup} below.

The high-level ideas for the two regimes will be the same:
	we use the same coupling,
	the same variable length path coupling type argument
and
	a burn-in phase with a similar flavour.
The stopping times used in the variable length path coupling will be different:
	when $\alpha < \alpha_c$, there are very few reroutings and we use the ending of calls to couple;
	when $\alpha > 1$, there will be a significant number of reroutings and we use failed attempts at rerouting to couple.
The flavour of the burn-in phase will be very similar in the two regimes:
	when $\alpha < \alpha_c$, we run until we are in a low-blocking state;
	when $\alpha > 1$, we run until we are in a high-blocking state.

From a holistic point of view, one should really think of this as ``the fast mixing case with two subcases'', rather than ``two fast mixing cases''.

\subsection{Variable Length Coupling Set-Up}
\label{sec:fast:var-coup}

Let $X,Y \sim \DAR_n(\alpha,K)$ be two dynamic alternative routing systems.
Let $\mfs \subseteq \Omega$, and write $\msg[0,t]$ for the event that $X$ and $Y$ are in $\mfs$ for the entire interval $[0,t]$, ie
\[
	\msg[0,t]
\cq
	\brb{ (X^s,Y^s) \in \mfs^2 \: \forall \, s \in [0,t] }.
\]
The definition of $\mfs$ will be different for the two regimes $\alpha < \alpha_c$ and $\alpha > 1$ (albeit of the same flavour).
The burn-in phase will run for sufficiently long so that both $X$ and $Y$ are `far enough inside' $\mfs$ so that they remain in $\mfs$ for a long time; this is, of course, made precise later.

We first define the coupling, then describe how to bound the coupling time using the variable length path coupling technique, given that $X$ and $Y$ remain in $\mfs$ throughout.

\subsubsection{Coupling}
\label{sec:fast:var-coup:coup}

The fact that we have obtained our model as a simplification, or approximation, of a rerouting scheme allows us to consider it in a more instructive (and intuitive) way than simply ``a variable rate Poisson arrival system with independent exponential departures'':
	we can set it up as a \emph{type} of rerouting scheme.
Recall that before (in the original model) two calls were added upon a rerouting, whereas we (in our model) only add one at a time.
To get the correct rates, we assume that the reroutings happen twice as fast as in the previous model.
Since the reroutings are accepted with probability approximately $(1-f)^2$, where $f$ is the current blocking level, it does not make sense to say they are accepted with probability $2(1-f)^2$, since this number may be larger than 1. Instead, we double the \emph{entire} arrival rate of the system and still reroute with probability $(1-f)^2$, but now choose a \emph{single} link (rather than a \emph{pair}) to have a call added; additionally, when a call arrives to a non-full link, it is only accepted with probability $\tfrac12$.
Equivalently, we could say that to non-full links there is a Poisson stream of rate $\lambda$ and to full links there is a Poisson stream of rate $2 \lambda$.

After that motivation, we can now give an explicit way to realise the system.

\begin{defn}
\label{def:realisation}
For arrivals, give to each link a Poisson stream (of arriving calls) of rate $2 \lambda$.
Upon a call's arrival to a link, $k$ say, we have the following procedure.
\begin{itemize}
	\item 
	If the link $k$ is not full, then toss a $\Bern(\tfrac12)$-coin:
	\begin{itemize}[noitemsep, topsep=0pt]
		\item if heads (ie `1'), then add a call to link $k$;
		\item if tails (ie `0'), then do nothing.
	\end{itemize}
	
	\item 
	If the link $k$ is full, then choose two links $i$ and $j$ uar (with replacement):
	\begin{itemize}[noitemsep, topsep=0pt]
		\item if both links $i$ and $j$ are not full, then add a call to link $i$;
		\item otherwise, ie if either link $i$ or link $j$ is full, do nothing.
	\end{itemize}
\end{itemize}

For departures, give to each call in the system an independent exponential-$1$ timer. Upon a timer's ringing, remove the corresponding call from the system.
\end{defn}

Since the probability that both $i$ and $j$ are not full is precisely $(1-f)^2$ when the system is in an $f$-blocking state, we see that this is a genuine realisation of the system.

%

\smallskip

It is this realisation, which is similar to the original dynamic alternative routing system, that we have in mind for the remainder of the paper. We speak of \textit{reroutings} with this interpretation.

Moreover, this realisation of the system lends itself very naturally to a coupling of two (or even more) systems---it will also extend (relatively) easily when we consider `multiple attempts at rerouting'
in \S\ref{sec:retries}.
Informally, we just use the same fair coin (to accept/reject calls which arrive to non-full links) and selection of $(i,j)$ (for reroutings) in each system.
Recall that $\lambda = \alpha K$.

\begin{defn}
\label{def:coupling}
For arrivals, give to each link a Poisson stream (of arriving calls) of rate $2 \lambda$.
Upon a call's arrival to a link, $k$ say, we have the following procedure.
\begin{itemize}
	\item
	Suppose $k$ is not full in either of $X$ or $Y$.
	Toss a $\Bern(\tfrac12)$-coin:\par
	\quad if heads, then add a call to link $k$ both in $X$ and in $Y$.
	
	\item
	Suppose $k$ is full in $X$ but not full in $Y$.
	Toss a $\Bern(\tfrac12)$-coin:\par
	\quad if heads, then add a call to link $k$ in $Y$.\par
	Also, independently, choose two links $i$ and $j$ uar (with replacement):\par
	\quad if both links $i$ and $j$ are not full in $X$, then add a call to link $i$ in $X$.

	\item
	Suppose $k$ is full in $Y$ but not full in $X$.
	Do analogously to the previous case.
	
	\item
	Suppose $k$ is full both in $X$ and in $Y$.
	Choose two links $i$ and $j$ uar (with replacement):\par
	\quad if both $i$ and $j$ are not full in $X$ (respectively $Y$), then add a call to $i$ in $X$ (respectively~$Y$).
\end{itemize}

For departures, use the same rate-1 departure clocks in $X$ as in $Y$ where possible, giving the `extra' calls (ie those in $X$ but not in $Y$ or vice versa) independent rate-1 departure clocks.
\end{defn}

\begin{rmkt*}
	By inspection, one can see that this is a genuine, Markovian coupling.
	When using this coupling and $(X^0,Y^0) = (x,y)$, we denote it $\mbp_{x,y}$.
	Write $\pr{\cdot} \cq \max_{(x,y) \in \Omega^2} \pr[x,y]{\cdot}$.
	Furthermore, it is a \textit{coalescent} coupling:
		writing $\tau_c \cq \inf\bra{ t \ge 0 \mid X^t = Y^t }$,
		we have $\bra{ X^t \ne Y^t } = \bra{ \tau_c > t }$.
\end{rmkt*}

\subsubsection{Variable Length Bound via Stopping Time}
\label{sec:fast:var-coup:bound}

The following is an adaptation of the \textit{variable length path coupling} of \textcite{HV:var-path-coup}.
It holds for any coalescent coupling; we always use the one from \cref{def:coupling}.
(We change the notation slightly, compared with \cite{HV:var-path-coup}, so as to not clash with our already-established notation.)

\begin{thm}[{cf \cite[Corollary 6]{HV:var-path-coup}}]
\label{res:var-length-coup}
Let $(X,Y)$ be a coalescent coupling of two realisations of the same Markov chain, with state space $\Omega$.
Let $\tau$ be a stopping time for the joint chain $(X,Y)$.
Let $\mfs \subseteq \Omega$ and, for $t \ge 0$, write
\[
	\msg[0,t] \cq \brb{ (X^s,Y^s) \in \mfs^2 \: \forall \, s \in [0,t] }.
\]
Let $d$ denote the graph distance on the graph on $\Omega$ induced by the permissible transitions of the Markov chain.
Write
\(
	S \cq \bra{ (x,y) \in \Omega^2 \mid d(x,y) = 1 }
\)
for the pairs of neighbours;
define
\[
	\gamma_0
\cq
	\MAX{(x,y) \in S}
	\exb[x,y]{ d(X^\tau,Y^\tau) \one{\msg[0,\tau]} }.
\]
Write $W$ for the supremum of $\max_{s \ge 0} d(X^s,Y^s) \one{s < \tau}$ over all possible evolutions $(X^s,Y^s)_{s\ge0}$ which have starting pair $(X^0,Y^0) \in S$.
Suppose that $M > 0$ satisfies
\[
	\MAX{(x,y) \in S}
	\prb[x,y]{ \tau > M }
\le
	\tfrac12 (1 - \gamma_0) / W.
\]
Then,
writing $\gamma \cq \tfrac12(1+\gamma_0)$,
for any $t \ge 0$,
we have
\[
	\MAX{x,y \in \Omega}
	\prb[x,y]{ X^t \neq Y^t, \, \msg[0,t] }
\le
	\diam \Omega \cdot \gamma^{t/M-1}.
\]
\end{thm}


\begin{Proof}
The coupling is coalescent, so if $X^t \neq Y^t$ then $X^s \neq Y^s$ for all $s \le t$.
Fix an $M \in \mbr_+$ and set $\sigma \cq \tau \wedge M$, which is deterministically bounded.
We split the interval
\(
	(0,kM]
\QUAD{into}
	(0,M] \cups ((k-1)M,M],
\)
and, for $\ell \in [k]$ consider the probability of coalescence by $\ell M$:
for $(x,y) \in \Omega^2$, let
\[
	p_\ell(x,y)
\cq
	\prb[x,y]{ X^{\ell M} \neq Y^{\ell M}, \, \msg[0,\ell M] };
\]
then write $p_\ell^*$ for the the maximum over adjacent pairs $(x,y)$, ie
\(
	p_\ell^*
\cq
	\maxt{(x,y) \in S}
	p_\ell(x,y).
\)
By the union bound along shortest paths,
we have
\(
	p_\ell(x,y) \le p_\ell^* d(x,y)
\)
for all $(x,y) \in \Omega^2$.

For $(x,y) \in S$ and $(u,v) \in \Omega^2$, define
\[
	Q(x,y; u,v)
\cq
	\prb[x,y]{ X^{\sigma} = u, \, Y^{\sigma} = v, \, \msg[0,M] }.
\]
Fix some $\ell \in \bra{1, ..., k}$.
We bound $p_\ell^*$ inductively:
for any $(x,y) \in S$, by the strong Markov property (applied at time $\sigma \le M$) and the fact that the coupling is coalescent,
we have
\[
	p_{\ell}(x,y)
&
\le
	\sumt{(u,v) \in \Omega^2}
	p_{\ell-1}(u,v) Q(x,y; u,v)
\\&
\le
	\sumt{(u,v) \in \Omega^2}
	p_{\ell-1}^* d(u,v) Q(x,y; u,v)
=
	p_{\ell-1}^* \, \exb[x,y]{ d(X^{\sigma},Y^{\sigma}) \one{\msg[0,M]} };
\]
note that the first relation is an inequality, rather than an equality, because coalescence may occur in the final $M - \sigma$ time units.
Hence, maximising over $(x,y) \in S$, we obtain
\[
	p_{\ell}^* \le \zeta p_{\ell-1}^*
\Quad{where}
	\zeta \cq \MAX{(x,y) \in S} \exb[x,y]{ d(X^{\sigma},Y^{\sigma}) \one{\msg[0,M]} }.
\]
Iterating this, we obtain $p_\ell^* \le \zeta^\ell$ as $p_0^* = 1$.
Hence
\[
	\prb{ X^{kM} \neq Y^{kM}, \, \msg[0,kM] }
\le
	\diam \Omega \cdot p_k^*
\le
	\diam \Omega \cdot \zeta^k.
\]

Next define the following quantities,
which we use to control $\zeta$,
recalling that $\sigma = \tau \wedge M$:
\[
	\gamma_0
\cq
	\MAX{(x,y) \in S}
	\exb[x,y]{ d(X^{\tau},Y^\tau) \one{\msg[0,\tau]}}
\Quad{and}
	\gamma_1
\cq
	\MAX{(x,y) \in S}
	\exb[x,y]{ d(X^M,Y^M) \cdot \one{\tau > M} };
\]
Then $\zeta \le \gamma_0 + \gamma_1$.
By definition of $W$,
we have
\(
	\gamma_1
\le
	W \cdot \maxt{(x,y) \in S} \prb[x,y]{\tau > M}.
\)
By the assumption on $M$,
the probability in the above display is at most $\tfrac12(1-\gamma_0)/W$.
Thus
\(
	\zeta
\le
	\gamma_0 + \gamma_1
\le
	\tfrac12(1+\gamma_0)
=
	\gamma.
\)
Noting that $\floor{t/M} \ge t/M - 1$, the final claim follows.
\end{Proof}

For our application, our stopping time $\tau$ for the pair $(X,Y)$ will be defined via a set of stopping rules; in the spirit of path coupling, these rules will require $X$ and $Y$ to initially be adjacent.

\subsection{Slow Arrivals: $\alpha < \alpha_c$ and Low-Blocking}
\label{sec:low}

In this subsection
we consider the `slow arrivals' regime, ie $\alpha < \alpha_c$.
We always use the coupling \mbpcoup from \cref{def:coupling}.
We prove the following theorem.

\begin{thm} \label{res:low}
	For all $\alpha < \alpha_c$,
	there exists a constant $C$ so that,
	for all $K$ and $n$ sufficiently large,
	if $X,Y \sim \DAR_n(\alpha,K)$, then,
	under the coupling \mbpcoup,
	for all $t \ge C \log n$,
	we have
	\[
		\MAX{(x,y) \in \Omega^2}
		\prb[x,y]{ X^t \neq Y^t }
	\le
		C / n
	=
		\oh1.
	\]
	Thus
	for all $K$ sufficiently large,
	all $\eps \in (0,1)$
	and all $n$ sufficiently large,
	we have
	$\tmix(\eps) \le C \log n$.
\end{thm}

Used throughout this subsection
repeatedly, will be the notation $\beta(f)$ and $p(f)$:
\begin{itemize}[itemsep = 0pt, topsep = \smallskipamount, label = \bcdot]
	\item 
	$\beta(f) = \alpha(1 + 2f(1-f))$ is the effective traffic intensity when a proportion $f$ of links are full;
	
	\item 
	$p(f) = 1 - 1/\beta(f)$ is the equilibrium probability that an $\Er(\beta(f),K)$ link is full as $\Kinf$.
\end{itemize}


We prove this theorem via a sequence of lemmas:
	in \S\ref{sec:low:burnin} we describe the burn-in phase;
	in \S\ref{sec:low:var-coupling} we describe and apply the variable length path coupling;
	finally we conclude in \S\ref{sec:low:conclusion}.

\subsubsection{Burn-In Phase}
\label{sec:low:burnin}

Since $p(f) < f$ for all $f$ when $\alpha < \alpha_c$, the system cannot `support' $f$-blocking for any $f$; this concept is elaborated on below.
The aim of this part
is to prove the following burn-in proposition.

\begin{prop}
\label{res:low:burnin}
	For all $\alpha < \alpha_c$,
	there exists a constant $C$ so that,
	for all $\eps \in (0,1)$,
	there exists $\xi_0 \in (0,1)$ so that,
	for all $\xi \in (\xi_0,1)$,
	setting
	$t \cq C \log n$,
	\[
		\mfs
	&\cq
		\brb{ e \in \Omega \midb \tfrac1n \abs{ \bra{ j \mid e_j > \xi K } } \le 2 \eps }
	\quad
		\text{and}
	\\
		\mfb
	&\cq
		\brb{ e \in \Omega \midb \pr[e]{X^s \in \mfs \: \forall \, s \le n} \ge 1 - 1/n },
	\]
	for all $K$ and $n$ sufficiently large (depending on $\alpha$ and $\eps$),
	if $X \sim \DAR_n(\alpha, K)$, then
	\[
		\MAX{x \in \Omega} \,
		\prb[x]{ X^t \notin \mfb }
	\le
		\tfrac15 C / n
	=
		\oh1.
	\]
\end{prop}

\begin{Proof}[Intuition]
\renewcommand{\qedsymbol}{\ensuremath{\triangle}}
By definition of $\alpha_c$, when $\alpha < \alpha_c$, for any $f \in (0,1)$, we have $p(f) < f$.
Thus the expected proportion of full links in an $\Er(\beta(f),K)^n$ system in equilibrum is \emph{less than} $f$ (in the limit $\Kinf$).
Intuitively, this suggests that for any $f$ the system cannot `support' $f$-blocking:
	if we start a $\DAR_n(\alpha, K)$ system from a state with proportion $f$ blocked then (typically) the proportion will initially decrease.
Thus if the proportion blocked initially, call it $f_0$, is non-negligible, then upon running the $\DAR_n(\alpha, K)$ system this proportion will decrease.
Contrast this with \cref{res:slow:high-to-low} where $\alpha > \alpha_c$ implied that there existed some $f$ which could be `supported'.

Initially stochastically dominate above by an $\Er(\beta(f_0),K)^n$ system.
If $\beta(f_0) > 1$, then next we choose $f_1$ with $p(f_0) < f_1 < f_0$ and stochastically dominate above by $\Er(\beta(f_1),K)^n$.
Iterating this, we eventually get the blocking level arbitrarily low (provided $K$ is sufficiently large).

We make the intuition above precise using \cref{alg:low:def,alg:low:stoch-dom}:
	the first makes definitions
and
	the second describes the stochastic domination and burn-in procedure.
	%
\end{Proof}



The following algorithm defines the sets and parameters used for stochastic domination.
Write
\[
	\eta
\cq
	\tfrac13 \inf\brb{ f - p(f) \midb \beta(f) \ge 1, \: f \in [0,1] }
\Quad{when}
	\alpha \ge \tfrac23.
\]
As $\max_{f \in [0,1]} \beta(f) = \beta(\tfrac12) = \tfrac32 \alpha$, there exists an $f \in [0,1]$ with $\beta(f) \ge 1$ if and only if $\alpha \ge \tfrac23$.
We are studying the regime $\alpha < \alpha_c$, so $p(f) < f$ for all $f \in [0,1]$ by definition of $\alpha_c$. Hence $\eta > 0$.


\begin{alg}[Definitions]
\label{alg:low:def}
	Assume that $\alpha \in (0,\alpha_c)$.
	Initialise $i \cq 0$.
	\begin{itemize}
		\item 
		Set
		$f_0 \cq \tfrac12$,
		$\alpha_0 \cq \beta(f_0)$,
		$\mfs_0 \cq \Omega$
		and
		$\mfb_0 \cq \Omega$.
		
		If $\alpha_0 \ge 1$, then increment $i \to i+1$ and \Proceed to the next step;
		otherwise \Stop.
		
		\item 
		Set
		$f_i \cq p(f_{i-1}) + 2 \eta$,
		$\alpha_i \cq \beta(f_i)$,
		\[
			\mfs_i
		&\cq
			\brb{ e \in \Omega \midb \varphie \le f_i }
		\quad	\text{and}
		\\
			\mfb_i
		&\cq
			\brb{ e \in \Omega \midb \prb[e]{ X^s \in \mfs_i \: \forall \, s \le n } \ge 1 - 1/n }.
		\]		
		If $\alpha_i \ge 1$, then increment $i \to i+1$ and \Repeat this step;
		otherwise \Stop.
	\end{itemize}
	If the algorithm terminates, then write $k$ for the number of steps it takes; otherwise set $k \cq \infty$.
	
	If $k < \infty$, then,
	for a parameter $\xi \in (0,1)$ to be chosen later,
	write
	\[
		\mfs_{k+1}
	\cq
		\mfs
	&\cq
		\brb{ e \in \Omega \midb \tfrac1n \abs{ \bra{ j \mid e_j > \xi K } } \le 2 \eps }
	\quad
		\text{and}
	\\
		\mfb_{k+1}
	\cq
		\mfb
	&\cq
		\brb{ e \in \Omega \midb \pr[e]{ X^s \in \mfs_{k+1} \: \forall \, s \le n } \ge 1 - 1/n }.
	\qedhere
	\]
\end{alg}

The following algorithm sets up the stochastic domination procedure, assuming $k < \infty$:
\begin{itemize}[noitemsep, topsep = \smallskipamount]
	\item 
	the sets and parameters used are from \cref{alg:low:def};
	
	\item 
	the legitimacy of the stochastic domination is provided by \cref{res:prelim:stoch-dom:gen}.
\end{itemize}
Observe that $\mfs_i$, and hence $\mfb_i$, is, for each $i$, a \textit{down-set}:
	if $x \le e$ and $e$ is in the set, then so is $x$.

\begin{alg}
\label{alg:low:stoch-dom}
Assume that $k < \infty$.
Set $T \cq 10 \log n$ and $t_i \cq i T$ for each $i = 0, ..., k+1$.

\medskip

\textit{Step $i$}, for $i \in \bra{1,...,k+1}$, takes the following form.

\begin{itemize}
	\item \emph{Stochastic Domination.}
	While in $\mfs_{i-1}$ stochastically dominate $X \lesssim E_{i-1} \sim \Er(\alpha_{i-1},K)^n$ and run for a time $T$, starting at time $t_{i-1} = (i-1)T$ and ending at time $t_i = i T$.
	
	\item \emph{Burn-In.}
	If $E_{i-1}^{t_i} \in \mfb_i$ (and hence $X^{t_i} \in \mfb_i$) and further $E_{i-1}^s \in \mfs_{i-1}$ (and hence $X^s \in \mfs_{i-1}$) for all $s \in [t_i, t_{k+1}]$, then continue;
	otherwise the burn-in phase fails and stop.
\end{itemize}
\end{alg}

\begin{Proof}[Outline]
	%
We upper bound the probability that $X^t \notin \mfb$ by the probability that the burn-in phase fails.
For $i = 1, ..., k$, each set $\mfs_i$ is defined by restricting the proportion of full links to be at most some value:
	this value is slightly large than the expected proportion for an $\Er(\alpha_{i-1}, K)^n$ system;
	namely,	the expected proportion is $p(f_{i-1})$ and the set requires a proportion at most $f_i = p(f_{i-1}) + 2 \eta$.
Finally, for step $k+1$, the dominating Erlang system is subcritical and so we can bring the proportion as low as we desire---we even impose slightly more.
In particular, $E_{i-1}$, and hence $X$, is highly likely to be in $\mfs_i$ in equilibrium.
\end{Proof}

The following two lemmas quantify this outline.
Their proofs are deferred to \cref{app:fast}.

\begin{lem}
\label{res:low:erlang_inv}
	For
		all $\beta \in (0,1)$,
		all $\eps \in (0,1)$
	and
		all $K$ sufficiently large,
	if $\pi_{\beta,K}$ is the invariant distribution of an $\Er(\beta, K)$ link, then
	\[
		\pi_{\beta, K}\rbb{ [\tfrac12(1 + \beta) K, K] }
	\le
		\eps.
	\]
\end{lem}


\begin{lem}
\label{res:low:high_inv-msre}
	For
		all $\alpha < \alpha_c$,
		all $K$ and $n$ sufficiently large
	and
		all $i = 1,...,k+1$,
	writing $\Pi_{i-1}$ for the invariant distribution of an $\Er(\alpha_{i-1}, K)^n$ system,
	we have
	\[
		\Pi_{i-1}\rbr{\mfb_i} \ge 1 - 1/n.
	\]
\end{lem}

\begin{Proof}[Proof of \cref{res:low:burnin}]
Observe that the statement is monotone in $\eps$: making $\eps$ larger can only decrease the probability.
Hence we may assume that $\eps$ is as small as we desire.

\smallskip

We check that the algorithm terminates when $\alpha < 1$.
If $\alpha < \tfrac23$, then it does in 0 steps.
Suppose $\alpha \ge \tfrac23$.
If $\alpha_{i-1} \ge 1$, then $p(f_{i-1}) \ge 0$; also $f_i = p(f_{i-1}) + 2\eta \le f_{i-1} - \eta$.
But $f_i \ge 0$ for all $i$.
Hence the algorithm does indeed terminate; further $k$ is a function only of $\alpha$ and satisfies $k \le \ceil{\tfrac12 \eta^{-1}}$.

Let $\eps > 0$, satisfying $2 \eps \le f_k \le \tfrac12$, so $\beta(2 \eps) \le \beta(f_k) = \alpha_k < 1$.
By definition of \cref{alg:low:def}, we have $\alpha_k < 1$.
Hence,
by \cref{res:low:erlang_inv},
for $K$ sufficiently large,
writing $\pi_k \cq (\pi_k(\ell))_{\ell=1}^K$ for the invariant distribution of a single $\Er(\alpha_k,K)$ link,
we have
\(
	\pi_k([\xi K, K]) \le \eps
\)
where $\xi \cq \tfrac12 (1 + \alpha_k) < 1$.

\smallskip

We now consider the probability that the burn-in phases succeeds.
For $i = 1,...,k+1$, write
\[
	1 - q_i
\cq
	\MIN{e \in \Omega}
	\prb[e]{ E_{i-1}^{T} \in \mfb_i }
\Quad{and}
	1 - q'_i
\cq
	\MIN{x \in \mfb_i}
	\prb[x]{ X^s \in \mfs_i \: \forall \, s \le t_{k+1} - t_i };
\]
note that $q'_{k+1} = 0$.
Write $Q \cq \max_{i=1,...,k+1} q_i$ and $Q' \cq \max_{i=1,...,k+1} q_i'$.
If the burn-in phase succeeds, then $X^{t_{k+1}} \in \mfb_{k+1}$.
Thus,
by using the Markov property, the stochastic domination and the union bound, taking worst-case scenarios at the start of each step,
we obtain
\[
	\MAX{x \in \Omega} \,
	\prb[x]{ X^{t_{k+1}} \notin \mfb_{k+1} }
\le
	(k+1) \max\bra{Q,Q'}.
\]
It remains to bound this maximum.
In particular, we set $C \cq 10(k+1)$, so then $t_{k+1} = C \log n$, and show that $\max\bra{Q,Q'} \le 2/n$; note also that $\mfb = \mfb_{k+1}$.
From this, the proposition follows.

By \cref{res:prelim:product-mixing}, the $1/n$ mixing time of an $\Er(\beta,K)^n$ system is at most $T = 10 \log n$ for all $\beta$, all $K \ge 4$ and all $n$ sufficiently large.
Write $\Pi_j$ for the invariant distribution of $\Er(\alpha_j,K)^n$, ie of $E_j$, for each $j = 0,...,k$.
By definition of the TV mixing time and \cref{res:low:high_inv-msre},
we obtain
\[
	q_i
=
	\MAX{e \in \Omega} \,
	\prb[e]{ E_{i-1}^{T} \notin \mfb_i }
\le
	1 - \Pi_{i-1}(\mfb_i) + 1/n
\le
	2/n
\Quad{for all}
	i = 1, ..., k+1.
\]
For all $i = 1, ..., k+1$, by definition of $\mfb_i$, we have $q'_i \le 1/n$ since $t_{k+1} - t_i \le t_{k+1} \asymp \log n \ll n$.
\end{Proof}

It remains to prove \cref{res:low:erlang_inv,res:low:high_inv-msre}.
These proofs are deferred to \cref{app:fast}.


\subsubsection{Variable Length Path Coupling}
\label{sec:low:var-coupling}

We now apply the variable length path coupling technique from \S\ref{sec:fast:var-coup}.
In this regime,
we have
\[
	\mfs
=
	\brb{ e \in \Omega \midb \tfrac1n \abs{ \bra{ j \mid e_j > \xi K } } \le 2 \eps },
\]
for some $\xi \in (0,1)$ sufficiently close to 1,
as in \cref{res:low:burnin},
with $\eps > 0$ to be specified later.

The reader is advised to recall the statement of the variable length path coupling result from \cref{res:var-length-coup}, as well as the notation and parameters defined therein.

In this part
we prove that the parameters satisfy the following properties.

\begin{prop}
	For
		$\eps \cq 10^{-3}$,
		all $\xi \in (\xi_0,1)$
	and
		all $K$ sufficiently large,
	in the scenario of \cref{res:var-length-coup},
	using the coupling from \cref{def:coupling} and $\mfs$ defined above,
	there exists a stopping time $\tau$ with
	$\gamma \le \tfrac23$ and $M \cq 6$ a valid choice.
\end{prop}

The remainder of this part is dedicated to proving this proposition.
See (below)
	\cref{def:low:st} for the definition of $\tau$
and
	\cref{res:low:gamma,res:low:M} for the bounds on $\gamma$ and $M$, respectively.

\medskip

For a link with different load on it in $X$ than $Y$, say the link is \textit{mismatched}; for a mismatched link, call the difference in load the \textit{mismatch distance}.
We first give an informal motivation for our stopping time and then the precise definition in \cref{def:low:st}.
We work on the event $\mfs$, which says that a proportion $1 - 2 \eps$ of the links are `well away from full'; call such links \textit{good}.

If the relative distance $d(X,Y) = 1$ and the mismatched link is good, then it is very likely that this additional call will end before the link becomes full: this is because a single link performs a type of mean-reverting random walk and $K$ is large.
For any bounded number of mismatched links with bounded, ie independent of both $n$ and $K$, mismatch distance, a similar result holds.

However, we cannot make this assumption on the first mismatched link: we need to choose a worst-case starting point, but from $\mfs$.
We simply wait for the first mismatched link to \textit{match}, ie stop being mismatched, and work on the event that any reroutings prior to this were to good links.
Given that the proportion of \textit{bad}, ie not good, links is at most $2 \eps$, this event is highly likely.
Further, in this time there will have been few reroutings.
If all this does happen, then we are in a state in which a bounded number of links are mismatched and all such links are good.
As above, all these links will match before any of them becomes full with high probability.

\medskip

We now make this precise and formal. We use the coupling from \cref{def:coupling}.
We first set out explicitly some terminology.
A link \textit{attempts to reroute} if it is full and a \textit{direct} call, ie one which has not been rerouted, arrives; this happens at rate $2\lambda$ for each full link.
Two other links are chosen uniformly and one of them accepts the call if both are not full, in which case the rerouting \textit{lands} in that link.
A call arrives \textit{indirectly} to a link if it lands in that link has a result of a rerouting.
A link \textit{accepts} both direct and indirect calls.

We define the stopping time $\tau$ to be the first time one of various events happens.
These are described in terms of stopping rules.
An event is \textit{triggered} when it happens.
For example, a stopping rule may be ``\Stop when the number of calls on link $1$ changes''. This event is triggered by a call's arrival to or departure from link $1$, before which the number of calls remains unchanged.

\begin{defn}
\label{def:low:st}
Consider $X, Y \sim \DAR_n(\alpha, K)$ using the coupling from \cref{def:coupling}.
Assume that $d(X^0, Y^0) = 1$ and $X^0 \ge Y^0$.
Consider the following stopping procedure.

\medskip

If at any point there is an arrival that is accepted in $Y$ but not in $X$, then \Stop; while this has not happened, use the following stopping rules.
\begin{ralist}
	\item \label{st:low:i}
	\Stop if one of the following occurs \emph{before} the original mismatched link matches:
	\begin{ralist}
		\item \label{st:low:i.a}
		the original mismatched link attempts to reroute 10 times;
		
		\item \label{st:low:i.b}
		one of the reroutings lands in a link with at least $\xi K$ calls on it in $X$;
		
		\item \label{st:low:i.c}
		a created mismatched link becomes full in either system.
	\end{ralist}
	
	\item \label{st:low:ii}
	Assume that the \Stop from \ref{st:low:i} is not triggered.
	\Stop when one of the following occurs:
	\begin{ralist}
		\item \label{st:low:ii.a}
		one of the created mismatched links becomes full in either system;
		
		\item \label{st:low:ii.b}
		the systems coalesce.
	\end{ralist}
\end{ralist}

Write $\tau$ for the time at which this procedure stops.
\end{defn}

We now bound this stopping time $\tau$ and determine the maximum relative distance $W$.
Write $\Gamma(k,1)$ for the $\Gamma$-distribution with shape $k$ and rate $1$.
Write $\mce(\mu)$ for the exponential distribution with rate $\mu$.
Then $\Gamma(k,1)$ is equivalent to the sum of $k$ independent $\mce(1)$-s.

\begin{lem}
\label{res:low:tau-W}
	We have $\tau \lesssim \Gamma(10,1)$ and $W = 10$.
\end{lem}

\begin{Proof}
While \ref{st:low:i.c} has not been triggered,
	no reroute can be accepted in $X$ but not in $Y$
and
	the only way a reroute can be accepted in $Y$ but not $X$ is if the original mismatched link is chosen.
Hence, in this case, the original mismatched link has mismatch distance at most 1.
This extra call (if it exists) departs at rate 1.
Hence the time taken for \ref{st:low:i} to be triggered is at most $\Gamma(1,1)$.
Further, the relative distance is at most 10 at any point before \ref{st:low:ii} is triggered;
if the process continues to \ref{st:low:ii}, then the relative distance is at most 9 at the time at which \ref{st:low:ii} starts.

Assume that the process continues to \ref{st:low:ii}, rather than stopping in \ref{st:low:i}.
By the same reasoning,
replacing \ref{st:low:i.c} with \ref{st:low:ii.a},
the set of full links is the same in $X$ as in $Y$,
and hence the relative distance cannot increase.
Since the relative distances is at most 9 initially,
it takes a time at most $\Gamma(9,1)$ for the systems to coalesce assuming that \ref{st:low:ii.a} has not been triggered.

In conclusion, $W = 10$ and $\tau \lesssim \Gamma(10,1)$, using the additive property of independent $\Gamma$-s.
	%
\end{Proof}

%

We now turn to bounding $\gamma$ and finding a suitable $M$, whose definitions we recall:
\begin{gather*}
	\gamma_0
=
	\MAX{(x,y) \in S} \,
	\exb[x,y]{ d(X^\tau,Y^\tau) \one{\msg[0,\tau]} }
\Quad{and}
	\gamma = \tfrac12(1+\gamma_0);
\\
	M
\Quad{is such that}
	\MAX{(x,y) \in S} \,
	\prb[x,y]{ \tau > M }
\le
	\tfrac12(1-\gamma_0)/W
\Quad{with}
	W = 10.
\end{gather*}
Recall that the event $\msg[0,\tau]$ means that $(X^s, Y^s) \in \mfs^2$ for all $s \le \tau$; in words, the number of links with at least $\xi K$ circuits in use is at most $2 \eps n$ in both $X$ and $Y$ for these times.

We first bound $\gamma_0$ and then use this to find a suitable $M$.

\begin{lem}
\label{res:low:gamma}
	For
		$\eps \cq 10^{-3}$,
		all $\xi \in (\xi_0,1)$
	and
		all $K$ sufficiently large,
	we have
	$\gamma \le \tfrac23$.
\end{lem}

\begin{Proof}
	%
By symmetry, without loss of generality we may assume that the originally mismatched link is link $1$ and that the extra call is in $X$: ie $d(X^0,Y^0) = 1$ and $X_1^0 = Y_1^0 + 1$.

Write $A$ for the complement of the event that $\tau$ is triggered by an arrival which is accepted to $Y$ but not to $X$, ie the first \Stop in \cref{def:low:st}. Then
\[
	\gamma_0 \le \gamma_0' + (W+1) \, \prb{ A^c, \, \msg[0,\tau] }
\Quad{where}
	\gamma_0'
\cq
	\MAX{(x,y) \in S} \,
	\exb[x,y]{ d(X^\tau,Y^\tau) \one{\msg[0,\tau]} \one{A} }.
\]

Suppose a link reroutes successfully in $Y$ but not in $X$. The choice of $i$ and $j$ in \cref{def:coupling} must then include the link $1$.
Write $F$ for the number of full links in $Y$ at this time.
The probability that the reroute lands in 1, rather than another non-full link of $Y$, is
\[
	\absb{ \brb{ (i,j) \subseteq [F]^2 \midb i = 1 } }
\big/
	\absb{ \brb{ (i,j) \subseteq [F]^2 \midb i = 1 \text{ or } j = 1 } }
=
	\tfrac12 \rbb{1 - \tfrac1{2F}}^{-1}
\ge
	\tfrac12.
\]
This is the conditional probability that a reroute is accepted to the original link in $Y$ but not in $X$ given that it is accepted to some link in $Y$ but not in $X$.
In this case, the relative distance decreases by 1; in the case that the reroute does not land in 1, the relative distance increases by 1.
Hence if the first \Stop is triggered,
then the expected change in relative distance is non-positive.

Write $\prref{st:low:i}$, $\prref{st:low:ii.a}$ and $\prref{st:low:ii.b}$ for the probabilities of the respective events in \cref{def:low:st},
	conditional on $A$ (ie that $\tau$ is not triggered by the first \Stop).
By inspection, if \tref{st:low:i} is triggered, then the relative distance is at most 11; if \tref{st:low:ii.a}, then also at most 11; if \tref{st:low:ii.b}, then 0.
Hence,
combined with the non-positivity on $A^c$,
the expected change in relative distance is at most
\[
	\rbb{
		\rbr{ 10 } \cdot \prrefb{st:low:i}
	+	\rbr{ 10 } \cdot \prrefb{st:low:ii.a}
	+	\rbr{ -1 } \cdot \prrefb{st:low:ii.b} }
\cdot
	\pr{ A }
+	\rbr{ 0 } \cdot \pr{A^c}
\]

Consider first \ref{st:low:i}. We consider the three subcases.
\begin{itemize}[leftmargin = 12mm]
	\item [\ref{st:low:i.a}]
	The matching can be caused by the extra call ending, or by a call arriving to the mismatched link and being accepted in $Y$ but not in $X$.
	Hence $\prref{st:low:i.a} \le 2^{-10}$.
	
	\item [\ref{st:low:i.b}]
	There are at most 10 reroutings before \ref{st:low:i.a} is triggered, and the probability that a given one of these lands on a link having at least $\xi K$ calls already on it is at most $2 \eps$, by definition of $\mfs$, as defined in \cref{res:low:burnin}.
	Hence $\prref{st:low:i.b} \le 20 \eps$ by the union bound.
	
	\item [\ref{st:low:i.c}]
	For neither \ref{st:low:i.a} nor \ref{st:low:i.b} to have been triggered, there must be at most 10 mismatched links and on each link the mismatch distance is at most 10.
	Write $q$ for the probability that an $\Er(\beta(2 \eps),K)$ link started from $\xi K$ hits $K$ before an independent $\Gamma(10,1)$ timer rings.
	Note also that no further mismatches can be made, because of the event $A$.
	Hence $\prref{st:low:i.c} \le 10 q$ by monotonicity and the union bound.
\end{itemize}
(For case \tref{st:low:i.a}, note that the event $A$ does not prohibit a call from being added to the original mismatched link in $Y$ but not in $X$; it only prohibits such events when the chosen link is not the original mismatched one.)
Combining these estimates with the union bound, we thus have
\[
	\prrefb{st:low:i} \le 2^{-10} + 20 \eps + 10 q.
\]

Now consider \ref{st:low:ii}.
By the same argument as for \ref{st:low:i.c}, we have
\[
	\prb{ \, \tref{st:low:ii.a} \midb \tref{st:low:ii} \, } \le 10 q,
\Quad{and hence}
	\prb{ \, \tref{st:low:ii.b} \midb \tref{st:low:ii} \, } \ge 1 - 10q.
\]

Next, by comparing a single Erlang link with a random walk on $\mbz$ which is biased towards its mean and using monotonicity we see that $q \to 0$ as $\Kinf$.

\smallskip

We must now consider $\pr{A}$.
Observe that if the set of full links is the same in $X$ as in $Y$, then any call added to $X$ is also added to $Y$ and vice versa. By definition of $\tau$, prior to $\tau$ there can be at most one link that is full in one system but not in the other; in particular, this is the original mismatched link, and it cannot be full in $Y$ but not in $X$ and if it is full in $X$ then it is one-off-full in $Y$ (prior to $\tau$).
While the mismatched link (link $1$) is full in $X$, $A$ is triggered if (and only if) a call arrives to a matched full link (ie any full link other than link $1$) and, in the notation of \cref{def:coupling}, $i$ is a non-full link other than link $1$ and $j = 1$.
Combining all this, we hence find that $A$ is triggered at rate $\lambda f(1-f) \le 2 \eps \lambda$, where $f$ is the current blocking level (and hence $f \le 2 \eps$).
Also, while it is full in $X$ and one-off-full in $Y$, link $1$ matches at rate $\lambda$, since if a call arrives directly and is accepted to $Y$ then it matches.
Hence, assuming $\eps \le \tfrac14$, we have
\[
	\prb{ A^c }
\le
	\prb{ \mce(4\eps\lambda) > \mce(\lambda) }
=
	4 \eps \lambda \big/ \rbb{ \lambda + 4\eps\lambda }
\le
	5 \eps.
\]

Combining all these, we find that the expected change in relative distance is at most
\[
	\rbb{ 10 \rbr{2^{-10} + 20\eps + 10q} + 100q + 10q - 1 } \rbb{ 1 - 5\eps }
+	(11) (5 \eps)
\le
	- \tfrac23
\]
for $K$ sufficiently large and $\eps$ sufficiently small, eg $\eps \cq 10^{-3}$.
Hence $\gamma_0 \le 1 - \tfrac23 = \tfrac13$ and $\gamma \le \tfrac23$.
	%
\end{Proof}

Given that $\gamma_0 \le \tfrac13$, we can now determine a permissible $M$.

\begin{lem}
\label{res:low:M}
	For
		$\eps \cq 10^{-3}$,
		all $\xi \in (\xi_0,1)$
	and
		all $K$ sufficiently large,
	we may take $M \cq 6$.
\end{lem}

\begin{Proof}
	From \cref{res:low:tau-W,res:low:gamma},
	we have
		$\tfrac12(1-\gamma_0)/W \ge \tfrac1{30}$
	and
		\[
			\prb{ \tau > M }
		\le
			\prb{ \Gamma(10,1) > M }
		\le
			10 \, \prb{ \mce(1) > M }
		=
			10 \, e^{-M}.
		\]
	It thus suffices to take any $M \ge \log(30 \cdot 10) = \log(300)$; since $\log(300) < 6$, we may take $M \cq 6$.
\end{Proof}

\subsubsection{Proof of Low-Blocking Mixing Theorem}
\label{sec:low:conclusion}

Now that we have defined the stopping time $\tau$, bounded $\gamma \le \tfrac23 \le e^{-1/3}$ and chosen $M \cq 6$, we can apply the variable length path coupling bound to prove our main theorem of the section, namely \cref{res:low}.
To this end, recall that
\(
	X, Y \sim \DAR_n(\alpha, K)
\)
under the coupling of \cref{def:coupling},
\[
	\mfb
=
	\brb{ e \in \Omega \midb \prb[e]{ X^s \in \mfs \: \forall \, s \le n } \ge 1 - 1/n }
\Quad{and}
	\msg[0,t]
=
	\brb{ (X^s, Y^s) \in \mfs^2 \: \forall \, s \in [0,t] }.
\]

\begin{Proof}[Proof of \cref{res:low}]
Plugging the expressions for $\gamma$ and $M$ from \cref{res:low:gamma,res:low:M}, respectively, into the variable length path coupling statement \cref{res:var-length-coup},
for $t \ge 0$,
we obtain
\[
	\prb{ X^t \neq Y^t, \, \msg[0,t] }
\le
	K n \rbb{\tfrac23}^{t/6 - 1}
=
	\tfrac32 K n e^{-t/18}.
\]
In particular, if we take $C_1 \cq 40$ and $t_1 \cq C_1 \log n$, then we obtain
\[
	\prb{ X^t \neq Y^t, \, \msg[0,t] }
=
	\tfrac32 K n^{-1-2/9}
\le
	1/n
\Quad{for all}
	t \ge t_1 = C_1 \log n.
\]
Next, \cref{res:low:burnin} gives us a constant $C_2$
depending only on $\alpha$
so that
\[
	\prb{ (X^{t_2},Y^{t_2}) \notin \mfb^2 }
\le
	\tfrac25 C_2 / n
\Quad{where}
	t_2 \cq C_2 \log n.
\]
by a union bound over $X$ and $Y$.
Finally, by definition of $\mfb$, we have
\[
	\MAX{(x,y) \in \mfb^2}
	\prb[x,y]{ \msg[0,t]^c }
\le
	2 / n
\Quad{where}
	t \cq t_1 + t_2 = (C_1 + C_2) \log n.
\]
Combining all these parts and applying the Markov property at time $t_1$ completes the proof.
\end{Proof}


\subsection{Fast Arrivals: $\alpha > 1$ and High-Blocking}
\label{sec:high}

In this subsection
we consider the `fast arrivals' regime, ie $\alpha > 1$.
We always use the coupling \mbpcoup from \cref{def:coupling}.
We prove the following theorem.

\begin{thm} \label{res:high}
	For all $\alpha > 1$,
	there exists a constant $C$ so that,
	for all $K$ and $n$ sufficiently large,
	if $X,Y \sim \DAR_n(\alpha,K)$, then,
	under the coupling \mbpcoup,
	for all $t \ge C \log n$,
	we have
	\[
		\MAX{(x,y) \in \Omega^2}
		\prb[x,y]{ X^t \neq Y^t }
	\le
		C / n
	=
		\oh1.
	\]
	Thus
	for all $K$ sufficiently large,
	all $\eps \in (0,1)$
	and all $n$ sufficiently large,
	we have
	$\tmix(\eps) \le C \log n$.
\end{thm}

The high-level idea of the proof will be the same as in the `slow arrivals' regime, but we shall require a different stopping time for the coupling.
Other than saying ``for all $\alpha > 1$'' rather than ``for all $\alpha < \alpha_c$'', the statements in this subsection
will be very similar to their counterparts in \S\ref{sec:low}; the subsection is even structured in a way that corresponding statements have the same number.

\smallskip

Used throughout this subsection, repeatedly, will be the notation $\beta(f)$ and $p(f)$:

\begin{itemize}[itemsep = 0pt, topsep = \smallskipamount]
	\item 
	$\beta(f) = \alpha(1 + 2f(1-f))$ is the effective traffic intensity when a proportion $f$ of links are full;
	
	\item 
	$p(f) = 1 - 1/\beta(f)$ is the equilibrium probability that an $\Er(\beta(f),K)$ link is full as $\Kinf$.
\end{itemize}


We prove this theorem via a sequence of lemmas:
	in \S\ref{sec:high:burnin} we describe the burn-in phase;
	in \S\ref{sec:high:var-coupling} we describe and apply the variable length path coupling;
	finally we conclude in \S\ref{sec:high:conclusion}.

\subsubsection{Burn-In Phase}
\label{sec:high:burnin}

Set $\varphi_c \cq 1 - 1/\sqrt2$; one can check directly that if $\alpha > 1$ then $p(\varphi_c) > \varphi_c$.
Thus when $\alpha > 1$ the system can `support' $\varphi_c$-blocking.
The aim of this part
is to prove the following burn-in proposition.

\begin{prop}
\label{res:high:burnin}
	For all $\alpha > 1$,
	there exist constants $C > 0$ and $\xi_0 \in (0,1)$ so that,
	for all $\xi \in (0,\xi_0)$,
	setting $t \cq C \log n$,
	\[
		\mfs &\cq \brb{ e \in \Omega \midb \varphie \ge \varphi_c + 2 \xi }
	\quad
		\text{and}
	\\
		\mfb &\cq \brb{ e \in \Omega \midb \pr[e]{ X^s \in \mfs \: \forall \, s \le n } \ge 1 - 1/n },
	\]
	for all $K$ and $n$ sufficiently large (depending on $\alpha$),
	if $X \sim \DAR_n(\alpha, K)$, then
	\[
		\MAX{x \in \Omega} \,
		\prb[x]{ X^t \notin \mfb }
	\le
		\tfrac15 C / n
	=
		\oh1.
	\]
\end{prop}

Since $\alpha > 1$, for any $f \in [0,1]$, we have $\beta(f) \ge \beta(0) = \alpha > 1$.
Also observe that
\[
	p(f) = 1 - 1/\beta(f) > \varphi_c = 1 - 1/\sqrt2
\Quad{if and only if}
	\beta(f) > \sqrt2.
\]
We say that an $\Er(\beta,K)$ link is \textit{very supercritical} if $\beta > \sqrt2$, and similarly for an $\Er(\beta,K)^n$ system.

\begin{Proof}[Intuition]
\renewcommand{\qedsymbol}{\ensuremath{\triangle}}
By definition of $\alpha_c$, when $\alpha > \alpha_c$, there exists an $f \in (0,1)$ so that $p(f) > f$.
Thus the expected proportion of full links in an $\Er(\beta(f), K)^n$ system in equilibrium is \textit{more than $f$} (in the limit $\Kinf$).
Intuitively, this suggests that the system can `support' such $f$-blocking:
	if we start a $\DAR_n(\alpha, K)$ system from a state with proportion $f$ blocked then (typically) the proportion initially increase.
Thus if the proportion blocked initially, call it $f_0$, satisfies $p(f_0) > f_0$, then upon running the $\DAR_n(\alpha, K)$ system this proportion will increase.

Initially we stochastically dominate from below by an $\Er(\beta(f_0), K)^n$ system.
If $f_0 < \varphi_c$, then next we choose $f_1$ with $f_0 < f_1 < p(f_0)$ and stochastically dominate below by $\Er(\beta(f_1), K)^n$.
Iterating this, we eventually get the blocking level above $\varphi_c$ (provided $K$ is sufficiently large).

We make the intuition above precise using \cref{alg:high:def,alg:high:stoch-dom}:
	the first makes definitions
and
	the second describes the stochastic domination and burn-in procedure.
We only use these when $\alpha \le \sqrt2$:
when $\alpha > \sqrt2$,
we immediately have $p(0) > \varphi_c$ and simply use $X \gtrsim \Er(\alpha,K)^n$.
\end{Proof}


We separate the proof into two cases:
	$\alpha > \sqrt2$ and $\alpha \le \sqrt2$;
	always $\alpha > 1$.
The former case is significantly easier; there is no need for an iterative burn-in period like there was for $\alpha < \alpha_c$.
Our target is to obtain a blocking level larger than $\varphi_c$; this is achieved by $\Er(\alpha, K)^n$, which trivially stochastically dominates below $\DAR_n(\alpha, K)$; recall \cref{res:prelim:stoch-dom:cor}.
When $\alpha \le \sqrt2$, we need an iterative burn-in period to get the blocking up to $\varphi_c$.

We give the proof for $\alpha > \sqrt2$ immediately.
The following lemma quantifies the probability that $\Er(\alpha, K)^n$ is in $\mfb$ in equilibrium.
Its proof is deferred to \cref{app:fast}.

\begin{lem}
\label{res:high:claim:a>sqrt2:inv-msre}
	For all $\alpha > \sqrt2$,
	there exists $\xi_0 \in (0,1)$ so that,
	for
		all $\xi \in (0, \xi_0)$
	and
		all $K$ and $n$ sufficiently large,
	we have
	\[
		\Pi_0(\mfb)
	\ge
		1 - 1/n.
	\]
\end{lem}

\begin{Proof}[Proof of \cref{res:high:burnin} when $\alpha > \sqrt2$]
Since $\alpha > \sqrt2$ implies that $p(0) > \varphi_c$, there exists an $\xi > 0$ so that $p(0) - 2 \xi \ge \varphi_c + 2 \xi$.
Recall the definitions
\[
	\mfs &\cq \brb{ e \in \Omega \midb \varphie \ge \varphi_c + 2 \xi }
\quad	\text{and}
\\
	\mfb &\cq \brb{ e \in \Omega \midb \prb[e]{ X^s \in \mfs \: \forall \, s \le n } \ge 1 - 1/n }.
\]

By \cref{res:prelim:stoch-dom:cor}, we can stochastically dominate $X \gtrsim E_0 \sim \Er(\alpha,K)^n$.
By \cref{res:prelim:product-mixing}, the $1/n$ mixing time of an $\Er(\beta,K)^n$ system is at most $t \cq 10 \log n$ for all $\beta$, all $K \ge 4$ and all $n$ sufficiently large.
Write $\Pi_0$ for the invariant distribution of $\Er(\alpha,K)^n$, ie of $E_0$.
By definition of the TV mixing time and \cref{res:high:claim:a>sqrt2:inv-msre},
we obtain
\[
	\MAX{x \in \Omega} \,
	\prb[x]{ X^t \notin \mfb }
\le
	\MAX{e \in \Omega} \,
	\prb[e]{ E_0^t \notin \mfb }
\le
	1 - \Pi_0(\mfb) + 1/n
\le
	2/n.
\]
Taking $C \cq 10$, this completes the proof when $\alpha > \sqrt2$.
\end{Proof}

The following algorithm defines the sets and parameters used for stochastic domination.
Write
\[
	\eta \cq \tfrac12 \inf\brb{ p(f) - f \midb f \in [0,\varphi_c] }
\Quad{when}
	\alpha \le \sqrt2.
\]
Note that $\alpha \le \sqrt2$ implies that
\(
	\max_{f \in [0,1]}
	\beta(f)
\le
	\beta(\tfrac12)
\le
	\tfrac32\sqrt2
\)
for all $f \in [0,1]$.
We are studying the regime $\alpha > 1$, so $p(f) > f$ for all $f \in [0,\varphi_c]$ by direct calculation. Hence $\eta > 0$.

\begin{alg}[Definitions]
\label{alg:high:def}
	Assume that $\alpha \in (1, \sqrt2]$.
	Set $\alpha_{-1} \cq \tfrac32 \sqrt2$.
	Initialise $i \cq 0$.
	\begin{itemize}
		\item 
		Set
		$f_0 \cq 0$,
		$\alpha_0 \cq \beta(f_0) = \alpha$,
		\[
			\mfs_0 &\cq \brb{ e \in \Omega \midb \varphie \le \tfrac23 }
		\quad	\text{and}
		\\
			\mfb_0 &\cq \brb{ e \in \Omega \midb \pr[e]{ X^s \in \mfs_0 \: \forall \, s \le n } \ge 1 - 1/n }.
		\]		
		(Note that $f_0 = 0$.)
		Increment $i \to i+1$ and \Proceed to the next step.
		
		\item 
		Set
		$f_i \cq \min\bra{p(f_{i-1}) - \eta, \tfrac13}$,
		$\alpha_i \cq \beta(f_i)$,
		\[
			\mfs_i &\cq \brb{ e \in \Omega \midb \varphie \ge f_i }
		\quad	\text{and}
		\\
			\mfb_i &\cq \brb{ e \in \Omega \midb \pr[e]{ X^s \in \mfs_i \: \forall \, s \le n } \ge 1 - 1/n }.
		\]	
		If $f_i \le \varphi_c = 1 - 1/\sqrt2$, then increment $i \to i+1$ and \Repeat this step; otherwise \Stop.
	\end{itemize}
	If the algorithm terminates, then write $k$ for the number of steps it takes; otherwise set $k \cq \infty$.
	
	If $k < \infty$, then choose $\xi > 0$ so that $f_k \ge \varphi_c + 4 \xi$, ie $f_k - 2 \xi \ge \varphi_c + 2 \xi$.
\end{alg}

The following algorithm sets up the stochastic domination procedure, assuming $k < \infty$:
\begin{itemize}[noitemsep, topsep = \smallskipamount]
	\item 
	the sets and parameters used are from \cref{alg:high:def};
	
	\item 
	the legitimacy of the stochastic domination is provided by \cref{res:prelim:stoch-dom:gen}.
\end{itemize}
Observe that $\mfs_i$, and hence $\mfb_i$, is, for each $i$, an \textit{up-set}:
	if $x \ge e$ and $e$ is in the set, then so is $x$.

\begin{alg}
\label{alg:high:stoch-dom}
Assume that $k < \infty$.
Set $T \cq 10 \log n$ and $t_i \cq (i+1) T$ for each $i = 0, ..., k$.

\medskip

\textit{Step $0$} takes the following form. (Note that $\beta(f) \le \alpha_{-1}$ for all $f$.)
\begin{itemize}
	\item \emph{Stochastic Domination.}
	Stochastically dominate $X \lesssim E_{-1} \sim \Er(\alpha_{-1},K)^n$ and run for a time $T$, starting at time $0$ and ending at time $t_0 = T$.
	
	\item \emph{Burn-In.}
	If $E_{-1}^{t_0} \in \mfb_0$ (and hence $X^{t_0} \in \mfb_0$) and further $E_{-1}^s \in \mfs_0$ (and hence $X^s \in \mfs_0$) for all $s \in [t_0,t_k]$, then continue;
	otherwise the burn-in phase fails and stop.
\end{itemize}

\textit{Step $i$}, for $i \in \bra{1, ..., k}$, takes the following form.
\begin{itemize}
	\item \emph{Stochastic Domination.}
	While in $\mfs_{i-1}$ stochastically dominate $X \lesssim E_{i-1} \sim \Er(\alpha_{i-1}, K)^n$ and run for a time $T$, starting at time $t_{i-1} = iT$ and ending at time $t_i = (i+1)T$.
	
	\item \emph{Burn-In.}
	If $E_{i-1}^{t_i} \in \mfb_i$ (and hence $X^{t_i} \in \mfb_i$) and further $E_{i-1}^s \in \mfs_{i-1}$ (and hence $X^s \in \mfs_{i-1}$) for all $s \in [t_i, t_k]$, then continue;
	otherwise the burn-in phase fails and stop.
	\qedhere
\end{itemize}
\end{alg}

\begin{Proof}[Outline]
We upper bound the probability that $X^t \notin \mfb$ by the probability that the burn-in phase fails.
Consider working in the restricted space $\mfs_0$, where the proportion of full links is at most $\tfrac23$.
For $i = 1, ..., k$, each set $\mfs_i$ is defined by restricting the proportion of full links to be at least some value:
	this value is slightly smaller than the expected proportion for an $\Er(\alpha_{i-1}, K)^n$ system;
	namely,	the expected proportion is $p(f_{i-1})$ and the set requires a proportion at least $f_i = \min\bra{p(f_{i-1}) - 2 \eta, \tfrac13}$.
In particular, $E_{i-1}$, and hence $X$, is highly likely to be in $\mfs_i$ in equilibrium.
\end{Proof}

The following lemma quantifies this outline.
Its proof is deferred to \cref{app:fast}.

\begin{lem}
\label{res:high:claim:a<sqrt2:inv-msre}
	For all $\alpha \in (1,\sqrt2]$
	and all $n$ sufficiently large,
	for each $i = 0, 1,...,k$,
	we have
	\[
		\Pi_{i-1}(\mfb_i) \ge 1 - 1/n.
	\]
\end{lem}

\begin{Proof}[Proof of \cref{res:high:burnin} when $\alpha \le \sqrt2$]
Observe that the statement is monotone in $\xi$: making $\xi$ larger can only decrease the probability.
Hence we may assume that $\xi$ is as small as we desire.

\smallskip

We check that the algorithm terminates when $\alpha \in (1, \sqrt2]$.
As $\alpha > 1$, we have $\beta(f) > 1$ and hence $p(f) > 0$ for all $f \in [0,1]$.
If $f_{i-1} \le \varphi_c$, then $f_i \ge p(f_{i-1}) - \eta$. But $p(f) - \eta \ge f + \eta$ for all $f \in [0,\varphi_c]$, and hence $f_i \ge f_{i-1} + \eta$.
Also, $f_i \le \tfrac13$ and $\varphi_c = 1 - 1/\sqrt2 < \tfrac13$.
Hence the algorithm does indeed terminate; further, $k$ is a function only of $\alpha$ and satisfies $k \le \ceil{\tfrac13 \eta^{-1}}$.

We now consider the probability that the burn-in phases succeeds.
For $i = 0,1,...,k$, write
\[
	1 - q_i
=
	\MIN{e \in \Omega} \,
	\prb[e]{ E_{i-1}^T \in \mfb_i }
\Quad{and}
	1 - q'_i
=
	\MIN{x \in \mfb_i} \,
	\prb[x]{ X^s \in \mfs_i \: \forall \, s \le t_k - t_i };
\]
note that $q'_k = 0$.
Write $Q \cq \max_{i=0,...,k} q_i$ and $Q' \cq \max_{i=0,...,k} q'_i$.
If the burn-in phase succeeds, then $X^{t_{k+1}} \in \mfb_{k+1}$.
Thus,
by using the Markov property, the stochastic domination and the union bound, taking worst-case scenarios at the start of each step,
we obtain
\[
	\MAX{x \in \Omega} \,
	\prb[x]{ X^{t_{k+1}} \notin \mfb_{k+1} }
\le
	(k+1) \max\bra{Q,Q'}.
\]
It remains to bound this maximum.
In particular, we set $C \cq 10(k+1)$, so then $t_k = C \log n$, and show that $\max\bra{Q,Q'} \le 2/n$; note also that $\mfb \supseteq \mfb_{k+1}$.
From this, the proposition follows.

By \cref{res:prelim:product-mixing}, the $1/n$ mixing time of an $\Er(\beta,K)^n$ system is at most $T = 10 \log n$ for all $\beta$, all $K \ge 4$ and all $n$ sufficiently large.
Write $\Pi_j$ for the invariant distribution of $\Er(\alpha_j,K)^n$, ie of $E_j$, for each $j = -1,0,...,k$.
By definition of the TV mixing time and \cref{res:high:claim:a<sqrt2:inv-msre}, we obtain
\[
	q_i
=
	\MAX{e \in \Omega}
	\prb[e]{ E_{i-1}^T \notin \mfb_k }
\le
	1 - \Pi_{i-1}(\mfb_i) + 1/n
\le
	2/n
\Quad{for all}
	i = 0,...,k.
\]
For all $i = 0, ..., k$, by definition of $\mfb_i$, we have $q_i \le 1/n$ since $t_k - t_i \le t_k \asymp \log n \ll n$.
\end{Proof}


%

It remains to prove \cref{res:high:claim:a>sqrt2:inv-msre,res:high:claim:a<sqrt2:inv-msre}.
These proofs are deferred to \cref{app:fast}.


\subsubsection{Variable Length Path Coupling}
\label{sec:high:var-coupling}

We now apply the variable length path coupling technique from \S\ref{sec:fast:var-coup}.
In this regime, we have
\[
	\mfs
=
	\brb{e \in \Omega \midb \varphie \in [\varphi_c + 2 \xi,\tfrac23] },
\]
for some $\xi \in (0,1)$ sufficiently close to 0,
as in \cref{res:high:burnin}

The reader is advised to recall the statement of the variable length path coupling result from \cref{res:var-length-coup}, as well as the notation and parameters defined therein.

In this part
we prove that the parameters satisfy the following properties.

\begin{prop}
	For
		all $\xi \in (0,\xi_0)$
	and
		all $K$ sufficiently large,
	in the scenario of \cref{res:var-length-coup},
	using the coupling from \cref{def:coupling} and $\mfs$ defined above,
	there exists a stopping time $\tau$ with
	$\gamma \le 1 - \xi$ and $M \cq \log(1/\xi)$ a valid choice.
\end{prop}

The remainder of this part is dedicated to proving this proposition.
See (below)
	\cref{def:high:st} for the definition of $\tau$
and
	\cref{res:high:gamma,res:high:M} for the bounds on $\gamma$ and $M$, respectively.

\medskip

For a link with different load on it in $X$ than $Y$, recall that we say the link is \textit{mismatched}; for a mismatched link, call the difference in load the \textit{mismatch distance}.
We first give an informal motivation for our stopping time, and then the precise definition (in \cref{def:high:st}).
We work on the event $\mfs$, which says that the proportion blocked is always at least $\varphi_c + 2 \xi$.

Start with $d(X^0,Y^0) = 1$.
Our stopping time is simple: wait for the first reroute attempt in the mismatched link.
The relative distance stays the same, decreases by 1 or increases by 1; we show that having a proportion blocked greater than $\varphi_c$ is sufficient for the expected distance to~decrease.

\medskip

We now make this precise and formal. We define $\tau$ via a set of stopping rules.

\begin{defn}
\label{def:high:st}
Consider $X, Y \sim \DAR_n(\alpha, K)$ using the coupling from \cref{def:coupling}.
Assume that $d(X^0, Y^0) = 1$ and $X^0 \ge Y^0$.
Consider the following stopping procedure.

\medskip

\Stop when one of the following events occurs:
\romannumbering
\begin{ralist}
	\item \label{st:high:i}
	the `extra' call ends or a call is added to the mismatched link in $Y$ but not $X$ via a rerouting;
	
	\item \label{st:high:ii}
	the original mismatched link attempts to reroute (in $X$, but not in $Y$);
	
	\item \label{st:high:iii}
	a link successfully reroutes in $Y$ but not in $X$.
	\qedhere
\end{ralist}
\end{defn}

We now bound this stopping time $\tau$ and determine the maximum distance $W$.
Write $\mce(\mu)$ for the exponential distribution with rate $\mu$.


\begin{lem}
\label{res:high:tau-W}
	We have $\tau \lesssim \mce(1)$ and $W = 1$.
\end{lem}

\begin{Proof}
While \ref{st:high:ii} and \ref{st:high:iii} have not been triggered,
	no reroute can be accepted in $X$ but not in $Y$
and
	the only way a reroute can be accepted in $Y$ but not $X$ is if the original mismatched link is chosen and added to.
Hence the original mismatched link has mismatch distances at most 1.
This extra call (if it exists) departs at rate 1.
Hence, in this case, the time taken for \ref{st:low:i} is at most $\mce(1)$.

For the same reasoning, $d(X^s,Y^s) = 1$ for all $s < \tau$. Thus $W = 1$.
\end{Proof}

%

We now turn to bounding $\gamma$ and finding a suitable $M$, whose definitions we recall:
\begin{gather*}
	\gamma_0
\cq
	\MAX{(x,y) \in S}
	\exb[x,y]{ d(X^\tau,Y^\tau) \one{\msg[0,\tau]} }
\Quad{and}
	\gamma \cq \tfrac12(1+\gamma_0);
\\
	M
\Quad{is such that}
	\MAX{(x,y) \in S}
	\prb[x,y]{ \tau > M }
\le
	\tfrac12(1-\gamma_0)/W
\Quad{with}
	W = 1.
\end{gather*}
Recall that the event $\msg[0,\tau]$ means that $(X^s, Y^s) \in \mfs^2$ for all $s \le \tau$; in words, the number of full links is at least $\varphi_c + 2 \xi$ in both $X$ and $Y$ for these times.

First we bound $\gamma_0$, and then use this to find a suitable $M$.

\begin{lem}
\label{res:high:gamma}
	For
		all $\xi \in (0,\xi_0)$
	and
		all $K$ sufficiently large,
	we have
	$\gamma \le 1 - \xi \le e^{-\xi}$.
\end{lem}

\begin{Proof}
	%
By symmetry, without loss of generality we may assume that the originally mismatched link is link $1$ and that the extra call is in $X$: ie $d(X^0,Y^0) = 1$ and $X_1^0 = Y_1^0 + 1$.

By inspection, if \tref{st:high:i} is triggered, then the systems coalesce, and so the change in relative distance is $-1$ necessarily (and hence $-1$ in expectation).

Suppose a link reroutes successfully in $Y$ but not in $X$. The choice of $i$ and $j$ in \cref{def:coupling} must then include the link $1$.
Write $F$ for the number of full links in $Y$ at this time.
The probability that the reroute lands in 1, rather than another non-full link of $Y$, is
\[
	\absb{ \brb{ (i,j) \subseteq [F]^2 \midb i = 1 } }
\big/
	\absb{ \brb{ (i,j) \subseteq [F]^2 \midb i = 1 \text{ or } j = 1 } }
=
	\tfrac12 \rbb{1 - \tfrac1{2F}}^{-1}
\ge
	\tfrac12.
\]
This is the conditional probability that a reroute is accepted to the original link in $Y$ but not in $X$ given that it is accepted to some link in $Y$ but not in $X$.
In this case, the relative distance decreases by 1; in the case that the reroute does not land in 1, the relative distance increases by 1.
Hence if \tref{st:high:iii} is triggered, then the expected change in relative distance is non-positive.

If \tref{st:high:ii} is triggered, then we have four cases:
\begin{ralist}[leftmargin = 12mm, label = \upshape(ii.\textit{\alph*})]
	\item \label{st:high:ii.a}
	the reroute is successful in $X$ and the call is not added in $Y$,
	giving $d(X^\tau,Y^\tau) = 2$;
	
	\item \label{st:high:ii.b}
	the reroute is successful in $X$ and the call is added in $Y$,
	giving $d(X^\tau,Y^\tau) = 1$;
	
	\item \label{st:high:ii.c}
	the reroute is unsuccessful in $X$ and the call is not added in $Y$,
	giving $d(X^\tau,Y^\tau) = 1$;
	
	\item \label{st:high:ii.d}
	the reroute is unsuccessful in $X$ and the call is added in $Y$, 
	giving $d(X^\tau,Y^\tau) = 0$.
\end{ralist}
Write $\prref{st:high:ii.a}$, $\prref{st:high:ii.b}$, $\prref{st:high:ii.c}$ and $\prref{st:high:ii.d}$ for the probabilities of the above events,
	conditional that $\tau$ is triggered by $\tref{st:high:ii}$.
The probability that the call is added to $Y$ is $\tfrac12$, independent of what happens in $X$.
Let $f \cq \varphie[X^\tau]$ denote the proportion blocked in $X$ at the time $\tau$.
Then
\[
	\prref{st:high:ii.a}
=
	\prref{st:high:ii.b}
=
	\tfrac12 (1-f)^2
\Quad{and}
	\prref{st:high:ii.c}
=
	\prref{st:high:ii.d}
=
	\tfrac12 \rbb{ 1 - (1-f)^2 }.
\]
If $X^\tau \in \mfs$ then $f \ge 1 - 1/\sqrt2 + 2 \xi$, by definition.
Hence the expected change in relative distance is 
\[
	\tfrac12(1-f)^2 - \tfrac12\rbb{1 - (1-f)^2}
=
	(1-f)^2 - \tfrac12
\le
	\rbb{1/\sqrt2 - 2 \xi}^2 - \tfrac12
=
	- 2 \sqrt2 \xi + 4 \xi^2
\le
	- \tfrac52 \xi,
\]
with the final inequality holding if $\xi$ is small enough ($\xi \le \tfrac1{40}$ is sufficient).

Combining the three cases, we see that the expected change in relative distance is at most
\[
	\rbr{ -1 } \cdot \prrefb{st:high:i}
+	\rbr{ -\tfrac52 \xi } \cdot \prrefb{st:high:ii}
+	\rbr{ 0} \cdot \prrefb{st:high:iii}
\le
	- \tfrac52 \xi \rbb{ \prrefb{st:high:i} + \prref{st:high:ii} }
=
	- \tfrac52 \xi \rbb{ 1 - \prrefb{st:high:iii} }.
\]
We now wish to upper bound $\prref{st:high:iii}$.
If the original mismatched link is both not full in $X$ and in $Y$, then the set of full links is the same in $X$ as in $Y$, and so \tref{st:high:iii} cannot occur: it cannot be full in $Y$ but not full in $X$ (prior to $\tau$).
While the mismatched link is full in $X$, \tref{st:high:ii} is triggered at rate $2 \lambda$ while \tref{st:high:iii} is at rate $r(f) = 2 \lambda f (1-f) \le \tfrac12 \lambda$. Hence we see that
\[
	\prrefb{st:high:iii}
\le
	\prb{ \mce(\tfrac12 \lambda) < \mce(2 \lambda) }
=
	\tfrac12 \big/ \rbb{2 + \tfrac12}
=
	\tfrac15.
\]
Hence the expected change in relative distance is at most $-\tfrac52 \xi \cdot \tfrac45 = - 2 \xi$, ie $\gamma_0 \le 1 - 2 \xi$.
	%
\end{Proof}

Given that $\gamma_0 \le 1 - 2 \xi$, we can now determine a permissible $M$.

\begin{lem}
\label{res:high:M}
	For
		all $\xi \in (0,\xi_0)$
	and
		all $K$ sufficiently large,
	we may take $M \cq \log(1/\xi)$.
\end{lem}

\begin{Proof}
	From \cref{res:high:tau-W,res:high:gamma},
	we have
		$\tfrac12(1-\gamma_0)/W \ge \xi$
	and
	\[
		\prb{ \tau > M }
	\le
		\prb{ \mce(1) > M }
	=
		e^{-M}.
	\]
	It thus suffices to take $M \cq \log(1/\xi)$.
\end{Proof}

\subsubsection{Proof of High-Blocking Mixing Theorem}
\label{sec:high:conclusion}

Now that we have defined the stopping time $\tau$, bounded $\gamma \le e^{-\xi}$ and chosen $M \cq \log(1/\xi)$, we can apply the variable length path coupling bound to prove our main theorem of the section, namely \cref{res:high}.
To this end, recall that
\(
	X, Y \sim \DAR_n(\alpha, K)
\)
under the coupling of \cref{def:coupling},
\[
	\mfb
=
	\brb{ e \in \Omega \midb \prb[e]{ X^s \in \mfs \: \forall \, s \le n } \ge 1 - 1/n }
\Quad{and}
	\msg[0,t] = \brb{ (X^s, Y^s) \in \mfs^2 \: \forall \, s \in [0,t] }.
\]

\begin{Proof}[Proof of \cref{res:high}]
Plugging the expressions for $\gamma$ and $M$ from \cref{res:high:gamma,res:high:M}, respectively into the variable length path coupling statement \cref{res:var-length-coup},
for $t \ge 0$,
we obtain
\[
	\prb{ X^t \neq Y^t, \, \msg[0,t] }
\le
	K n e^{-\xi(t/\log(1/\xi) - 1)}
\le
	2 K n e^{-\xi t / \log(1/\xi)},
\]
assuming $\xi \le \log 2$.
Here $\xi \in (0,1)$ is a constant which is sufficiently small, in a manner depending only on $\alpha$.
In particular, if we take $C_1 \cq 2 \xi^{-1} \log(1/\xi)$ and $t_1 \cq C_1 \log n$, then we obtain
\[
	\prb{ X^t \neq Y^t, \, \msg[0,t] }
\le
	2 K n^{-2}
\le
	n^{-1}
\Quad{for all}
	t \ge t_1 = C_1 \log n.
\]
Next, \cref{res:high:burnin} gives us a constant $C_2$
depending only on $\alpha$
so that
\[
	\prb{ (X^{t_2}, Y^{t_2}) \notin \mfb^2 }
\le
	\tfrac25 C_2 / n
\Quad{where}
	t_2 \cq C_2 \log n,
\]
by a union bound over $X$ and $Y$.
Finally, by definition of $\mfb$, we have
\[
	\MAX{(x,y) \in \mfb^2}
	\prb[x.y]{ \msg[0,t]^c }
\le
	2 / n
\Quad{where}
	t \cq t_1 + t_2 = (C_1 + C_2) \log n.
\]
Combining all these parts and applying the Markov property at time $t_1$ completes the proof.
\end{Proof}


\section{Dynamic Alternative Routing with Retries}
\label{sec:retries}

\subsection{Introduction}

We now consider a generalisation of the high-level model, as described at the start of the paper.
We call it \textit{dynamic alternative routing with retries}; it is also known as \textit{multiple alternatives}.

In the original model, if a call arrives asking for the link connecting stations $\alpha$ and $\beta$, if this link is full (ie at capacity), then a third station $\gamma$ is selected uniformly at random amongst the remaining stations:
	if there is free capacity on both $\alpha\gamma$ and $\gamma\beta$ then the call is held on these two links simultaneously;
	otherwise the call is simply declined (ie lost).
We call that act of picking a third station and attempting to route via it a \textit{retry}. So the above model has one retry.

We now generalise this model: instead of declining (losing) the call if the first reroute attempt fails, we allow $\rho$ (independent) attempts, where $\rho \in \mbn$; if all $\rho$ retries fail, then the call is declined (ie lost).
(Of course, if, say, the 3rd retry is successful, then we accept the call and stop: we do not do the remaining $\rho-3$ retries.)
Taking $\rho = 1$ reduces to the original model.

\smallskip

We show that the overall behaviour of this system, for general $\rho$, exhibits the same properties as for $\rho = 1$; in particular, we have an interim regime $\alpha \in (\alpha_c,1)$ with metastability.
The reasons for this are the same as in the $\rho = 1$ case:
	even with the traffic intensity $\alpha < 1$,
	the \emph{effective} traffic intensity (taking into account the fact that rerouted calls hold \emph{two} circuits) may be larger than $1$.

We denote this system by $\DAR_n^\rho(\alpha,K)$, for $n$ links, each of capacity $K$, traffic intensity $\alpha$ and $\rho$ rerouting attempts; we also write $\lambda \cq \alpha K$.

\subsubsection{Model Set-Up and Main Theorem}

We now describe our model and then motivate the details after (as we did in the introduction).
Fix $\rho \in \mbn$.
Suppose the system is in state $x$; write $f = \tfrac1n \sum_{j=1}^n \one{ x_j = K }$ for the proportion of links which are full.
To the full links, no calls arrive.
To the non-full links, calls arrive (independently amongst links) at rate $\beta_\rho(f) K$ where
\[
	\beta_\rho(f) \cq \alpha\rbb{1 + r_\rho(f) }
\Quad{and}
	r_\rho(f) \cq 2 f \rbb{ 1 - \rbb{1 - (1-f)^2}^\rho } / (1 - f).
\]
Also write $p_\rho(f)$ for the equilibrium probability that a single $\Er(\beta_\rho(f),K)$ is full when $\beta_\rho(f) \ge 1$:
\[
	p_\rho(f) \cq 1 - 1/\beta_\rho(f).
\]
Again, it is not immediately clear that this is the correct rate function for the model, in the same way that it was not clear in \S\ref{sec:intro:model}; we show in the next part
that this is the correct rate.
We also define the critical $\alpha$ analogously to before:
\[
	\alpha_c(\rho)
\cq
	\hphantom{\sup}
	\mathllap{\inf}
	&\brb{ \alpha \in (0,\infty) \midb \exists \, f \in [0,1] \ST p_\rho(f) > f }
\\
=
	\sup
	&\brb{ \alpha \in (0,\infty) \midb p_\rho(f) < f \: \forall \, f \in [0,1] }.
\]

\medskip

The main theorem of this section is the following mixing time result.
(It was stated as \cref{res:intro:retry} in the introduction; we recall it here for convenience.)
Over the next three subsections
we \emph{explain} how to prove it.
We do not give all the details, but rather explain which details differ from the no-retries ($\rho = 1$) case, and explain how to overcome these.

\begin{thm}
\label{res:retry:mainthm}
	Let $\alpha \in (0,\infty)$ and let $\rho \in \mbn$.
	\begin{itemize}[itemsep = 0pt, topsep = \smallskipamount, label = \ensuremath{\bcdot}]
		\item 
		\emph{Fast Mixing.}
		Suppose $\alpha < \alpha_c(\rho)$ or $\alpha > 1$.
		Then there exists a constant $C$ so that,
		for all $K$ sufficiently large,
		all $\eps \in (0,1)$
		and all $n$ sufficiently large,
		we have
		$\tmix(\eps) \le C \log n$.
		
		\item 
		\emph{Slow Mixing.}
		Suppose $\alpha_c(\rho) < \alpha < 1$.
		Then there exists a constant $c$ so that,
		for all $K$ sufficiently large,
		all $\eps \in (0,\tfrac12)$
		and all $n$ sufficiently large,
		we have
		$\tmix(\eps) \ge e^{cn}$.
	\end{itemize}
	Further the map $\rho \mapsto \alpha_c(\rho) : \mbn \to (0,1)$ is strictly decreasing.
\end{thm}

\begin{rmkt*}
	Algebraic manipulations give
	$\alpha_c(1) = \tfrac13(5\sqrt{10} - 13) \approx 0.9732$,
	but for $\rho \ge 2$ we cannot solve symbolically. 
	Numerical calculations gives
		$\alpha_c(2) \approx 0.8662$,
		$\alpha_c(3) \approx 0.8191$
		and
		$\alpha_c(4) \approx 0.7858$.
	
	What we can see, however, is that $\rho \mapsto \alpha_c(\rho) : \mbn \to (0,1)$ is decreasing; this follows easily from the fact that, for each $\alpha$ and $f$, the map $\rho \mapsto r_\rho(f) : \mbn \to \mbr$ is increasing.
	This says that the interim slow-mixing region \emph{grows} with $\rho$.
	This should not be surprising:
		it is easier to accept a rerouted call when $\rho$ is larger;
		thus we do not need such a large blocking $f$ to obtain $\beta_\rho(f) > 1$.
\end{rmkt*}
	
%

\subsubsection{Motivation and Preliminary Properties}
\label{sec:retries:intro:motiv-prelim}

We now consider approximations to the high-level model, similar to those made before.
Consider the approximation for $\rho = 1$ which is exchangeable and a Markov process on the number of calls in each link (it does not differentiate between direct and rerouted calls).
In this model instead of picking a third station ($\gamma$) and attempting to route via this station (ie using $\alpha\gamma$ and $\gamma\beta$), we choose a pair of links uniformly at random: if there is free capacity on both then the call is accepted on these links (and the parts on the two links are released independently).
For general $\rho \in \mbn$, make the same approximation, choosing pairs (up to) $\rho$ times.

We now determine the rate, in the exchangeable model, at which calls arrive indirectly (ie via rerouting), to a specific route.
It only depends on the state of the system via the proportion $f$ of full links; we denote it $\lambda \cdot \tilde r_\rho(f)$.
The total arrival rate (ie direct and indirect) to a specific route is then $\lambda \rbr{ 1 + \tilde r_\rho(f) }$.
We show that $\tilde r_\rho(f) \to r_\rho(f)$ as \ninf (for $f$ independent of $n$).
The rates, in the limit \ninf, are thus the same in our model as in the exchangeable model.

The rate at which reroutings are attempted is $\lambda f n$, since there are $fn$ links that are full.
We observe that the probability that a reroute call is accepted (in one of the $\rho$ tries) is
\[
	1 - \rbb{ 1 - {\binomt{(1-f)n}2 / \binomt{n}2} }^\rho
\to
	1 - \rbb{ 1 - (1-f)^2 }^\rho.
\]
By symmetry (using the exchangeability), if a call is successfully rerouted (with some number of tries) then the pair chosen is uniform amongst all pairs of non-full links. Given a specific non-full link, the number of such pairs including this specific link is $(1-f)n-1$. Hence we see that
\[
	\tilde r_\rho(f)
&
=
	fn \cdot \rbb{ 1 - \rbb{ 1 - {\binomt{(1-f)n}2 / \binomt{n}2} }^\rho } \cdot \rbb{(1-f)n-1}/{\binomt{(1-f)n}2 }
\\&
\to
	2f \rbb{ 1 - \rbb{1 - (1-f)^2}^\rho } / (1-f)
=
	r_\rho(f).
\]
As a sanity check, observe that when $\rho = 1$ we do indeed get the same expression as in \S\ref{sec:motivation-comparison}.

Using these calculations, our additional approximation (ie going from the exchangeable model to our model) can be justified in the same way as in \S\ref{sec:motivation-comparison}.

\medskip

We now consider some properties of the polynomial $r_\rho$ that we are going to need.
The following claims, namely
\cref{res:retry:unique-maxima,res:retry:fsp>1/2,res:retry:num-solns},
are proved in \cref{app:poly-calcs}.
\begin{clm}
\label{res:retry:unique-maxima}
	There exists a unique stationary point $f^\SP_\rho \in [0,1]$ with $r_\rho'(f^\SP_\rho) = 0$.
\end{clm}

For example, algebraic manipulations (with a linear polynomial) give $f^\SP_1 = \tfrac12$, and solving a cubic numerically gives $f^\SP_2 \approx 0.5600$.
The next claim is used only in \S\ref{sec:retries:high}.
\begin{clm}
\label{res:retry:fsp>1/2}
	We have
	\(
		f^\SP_\rho \ge \tfrac12 \QUAD{for all} \rho \in \mbn.
	\)
\end{clm}

Similarly to before, the condition $p_\rho(f) > f$ is equivalent to the $h_\rho(f) > 0$ where
\[
	h_\rho(f)
\cq
	f \rbb{1 - 2 \rbb{1 - (1-f)^2}^\rho} + 1 - 1/\alpha.
\]
Recall that there are no solutions to $h_\rho(f) = 0$ for $\alpha < \alpha_c(\rho)$, by definition of $\alpha_c$.
\begin{clm}
\label{res:retry:num-solns}
	The following hold:
	\begin{itemize}[noitemsep, topsep = \smallskipamount]
		\item 
		for $\alpha \in (\alpha_c(\rho),1)$, there are precisely two (distinct) solutions $f \in (0,1)$ to $h_\rho(f) = 0$;
		
		\item 
		for $\alpha > 1$, there is a unique solution $f \in (0,1)$ to $h_\rho(f) = 0$.
	\end{itemize}
\end{clm}

\smallskip

Throughout this section
terms like $f^\SP_\rho$, $\alpha_c(\rho)$ or the zeros of $h_\rho$ will depend on $\rho$, but for notational ease we may sometimes drop the $\rho$ from the notation, unless it is explicitly needed.

\subsection{Coupling}
\label{sec:retries:intro:coup}

In this part
we give the coupling that we use for the fast mixing cases.
As previously, it will be `natural'; it is, in essence, the same as in the original ($\rho = 1$) case.

First, we give the analogous version of \cref{def:realisation} for the realisation of a single system.

\begin{defn}
\label{def:retry:realisation}
For arrivals, to each link give a Poisson stream (of arriving calls) of rate $2 \lambda$. Upon a call's arrival to a link, $k$ say, we have the following procedure.
\begin{itemize}
	\item 
	If the link $k$ is not full, then toss a $\Bern(\tfrac12)$-coin:
	\begin{itemize}[noitemsep, topsep=0pt]
		\item if heads (ie `1'), then add a call to link $k$;
		\item if tails (ie `0'), then do nothing.
	\end{itemize}
	
	\item
	If the link $k$ is full, then set $R = 1$ and run the following algorithm:
	\begin{ralist}[noitemsep, topsep = 0pt]
		\item 
		choose two links $i$ and $j$ uar (with replacement);
		
		\item 
		if both links $i$ and $j$ are not full, then add a call to link $i$ and \Stop;
		
		\item 
		if $R = \rho$, then \Stop; otherwise, increment $R \to R+1$ and \Return to \emph{Step (i)}.
	\end{ralist}
\end{itemize}

For departures, give each call in the system an independent exponential-$1$ clock. Upon a clock's ringing, remove the corresponding call from the system.
\end{defn}

Note that when $\rho = 1$ this reduces to the realisation of \cref{def:realisation}.
We also note that the algorithm described above is equivalent to the following one.
\begin{quote}
	If the link is full, then choose $i_1, ..., i_\rho$ and $j_1, ..., j_\rho$ uar
	and run the following algorithm:
	\begin{quotation}
		\noindent
		\FOR $R$ in $(1, ..., \rho)$ sequentially
		
		\noindent\quad
		\IF both links $i_R$ and $j_R$ are not full
		
		\noindent\qquad
		add a call to link $i_R$ and \Stop the loop
		
		\noindent\quad
		\END
		
		\noindent
		\END
	\end{quotation}
\end{quote}
While the algorithm described in \cref{def:retry:realisation} is perhaps easier to digest, it is the second algorithm that will generalise more easily to two (or more) systems.

We now give the analogous version of \cref{def:coupling} for the coupling of two systems.

\begin{defn}
\label{def:retry:coupling}
For arrivals, give to each link a Poisson stream (of arriving calls) of rate $2 \lambda$.
Upon a call's arrival to a link, $k$ say, we have the following procedure.
\begin{itemize}
	\item 
	Suppose $k$ is not full in either of $X$ or $Y$.
	Toss a $\Bern(\tfrac12)$-coin:
	\begin{itemize}[noitemsep, topsep = 0pt, label = \ensuremath{*}]
		\item if heads, then add a call to link $k$ both in $X$ and in $Y$.
	\end{itemize}
	
	\item 
	Suppose $k$ is full in $X$ but not full in $Y$.
	Toss a $\Bern(\tfrac12)$-coin:
	\begin{itemize}[noitemsep, topsep = 0pt, label = \ensuremath{*}]
		\item if heads, then add a call to link $k$ in $Y$.
	\end{itemize}
	Also, independently, choose $i_1, ..., i_\rho$ and $j_1, ..., j_\rho$ uar and run the following algorithm:
	\begin{quotation}
		\noindent
		\FOR $R$ in $1, ..., \rho$,
		\quad
		\IF both links $i_R$ and $j_R$ are not full in $X$,
		
		\noindent\qquad
		add a call to link $i_R$ in $X$ and \Stop the loop
		
		\noindent
		\END\quad \END
	\end{quotation}
	
	\item 
	Suppose $k$ is full in $Y$ but not full in $X$.
	Do analogously to the previous case.
	
	\item 
	Suppose $k$ is full both in $X$ and in $Y$.
	Choose $i_1, ..., i_\rho$ and $j_1, ..., j_\rho$ uar and run the following algorithms (once for $X$ and once for $Y$):
	\begin{quotation}
		\noindent
		\FOR $R$ in $1, ..., \rho$,
		\quad
		\IF both links $i_R$ and $j_R$ are not full in $X$ (respectively in $Y$),
		
		\noindent\qquad
		add a call to link $i_R$ in $X$ (respectively in $Y$) and \Stop the loop.
		
		\noindent
		\END\quad \END
	\end{quotation}
%
\end{itemize}

For departures, use the same rate-1 departure clocks in $X$ as in $Y$ where possible, giving the `extra' calls (ie those in $X$ but not in $Y$ or vice versa) independent rate-1 departure clocks.
\end{defn}

\begin{rmkt*}
	By inspection, one can see that this is a genuine, Markovian coupling.
	When using this coupling and $(X^0,Y^0) = (x,y)$, we denote it $\mbp_{x,y}$.
	Write $\mbp \cq \max_{(x,y) \in \Omega^2} \mbp_{x,y}$.
	Furthermore, it is a \textit{coalescent} coupling:
		we have $\bra{ X^t \ne Y^t } = \bra{ \tau_c > t }$,
		recalling that $\tau_c = \inf\bra{ t \ge 0 \mid X^t = Y^t }$.
\end{rmkt*}

This is a `natural' coupling, and can be applied whatever the state of the pair $(X,Y)$.
However, for the fast arrivals case $\alpha > 1$, we actually require a slightly more refined coupling, which will only work when the set of full links in $Y$ is a subset of those in $X$, or vice versa.
This is markedly different to the $\rho = 1$ case---although, in some sense, this definition will be an extension of the $\rho = 1$ case.
That coupling definition is deferred until it is required; it is given in \cref{def:retry:refined}.

%
%
%
%

\subsection{Slow Mixing in Interim Regime: $\alpha_c < \alpha < 1$}
\label{sec:retries:slow}

In this subsection
we consider the interim regime, $\alpha \in (\alpha_c(\rho), 1)$; we show slow mixing.
The statement is the natural extension of \cref{res:slow}; we sketch the argument, giving references to the $\rho = 1$ case, given in \S\ref{sec:slow}.

\begin{thm}
\label{res:retries:slow}
	For all $\rho \in \mbn$ and all $\alpha \in (\alpha_c(\rho), 1)$,
	there exists a constant positive $c$ so that,
	for all $K$ and $n$ sufficiently large,
	for all $t \le e^{cn}$,
	we have
	\[
		\MAX{x \in \Omega} \,
		\tvb{ \prb[x]{ X^t \in \cdot } - \Pi }
	\ge
		\tfrac12 - e^{-cn}.
	\]
	Thus
	for all $K$ sufficiently large,
	all $\eps \in (0,\tfrac12)$
	and all $n$ sufficiently large,
	we have
	$\tmix(\eps) \ge e^{cn}$.
\end{thm}

\begin{Proof}[Sketch of Proof]
The way we proved the slow-mixing case was to show that it takes exponentially long to move from a stable high-blocking state to a stable low-blocking state, and vice versa; call these `going down' and `going up', respectively.
We verify that these statements holds for general $\rho$.
We then deduce \cref{res:retries:slow} exactly as \cref{res:slow} was deduced from \cref{res:slow:high-to-low,res:slow:low-to-high}.

\smallskip

The proof that `going down' takes exponentially long required a general result on the rate $r$. We both upper bounded the system by $\Er(\beta(f^\SP_\rho),K)^n$ and lower bounded by $\Er(\beta(\delta),K)^n$, where $\delta$ is such that $p_\rho(\delta) > \delta$ and $\delta < f^\SP_\rho$.
We also used that $f \mapsto r_\rho(f)$ is strictly increasing on $[0,f^\SP_\rho)$ and strictly decreasing on $(f^\SP_\rho,1]$.
We proved this directly for $\rho = 1$; for general $\rho$, it follows from \cref{res:retry:unique-maxima}.
The proof then follow as previously, ie as in \cref{res:slow:high-to-low}.


The proof that `going up' takes exponentially long relied only on the fact that $f \mapsto r_\rho(f)$ is increasing on some interval of $[0,1]$ which includes 0.
While this is implied by \cref{res:retry:unique-maxima}, that claim is much stronger:
	it discusses the global behaviour of $r_\rho$;
	here we only need local behaviour near 0.
This local claim follows immediately from the fact that
	$r_\rho(0) = 0$ and $r_\rho(f) > 0$ for all $f \in (0,1)$
as well as the fact that
	$r_\rho$ is a polynomial (so has finitely many turning points).
The proof then follows as previously, ie as in \cref{res:slow:low-to-high}.
\end{Proof}

\subsection{Fast Mixing with Slow Arrivals: $\alpha < \alpha_c$}
\label{sec:retries:fast:low}

In this subsection
we consider the slow arrivals regime, $\alpha < \alpha_c(\rho)$; we show fast mixing.
The statement is the natural extension of \cref{res:low}; we sketch the argument, giving references to the $\rho = 1$ case, given in \S\ref{sec:low}.

\begin{thm}
	For all $\rho \in \mbn$ and all $\alpha < \alpha_c(\rho)$,
	there exists a constant $C$ so that,
	for all $K$ and $n$ sufficiently large,
	if $X,Y \sim \DAR_n^\rho(\alpha,K)$, then under the coupling \mbpcoup,
	for all $t \ge C \log n$,
	we have
	\[
		\MAX{(x,y) \in \Omega^2}
		\prb[x,y]{ X^t \neq Y^t }
	\le
		C / n
	=
		\oh1.
	\]
	Thus
	for all $K$ sufficiently large,
	all $\eps \in (0,1)$
	and all $n$ sufficiently large,
	we have
	$\tmix(\eps) \le C \log n$.
\end{thm}

The adaptation from $\rho = 1$ to general $\rho \in \mbn$ is straightforward here.

\begin{Proof}[Sketch of Proof]
For the burn-in phase, all that we used was that $p(f) < f$ for all $f \in [0,1]$. This result still holds here (by definition of $\alpha_c$); hence the burn-in phase proof is identical.

For the variable length coupling, we use the same stopping time. Observe that if the blocking is at $f \le 2 \eps$ with $\eps$ sufficiently small, then
\[
	r_\rho(f) \le 4 \eps / (1 - 2 \eps) \le 5 \eps.
\]
(This does not require any `unique local maximum' property, or anything like this; it simply uses the fact that $1 - (1 - (1 - f)^2)^\rho \in (0,1)$ for $f \in (0,1)$.)
Hence the same proof for the coupling works also, up to changing some constants.
\end{Proof}

\subsection{Fast Mixing with Fast Arrivals: $\alpha > 1$}
\label{sec:retries:high}

In this subsection
we consider the fast arrivals regime, $\alpha > 1$; we show fast mixing.
The statement is the natural extension of \cref{res:high}; we sketch the argument, giving references to the $\rho = 1$ case, given in \S\ref{sec:high}.

\begin{thm}
	For all $\rho \in \mbn$ and all $\alpha > 1$,
	there exists a constant $C$ so that,
	for all $K$ and $n$ sufficiently large,
	if $X,Y \sim \DAR_n^\rho(\alpha,K)$, then,
	under the coupling given by \mbpcoup,
	for all $t \ge C \log n$,
	we have
	\[
		\MAX{(x,y) \in \Omega^2}
		\prb[x,y]{ X^t \neq Y^t }
	\le
		C/n.
	\]
	Thus
	for all $K$ sufficiently large,
	all $\eps \in (0,1)$
	and all $n$ sufficiently large,
	we have
	$\tmix(\eps) \le C \log n$.
\end{thm}

This regime is rather harder to prove in the general-$\rho$ case; in particular, it requires more detailed knowledge of the high-degree polynomials in question, but also we need to introduce a new coupling.
While some of the argument below will be sketched, similarly to in the previous two proofs, anything new will be explained fully and rigorously.

\smallskip

First, we describe the new coupling needed; see \cref{def:retry:refined} below.
Using the $\rho$-retries coupling of \cref{def:retry:coupling},
if a rerouting happens in two systems $X$ and $Y$ from the same link, then the it can be accepted in both $X$ and $Y$, but onto different links.
This was not possible using \cref{def:coupling} when $\rho = 1$, since there was only one rerouting attempt and the same links were chosen in $X$ as in $Y$.
So we see that the relative distance can actually increase by 2.
This behaviour makes controlling the `difference' between the two systems difficult; it is highly undesirable.

We now give a more refined coupling, fixing this issue.
Write $f_X$ for the proportion blocked in $X$ and $f_Y$ for the proportion blocked in $Y$.

\begin{defn}
\label{def:retry:refined}
Let $X, Y \sim \DAR_n^\rho(\alpha,K)$, with the pair $(X,Y)$ in a state with $d(X,Y) = 1$.

If the mismatched link is not full in either $X$ or $Y$, then use the (original) description given in \cref{def:retry:coupling}.
Suppose then, without loss of generality, that $X_1 = K > Y_1$.

For arrivals, give to each link a Poisson steam (of arriving calls) of rate $2 \lambda$.
Upon a call's arrival to a link, $k$ say, we have the following procedure.
\begin{itemize}
	\item 
	Suppose $k$ is not full in either of $X$ or $Y$.
	Toss a $\Bern(\tfrac12)$-coin:
	\begin{itemize}[noitemsep, topsep=0pt]
		\item 
		if heads, then add a call to link $k$ both in $X$ and in $Y$.
	\end{itemize}
	
	\item 
	Suppose $k$ is full in $X$ but not full in $Y$ (ie $k = 1$).
	Toss a $\Bern(\tfrac12)$-coin:
	\begin{itemize}[noitemsep, topsep=0pt]
		\item 
		if heads, then add a call to link $k$ in $Y$.
	\end{itemize}
	Also, independently, choose pairs $(i_1,j_1), (i_2,j_2), ... \in [n]^2$ uar (with replacement) until a pair $(i,j)$ has neither $i$ nor $j$ full in $Y$ or until $\rho$ have been chosen; call the final pair $(i,j)$.
	Add a call to link $i$ in $Y$ if neither $i$ nor $j$ are full in $Y$.
	
	\item 
	Suppose $k$ is full both in $X$ and in $Y$.
	Choose pairs $(i_1,j_1), (i_2,j_2), ... \in [n]^2$ uar (with replacement) until a pair $(i,j)$ has neither $i$ nor $j$ full in $Y$ or until $\rho$ have been chosen; call the chosen pairs $(i_1,j_1), ..., (i_R,j_R)$.
	If the above selection is terminated by choosing $\rho$ `unsuitable' pairs, ie pairs with at least one full in $Y$,
	then do nothing.
	Otherwise,
	perform the following procedure.
	\begin{itemize}[noitemsep, topsep=0pt]
		\item 
		Add a call to link $i_R$ in $Y$.
		\item 
		If $1 \notin \bra{i_R,j_R}$,
		add a call to link $i_R$ in $X$.
		
		\item 
		Now suppose that $1 \in \bra{i_R,j_R}$,
		and perform the following procedure:
		
		\begin{itemize}[noitemsep, topsep = 0pt, label = \bcdot]
			\item 
			with probability $1 - (1 - (1 - f_X)^2)^{\rho-R}$,
			if $1 \neq i_R$,
			then add a call to link $i_R$ in $X$,
		and
			otherwise (ie if $1 = i_R$),
			independently choose $i$ uar from the set of non-full links in $X$ and add a call to link $i$ in $X$;
			
			\item 
			with probability $(1 - (1 - f_X)^2)^{\rho-R}$,
			do nothing
			(ie do not add a call to $X$).
		\end{itemize}
	\end{itemize}
\end{itemize}

Couple departures, with the same rate-1 clocks in $X$ and in $Y$, as before.
\end{defn}

Recall that in the previous coupling it was possible to add to two different links (one in $X$ and one in $Y$). The above coupling mitigates this issue:
	now it can be the case that link $1$ is mismatched, with $X_1 = K > Y_1$, and a reroute pick link $1$ in $Y$ and a different link in $X$;
	hence the mismatched link can change, but in a way that keeps the relative distance 1.
(Note that if a reroute picks a link other than 1 in $Y$, then the same link is picked in $X$.

\smallskip

Our stopping time will be such that if $d(X^0,Y^0) = 1$ then $d(X^s,Y^s) = 1$ for all $s < \tau$, and so we shall be able to use the explicit formulation of the coupling given above.

\begin{rmkt*}
	It is not difficult to check, and we do so below, that this is a genuine, Markovian coupling.
	When using this coupling and $(X^0,Y^0) = (x,y)$, we denote it $\mbp_{x,y}$.
	Write $\mbp \cq \max_{(x,y) \in \Omega^2} \mbp_{x,y}$.
	Furthermore, it is a \textit{coalescent} coupling:
		we have $\bra{ X^t \ne Y^t } = \bra{ \tau_c > t }$.
\end{rmkt*}

\begin{rmkt*}
	This is, in some sense, an extension of the single-try coupling given in \cref{def:coupling}.
	There if we had $(i,j)$ with neither $i$ nor $j$ full in $Y$ but one full in $X$, then we did not have another try to pick another choice for $X$.
	We always had $R = \rho = 1$ in the notation above.
\end{rmkt*}

\begin{Proof}[Validity of Coupling]
The required independence structure of the links in each system is immediate and the departures are as required.
Consider arrivals with $(X,Y)$ in a state with $d(X,Y) = 1$.
Without loss of generality, assume that $X_1 = k > Y_1$.

For arrivals to a link that is not full in $Y$ (but may or may not be in $X$), the arrival rate is as required.
When the link is full both in $X$ and in $Y$, the reroute pair $(i_R,j_R)$ for $Y$ is chosen by the usual procedure; if $1 \notin \bra{i_R,j_R}$,
then the disparity between $X$ and $Y$ has not played a role and a call is added to $i_R$ in both $X$ and $Y$, as required.

Now suppose that $1 \in \bra{i_R,j_R}$; we then need further retries to add a call to $X$. If a call is added to $X$, then which link is chosen must be uniform amongst the non-full links in $X$.
As some notation, if we let $\mcf_Y$ be the set of full links in $Y$ and $\mcf_X$ in $X$, then we see that $\mcf_X = \mcf_Y \cup \bra{1}$ (where $X_1 = K > Y_1$); so $\mcf_Y^c = \mcf_X^c \cup \bra{1}$.
Note that $i_R$ is chosen uar from $\mcf_Y^c$. We wish to choose $i$ uar from $\mcf_X^c$.
Conditional on $i_R \in \mcf_X^c$, ie $1 \neq i_R$, we see that $i_R$ is uniform over $\mcf_X^c$; on this event we set $i = i_R$.
On the complementary event, ie $1 = i_R$, we independently choose $i$ uar from $\mcf_X^c$. Hence, if a call is added, then which link is chosen is uniform amongst the available links, as required.
It is immediate that the probability with which this step is taken is correct.
	%
\end{Proof}

We now sketch the proof of \cref{res:retry:mainthm}. We outline all the ideas, but omit some details.

\begin{Proof}[Sketch of Proof]
Define $\varphi_\rho$ so that the probability a reroute is accepted is precisely $\tfrac12$:
\[
	1 - \rbb{ 1 - (1 - \varphi_\rho)^2 }^\rho = \tfrac12,
\Quad{and so}
	\varphi_\rho \cq 1 - \sqrt{1 - 1/2^\rho}.
\]
(Previously we had $\varphi_c = \varphi_1 = 1 - 1/\sqrt2$.)
Note that $\varphi_\rho$ satisfies
\[
	h_\rho(\varphi_\rho)
=
	\varphi_\rho \cdot \rbb{ 1 - 2 \cdot \tfrac12 } + 1 - 1/\alpha = 1 - 1/\alpha,
\]
and hence $h_\rho(\varphi_\rho) > 0$ if and only if $\alpha > 1$.
Note that $\varphi_\rho$ is decreasing in $\rho$; thus $\varphi_\rho \le \varphi_1 = 1 - 1/\sqrt 2$ for all $\rho$.
By \cref{res:retry:fsp>1/2}, we have $f^\SP_\rho \ge \tfrac12$ for all $\rho$.
We thus deduce that $\varphi_\rho < f^\SP_\rho$ for all $\rho \in \mbn$.

Also note that
\(
	p(0) > \varphi_\rho
\)
if and only if
\(
	1 - 1/\alpha > 1 - \sqrt{1 - 2^{-\rho}}
\),
which in turn holds if and only if
\(
	\alpha > (1 - 2^{-\rho})^{-1/2},
\)
and note that this lower bound is decreasing in $\rho$ (and decreases down to 1);
previously (when $\rho = 1$), we had $\alpha > \sqrt2$.
We split into two subregimes:\par
{\centering \smallskip
	$\alpha \le (1 - 2^{-\rho})^{-1/2}$
and
	$\alpha > (1 - 2^{-\rho})^{-1/2}$;\par \smallskip
}\noindent%
this is analogous with the $\rho = 1$ case where, for $\alpha > 1$, we separated $\alpha \le \sqrt2$ and $\alpha > \sqrt2$.

As in the $\rho = 1$ case, for the regime with $\alpha > (1 - 2^{-\rho})^{-1/2}$, we stochastically dominate our system from below by an $\Er(\alpha,K)^n$ system, which has expected proportion of full links greater than $\varphi_\rho$, ie $p(0) > \varphi_\rho$.
We use a burn-in of length $10 \log n$ to reach a state where the proportion blocked is strictly greater than $\varphi_\rho$.

For $1 < \alpha \le (1 - 2^{-\rho})^{-1/2}$, we again use an algorithmic stochastic domination procedure, analogous to that used for $\alpha \le \sqrt2$ in \S\ref{sec:high:burnin}.
Since $\alpha > 1$, by \cref{res:retry:num-solns} the polynomial $h_\rho(\cdot)$ has precisely one zero, and also $h_\rho(0) > 0$ and $h_\rho(\varphi_\rho) > 0$;
hence there exists an $\xi > 0$ so that $h_\rho(f) > 0$, ie $p_\rho(f) > f$, for all $f \in [0,\varphi_\rho + 2\xi]$.
Using the same algorithmic procedure as before, we are able to get the proportion blocked to be strictly greater than $\varphi_\rho$.
We need to be slightly careful, though. Recall that, in the $\rho = 1$ case, we first dominated above by $\Er(\beta(f^\SP_1),K)^n$, and said that this meant the proportion blocked was at most $\tfrac23$.
However, we could have been more restrictive and used anything larger than
\[
	f_\text{up}
\cq
	p(f^\SP_\rho)
=
	1 - \alpha^{-1} \rbb{1 + r_\rho(f^\SP_\rho)}^{-1};
\Quad{when $\rho = 1$,}
	f_\text{up} = 1 - \tfrac23/\sqrt2 \approx 0.5286.
\]
The key is that
$r_\rho(f_\text{up}) > r_\rho(f)$ for all $f \in [0,\varphi_\rho]$,
as we justify now:
some crude bounds give
\[
	r_\rho(\varphi_\rho) \le 2^{-\rho}
\Quad{and}
	r_\rho(f_\text{up}) > \tfrac13 > 2^{-\rho} \text{ for } \rho \ge 2;
\]
the case $\rho = 1$ was the original case, and in it we had
\[
	0.4142 \approx \sqrt2 - 1 = r_\rho(\varphi_\rho) < r_\rho(f_\text{up}) \approx 0.4984.
\]
Hence we may apply the stochastic domination procedure to get strictly greater than $\varphi_\rho$ blocking.


As in \S\ref{sec:high:burnin}, the burn-in phase will be of length $C \log n$ with $C$ independent of $K$ and $n$. After this phase, with high probability the system will be in a state that has blocking proportion strictly greater than $\varphi_\rho$ and will keep this property for a long while (say at least for time $n$).

Observe that the \textit{acceptance probability} $1 - (1 - (1 - f)^2)^\rho$ is strictly decreasing in $f$. Hence if the blocking is at level greater than $\varphi_\rho$, then the acceptance probability is strictly less than $\tfrac12$.
When the reroute is declined, the systems coalesce; when it is accepted, the relative distance increases from 1 to 2. Hence the desire to have the acceptance probability strictly less than $\tfrac12$ so that the expected relative distance decreases.

\smallskip

As noted above, which link is mismatched can change before the relative distance changes.
However, we can still use exactly the same stopping time as previously, as given in \cref{def:high:st}, and the same bounds apply; we just need to change the definition to allow for which link is mismatched to change.
Hence we make the same deductions, completing the proof as before.
\end{Proof}

\section{Dynamic Alternative Routing with Trunk Reservation}
\label{sec:tr}

\subsection{Introduction}

We add \textit{trunk reservation}.
We reserve some amount $\sigma$ of the capacity of a link for direct arrivals only:
	if a call tries to be rerouted, then it can only be held on links with current load strictly less than $K-\sigma$, ie strictly more than $\sigma$ free circuits.
	If we set $\sigma = 0$, then we recover our previous model.
In a similar vein to the previous sections, we consider an approximation to this model via a type of state-dependent Poisson arrival process with independent departures.

One could then allow multiple rerouting attempts, as in the previous section.
For clarity of exposition, we consider only one rerouting attempt.

We show that when $\sigma$ is chosen appropriately, we \emph{do not} get metastability:
	that is,
	we do not have an intermediary slow-mixing regime, for $\alpha \in (\alpha_c,1)$, between a high-blocking and low-blocking regime;
	rather we have fast mixing for all $\alpha < 1$.
Importantly, we can choose $\sigma$ to be some fixed number, depending only on $\alpha$, but not on $K$ or $n$.
This is very important from a network engineering point of view:
	the system is \emph{scalable} in the sense that as $\ninf$ and $\Kinf$,
	the number of circuits to be reserved does not grow.

\smallskip

\begin{Proof}[Intuition]
\renewcommand{\qedsymbol}{\ensuremath{\triangle}}
We now make a few comments on why we do not have metastability (for $\sigma$ sufficiently large) for the regime $\alpha \in (0,1)$.
To be in a `stable' high-blocking regime (ie one where the proportion of full links remains bounded away from 0 for a long time), trunk reservation requires the system not only to `support' the current number of full links, but not have too many `nearly full' (ie capacity between $K-\sigma$ and $K$) links, as `nearly full' links do not accept reroutings.
We show that these will be (approximately) mutually exclusive events:
	the invariant distribution of a single supercritical Erlang link concentrates very tightly at the capacity, so if we look at $\bra{K-\sigma, ..., K}$ then this will contain almost all of the invariant mass when $\sigma$ is sufficiently large;
	so a large number of `nearly full' links implies a large number of full links and vice versa.

Our analysis of the original low-blocking regime, ie $\alpha < \alpha_c$, relied on the fact that $p(f) < f$ for all $f \in [0,1]$. This allowed us to repeatedly stochastically dominate from above until we reached a suitably low-blocking set.
We use this same idea with trunk reservation when $\alpha < 1$.
	%
\end{Proof}

\subsubsection{Model Set-Up and Main Theorem}
\label{sec:tr:model}

Similarly to in previous sections, we first describe explicitly our model and then, after, explain why this is the right model.
Fix $\sigma \in \mbn$, independent of $K$ and $n$. Given $x \in \Omega$, write
\[
	\mcf \cq \brb{ j \in [n] \midb x_j = K }
\Quad{and}
	\mcg \cq \brb{ j \in [n] \midb x_j \ge K-\sigma };
\]
also write $f \cq \abs \mcf/n$ and $g \cq \abs \mcg/n$, and note that $\mcf \subseteq \mcg$ so $f \le g$.
As before, $f$ is the proportion of full links; the new variable, $g$, is the proportion of `almost full or full links', ie ones that will not accept a rerouting.
As before, we write $\lambda \cq \alpha K$.

Suppose the system is in state $x$.
To the links in $\mcf$,
no calls arrive.
To the links in $\mcg \setminus \mcf$,
calls arrive at rate $\lambda$.
To the links in $\mcg^c$,
calls arrive at rate
\(
	\lambda \rbr{ 1 + 2 f (1-g) }.
\)
Calls depart at rate-1 independently.
We denote this system $\DAR_n^\sigma(\alpha,K)$.

This is now a `mixture of state-dependent Poisson arrival processes':
	one for the links in $\mcg \setminus \mcf$ of rate $\lambda$
and
	one for the links in $\mcg^c$ of rate $\lambda \rbr{ 1 + 2 f (1-g) }$;
	the links in $\mcf$ receive no arrivals.

\medskip

We require $\sigma \ge \sigma_*(\alpha)$, for some constant $\sigma_*$ that depends only on $\alpha$, not on $K$ or $n$.

\begin{defn}
\label{def:tr:sigma*}
	For $\alpha \in (0,\infty) \setminus \bra{1}$,
	define $\sigma_*(\alpha)$ as follows:
	\[
		\sigma_*(\alpha)
	\cq
	\begin{cases}
		\tfrac1{1-\alpha} \rbb{ 2 \log\rbb{ \tfrac1{1-\alpha} } + 14 }
			&\text{when}\quad
		\alpha < 1;
	\\
		\log 4 / \log \alpha
			&\text{when}\quad
		\alpha > 1.
	\end{cases}
	\]
	We reserve the notation $\sigma_*$ for this parameter.
\end{defn}


When considering $\alpha < 1$, ie slow arrivals, we only really need to apply this for $\alpha \ge \alpha_c$, where $\alpha_c$ is the original critical threshold:
	for $\alpha < \alpha_c$, we already had fast mixing with the additional property that an insignificant proportion of links were blocked in equilibrium.
(For $\alpha > 1$, calls arrive faster than they could possibly be processed, regardless of the state of the system, so it is always the case that a significant proportion of calls are lost.)
The idea, from an application point of view, is that a single system should be used always, and be able to handle busy periods without metastability; for this reason, we analyse $\alpha < \alpha_c$ and $\alpha > 1$ as well as $\alpha_c < \alpha < 1$.

The main theorem of this section is the following mixing time result.
(It was stated as \cref{res:intro:tr} in the introduction; we recall it here for convenience.)
There will be two regimes: $\alpha < 1$ and $\alpha > 1$; the proofs will be in a similar vein to the previous fast mixing proofs.

\begin{thm}
	For all $\alpha \in (0,\infty) \setminus \bra{1}$,
	if $\sigma \ge \sigma_*(\alpha)$ (as given in \cref{def:tr:sigma*}), then
	for
		all $\eps \in (0,1)$
	and
		all $K$ and $n$ sufficiently large,
	we have
	$\tmix(\eps) \le 60 \log n$.
	
	Furthermore, with these parameters, when $\alpha < 1$, the proportion of links which are full, in equilibrium, may be made as small as desired by taking $K$ sufficiently large (independently of $n$).
\end{thm}

\begin{rmkt*}
	No effort has been made to optimise the particular choice of $\sigma_*$.
	The key is not the particular value, but that it can be chosen independently of $K$ and $n$, depending only on $\alpha$.
\end{rmkt*}

\begin{Proof}[Proof References]
\renewcommand{\qedsymbol}{\ensuremath{\triangle}}
	See \cref{res:tr:low,res:tr:high} for the cases $\alpha < 1$ and $\alpha > 1$, respectively.
\end{Proof}


%
%

\subsubsection{Motivation, Notation and Interpretation}

Consider the exchangeable model with trunk reservation:
	calls arrive to each link as a Poisson process with rate $\lambda$;
	if the link is full on arrival, then
		a rerouting is attempted,
		two links are chosen uniformly at random
	and
		the call is accepted (and held on both links) if and only if both have more than $\sigma$ spare capacity.
We determine the rate at which calls are accepted to different links in this model.
We show that these rates are asymptotically equivalent to the rates in the model we described in the previous part.
First note that
	if a link is full then no calls are accepted to it
and
	if it is not full but has at most $\sigma$ spare capacity then calls arrive as a Poisson process of rate $\lambda$;
this is the same in our model.

Now consider links with fewer than $K - \sigma$ calls on them, ie links in $\mcg^c$.
Suppose the system is in state $x$ with proportions $f$ and $g$ defined previously.
In the same way as before, attempted reroutings happen at rate $\lambda f n$; these are accepted if both links chosen have current load less than $K-\sigma$, and hence the probability of being accepted is $\binom{(1-g)n}2 / \binom{n}2 \approx (1-g)^2$.
By symmetry, if a call is successfully rerouted then the pair chosen is uniform amongst all pairs of links with fewer than $K-\sigma$ calls on them. Given a specific such link, the number of such pairs including this specific link is $(1-g)n - 1$.
Hence we see that the arrival rate due to reroutings is
\[
	fn \cdot \rbB{ \binomt {(1-g)n}2 / \binomt n2 } \cdot \rbB{ \rbb{ (1-g)n - 1 } / \binomt {(1-g)n}2 }
\to
	2 f (1-g).
\]
As before, these calculations (go some way to) justify our approximation.

\medskip

For $f,g \in [0,1]$,
analogously to without trunk reservation, we set
\begin{gather*}
	\beta(f,g) \cq \alpha\rbb{1+2f(1-g)}
\Quad{and}
	p(f,g) \cq 1 - 1/\beta(f,g);
\\
	\text{also define}
\quad
	q(f,g) \cq 1 - 1/\beta(f,g)^{\sigma+1},
\end{gather*}
which (direct calculation shows) is the ($\Kinf$)-limit of the equilibrium probability that an $\Er(\beta(f,g),K)$ link has at least $K-\sigma$ calls on it.
(Recall that $\sigma$ is independent of $K$, so $q$ is the genuine limit; it is not an asymptotic statement or approximation.)

Observe that $\beta(f,g)$, $p(f,g)$ and $q(f,g)$ are increasing in $f$ and decreasing in $g$ when the other argument is held fixed. Considering the interpretation of $f$ and $g$, we are only interested in these when $f \le g$. In particular, subject to $f \le g$, all three functions are maximised when $f = g = \tfrac12$.

\subsection{Preliminaries}

\subsubsection{Coupling}

As previously, we use the `natural' realisation of a single system and extend this to a coupling of two systems.
Since there is no slow-mixing regime here, we use this coupling for all the regimes.

We first give the realisation of a single system; cf \cref{def:realisation}.

\begin{defn}
\label{def:tr:realisation}
For arrivals, to each link give a Poisson stream (of arriving calls) of rate $2 \lambda$. Upon a call's arrival to a link, $k$ say, we have the following procedure.

\begin{itemize}
	\item 
	If the link is not full, then toss a $\Bern(\tfrac12)$-coin:
	\begin{itemize}[noitemsep, topsep=0pt]
		\item if heads (ie `1'), then add a call to link $k$;
		\item if tails (ie `0'), then do nothing.
	\end{itemize}
	
	\item If the link is full, then choose two links $i$ and $j$ uar (with replacement):
	\begin{itemize}[noitemsep, topsep=0pt]
		\item if both links $i$ and $j$ have fewer than $K - \sigma$ calls on them, then add a call to link $i$;
		\item otherwise, ie if either link $i$ or link $j$ has at least $K - \sigma $ calls on it, do nothing.
	\end{itemize}
\end{itemize}

For departures, give each call in the system an independent exponential-$1$ clock. Upon a timer's ringing, remove the corresponding call from the system.
\end{defn}

This realisation naturally extends to a couple of two systems, $X$ and $Y$; cf \cref{def:coupling}.
Recall the definition of $\mcf$ and $\mcg$ as the set of full and full or almost full links, respectively.
Write $\mcf_X$ and $\mcg_X$ to denote the relevant quantities in the $X$-system and $\mcf_Y$ and $\mcg_Y$ for the $Y$-system.

\begin{defn}
\label{def:tr:coupling}
For arrivals, give to each link a Poisson stream (of arriving calls) of rate $2 \lambda$.
Upon a call's arrival to a link, $k$ say, we have the following procedure.
\begin{itemize}
	\item 
	Suppose $k$ is not full in either of $X$ or $Y$.
	Toss a $\Bern(\tfrac12)$-coin:\par
	\quad if heads, then add a call to link $k$ both in $X$ and in $Y$.
	
	\item 
	Suppose $k$ is full in $X$ but not full in $Y$.
	Toss a $\Bern(\tfrac12)$-coin:\par
	\quad if heads, then add a call to link $k$ in $Y$.\par
	Also, independently, choose two links $i$ and $j$ uar (with replacement):\par
	\quad if $i,j \in \mcg^c_X$, then add a call to link $i$ in system $X$.
	
	\item 
	Suppose $k$ is full in $Y$ but not in $X$.
	Do analogously to the previous case.
	
	\item 
	Suppose $k$ is full both in $X$ and in $Y$.
	Choose two links $i$ and $j$ uar (with replacement):\par
	\quad if $i,j \in \mcg^c_X$ (respectively $i,j \in \mcg_Y^c$), then add a call to link $i$ in system $X$ (respectively~$Y$).
\end{itemize}

For departures, use the same rate-1 departure clocks in $X$ as in $Y$ where possible, giving the `extra' calls (ie those in $X$ but not in $Y$ or vice versa) independent rate-1 departure clocks.
\end{defn}

\begin{rmkt*}
	By inspection, one can see that this is a genuine, Markovian coupling.
	When using this coupling and $(X^0,Y^0) = (x,y)$, we denote it $\mbp_{x,y}$.
	Write $\mbp \cq \max_{(x,y) \in \Omega^2} \mbp_{x,y}$.
	Furthermore, it is a \textit{coalescent} coupling:
		we have $\bra{ X^t \ne Y^t } = \bra{ \tau_c > t }$,
		recalling that $\tau_c = \inf\bra{ t \ge 0 \mid X^t = Y^t }$.
	
	This can be extended from `single try to `$\rho$ retries' similarly to how we did in \S\ref{sec:retries}, if desired.
\end{rmkt*}

%
%
%
%

\subsubsection{Erlang Link with Trunk Reservation}
\label{sec:tr:trunked-erlang}

We also need properties of the following Markov chain, which we call a `trunk-reserved Erlang link' and denote by $\ET(\alpha,\beta,\sigma,K)$; the parameters are explained below.
It has capacity $K$ with the following dynamics:
	calls arrive at rate $\beta K$ when the current load is less than $K-\sigma$,
	at rate $\alpha K$ when the load is at least $K-\sigma$ but the link is not full
and
	no calls arrive when the link is full;
	calls depart independently at rate 1.
Formally, it is a Markov chain on $\bra{0,...,K}$ with non-zero transition rates $(r_{i,j})_{i,j \in [K]}$ given by
\[
	r_{i,i-1} = i
\Quad{for}
	i = 1,...,K
\Quad{and}
	r_{i,i+1}
=
	\begin{cases}
		\alpha K	& \text{for}\quad i = 0, ..., K-\sigma-1,\\
		\beta K		& \text{for}\quad i = K-\sigma, ..., K-1.
	\end{cases}
\]
We often write $\lambda \cq \alpha K$.
We always consider $\alpha \le \beta$. 
We now investigate typical behaviour of the link in different regimes determined by $\alpha$ and~$\beta$.
Let $F \sim \ET(\alpha,\beta,\sigma,K)$.


\begin{Proof}[Consider $\alpha \le \beta < 1$]
\qedtriangle
We can stochastically dominate $F \lesssim E \sim \Er(\beta,K)$, which is a subcritical Erlang system.
Then $\pi_F(A) \le \pi_E(A)$ for any up-set $A$.
In particular, we have
\[
	\pi_F\rbb{ \bra{K-\sigma, ..., K} }
\le
	\pi_E\rbb{ \bra{K-\sigma, ..., K} },
\quad
	\text{which decays exponentially in $K$}.
\qedhere
\]
	%
\end{Proof}


\begin{Proof}[Consider $1 < \alpha \le \beta$]
\qedtriangle
We can stochastically dominate $F \gtrsim E \sim \Er(\alpha,K)$, which is a supercritical Erlang system.
Then $\pi_F(A) \le \pi_E(A)$ for any down-set $A$.
By direct calculation, we have
\[
	\pi_E\rbb{ \bra{K-\sigma, ..., K} } \to 1 - \alpha^{-(\sigma+1)}
\Quad{as}
	\Kinf,
\]
recalling that $\sigma$ depends only on $\alpha$ (not on $K$),
and hence if
\[
	\sigma \ge \log_{\tilde\alpha}(1/\eps) = \log(1/\eps)/\log \tilde\alpha
\]
for some $\tilde\alpha < \alpha$, eg $\tilde\alpha \cq \tfrac12(1 + \alpha)$, then $\pi_F(\bra{K-\sigma,...,K}) \ge 1 - \eps$ for large enough $K$.
(Observe that the complement of a down-set is an up-set.)
\end{Proof}


\begin{Proof}[Consider $\alpha < 1 \le \beta$]
\qedtriangle
We can stochastically dominate $F \lesssim L + (K-\sigma)$ where $L$ is a biased simple random walk on $\bra{0,...,\sigma}$ with up-rate $\lambda$ and down-rate $K-\sigma$; note that $\alpha = \lambda/K < 1$ implies $\alpha' \cq \lambda/(K-\sigma) < 1$ for large enough $K$ (as $\sigma$ is independent of $K$).
(This stochastic domination basically says that if $F$ is in $\bra{0,...,K-\sigma}$ then we assume it is at the highest point, ie $K-\sigma$; since $\beta \ge 1$, we have concentration of an Erlang link with arrival rate $\beta K$ and capacity $K-\sigma$ about its capacity, so this domination is not as wasteful as it may initially appear.)
Then the invariant distribution of $L$, which we denote $\pi_L$, decays exponentially, ie $\pi_L(j) \le \tilde\alpha^j$ for some $\tilde\alpha \in (0,1)$. In fact, it is easy to see that we may set $\tilde\alpha = \tfrac12(\alpha+1) < 1$ providing $K$ is large enough.
From this we obtain $\pi_F(K) \le \pi_L(\sigma) \le \tilde\alpha^\sigma$.
We want this to be at most $\eps$, so we require
\[
	\sigma \ge \log_{\tilde\alpha}(\eps) = \log(1/\eps)/\log(1/\tilde\alpha),
\]
with $\eps > 0$ to be specified later.
Note that this domination is independent of $\beta \ge 1$.
\end{Proof}



In \cref{res:prelim:product-mixing}, we stated that the mixing time of an Erlang system of $n$ links is at most $10 \log n$; this is proved in
\cref{res:er-mix:multi-link:mixing}.
In \cref{app:er-mix},
we actually consider the trunked links, which are a generalisation of the standard links. (Take $\alpha \cq \beta$ and $\sigma \cq 0$ to reduce $\ET(\alpha,\beta,\sigma,K)$ to $\Er(\beta,K)$.)
We prove the same mixing result for trunked links as for standard links; see
\cref{res:er-mix:multi-link:mixing}.


\begin{lem}
\label{res:tr:trunked-erlang:product-mixing}
	For
		all $\alpha, \beta \in (0,\infty)$ with $\alpha \le \beta$,
		all $\sigma, K \in \mbn_0$ with $K \ge \max\bra{\sigma,1}$
	and
		all $n$ sufficiently large,
	we have
	\[
		\tmix^{\alpha,\beta,\sigma,K;n}(1/n) \le 10 \log n.
	\]
\end{lem}

\subsection{Fast Mixing with Slow Arrivals: $\alpha < 1$ and Low-Blocking}
\label{sec:tr:low}

In this subsection
we consider the `slow arrivals' regime, ie $\alpha < 1$.
We always use the coupling \mbpcoup from \cref{def:tr:coupling}.
The aim is to prove the following theorem.
Recall~that
\[
	\sigma_*(\alpha)
=
	\rbr{ - 2 \log(1 - \alpha) + 14 } / (1 - \alpha)
\Quad{for}
	\alpha < 1.
\]

\begin{thm}
\label{res:tr:low}
	For
		all $\alpha < 1$,
		all $\sigma \ge \sigma_*(\alpha)$
	and
		all $K$ and $n$ sufficiently large,
	if $X, Y \sim \DAR_n^\sigma(\alpha,K)$,
	under the coupling \mbpcoup,
	for all $t \ge 60 \log n$,
	we have
	\[
		\MAX{(x,y) \in \Omega^2}
		\prb[x,y]{ X^t \neq Y^t }
	\le
		9/n
	=
		\oh1.
	\]
	Thus
	for all $K$ sufficiently large,
	all $\eps \in (0,1)$
	and all $n$ sufficiently large,
	we have
	$\tmix(\eps) \le 60 \log n$.
	
	Furthermore, with these parameters, the proportion of links which are full, in equilibrium, may be made as small as desired by taking $K$ sufficiently large (independently of $n$).
\end{thm}

The proof of this proposition will be very similar to that of \cref{res:low}, except that instead of using a careful sequence of stochastic dominations for the burn-in the trunk reservation will allow us to do only one stochastic domination.
Note that the result has a key difference beyond allowing $\alpha$ to get arbitrarily close to 1: previously, the upper bound on the mixing time depended on $\alpha$; now it does not.

\medskip

We prove the following burn-in proposition.

\begin{prop}
\label{res:tr:low:burn-in}
	For
		all $\alpha < 1$
	and
		all $\sigma \ge \sigma_*(\alpha)$,
	setting $t \cq 20 \log n$,
	\[
		\mfs &\cq \brb{ e \in \Omega \midb \tfrac1n \abs{ \bra{ j \mid e_j > \tfrac13(2\alpha + 1) K } } \le 2 \min\brb{ \tfrac18 (1/\alpha - 1), \: 10^{-3} } }
	\quad	\text{and}
	\\
		\mfb &\cq \brb{ e \in \Omega \midb \prb[e]{ X^s \in \mfs \: \forall \, s \le n } \ge 1 - 1/n },
	\]
	for all $K$ and $n$ sufficiently large (depending on $\alpha$ and $\eps$),
	if $X \sim \DAR_n^\sigma(\alpha, K)$, then
	\[
		\MAX{x \in \Omega} \,
		\prb[x]{ X^t \notin \mfb^2 }
	\le
		4 / n
	=
		\oh1.
	\]
\end{prop}

\begin{Proof}
Recall that $\beta(f,g) = \alpha( 1 + 2 f(1-g) )$ is the effective traffic intensity.

In the trunk reservation system $\DAR_n^\sigma(\alpha,K)$, we have independence between links while the arrival rate remains unchanged.
Thus, analogously to \cref{res:prelim:stoch-dom:cor}, we may stochastically dominate $X \lesssim F \sim \ET(\alpha, \tfrac32\alpha, \sigma, K)^n$, recalling that $\beta(f,g) \le \beta(\tfrac12,\frac12) = \tfrac32 \lambda
$:
\begin{itemize}[noitemsep, topsep = \smallskipamount, label = \bcdot]
	\item 
	the links which are full (ie have $K$ calls on them) do not have any arrivals;
	
	\item 
	the links with number of call on them in $[K-\sigma,K)$ have arrivals at rate $\lambda = \alpha K$;
	
	\item 
	the links with fewer than $K-\sigma$ calls on them have arrivals at rate at most $\tfrac32\alpha$.
\end{itemize}
We have now dominated $X \lesssim F$; write $\Pi_F$ for the invariant distribution of $F$.

By considering the two cases with $\alpha < 1$ in \S\ref{sec:tr:trunked-erlang}, ie $\alpha \le \beta < 1$ and $\alpha < 1 \le \beta$, we see that the invariant measure of $K$ under $\ET(\alpha, \tfrac32 \alpha, \sigma, K)$ is at most $\tilde\alpha^\sigma$ where $\tilde\alpha = \tfrac12(\alpha + 1) < 1$, for $K$ sufficiently large.
Let $\eps > 0$ and suppose that
\[
	\sigma
\Quad{is such that}
	\tilde\alpha^\sigma \le \eps
\Quad{ie}
	\sigma \ge \log_{\tilde\alpha}(\eps) = \log(1/\eps)/\log(1/\tilde\alpha);
\]
Let $\eps \cq \min\bra{ \tfrac18 (1/\alpha-1), 10^{-3} }$.
Define
\[
	\mfs_0 &\cq \brb{ e \in \Omega \midb \varphie \le 2 \eps }
\quad	\text{and}
\\
	\mfb_0 &\cq \brb{ e \in \Omega \midb \prb[e]{ X^s \in \mfs_0 \: \forall \, s \le n } \ge 1 - 1/n }.
\]
Then $\Pi_F(\mfb_0) \ge 1 - 1/n$, using the concentration of the Binomial (as done before, eg \cref{res:high:claim:a>sqrt2:inv-msre}).

By \cref{res:tr:trunked-erlang:product-mixing}, the $1/n$ mixing time of an $\ET(\alpha, \beta, \sigma, K)^n$ system is at most $T \cq 10 \log n$ for any $\alpha, \beta, \sigma$ and $K$ providing $\alpha \le \beta$ and $K$ and $n$ are sufficiently large.
So by running the stochastic domination $X \lesssim F$ for a time $T$, we see that
\(
	\pr{X^T \in \mfb_0}
\ge
	1 - 2/n
\)
as $\Pi_F(\mfb_0) \ge 1 - 1/n$.

For all $x \in \mfs_0$, we have $\beta(f,g) \le \alpha(1+4\eps)$.
By choice of $\eps$, we have $\alpha(1+4\eps) \le \tfrac12(\alpha+1)$.
Then while $X$ is in $\mfs_0$ we can stochastically dominate $X \lesssim E \sim \Er(\tilde\alpha,K)^n$, which is a subcritical system since $\tilde\alpha = \tfrac12 (\alpha + 1) < 1$.
(We could dominate by a product of Erlang links with trunk reservation, but since the link is subcritical we may as well just use the easier formula for the~usual~Erlang~system.)

Given that we have stochastically dominated by a subcritical system, by running for a further $T$ time units, we can now get the proportion of links with more than $\tfrac13 (2 \alpha + 1) K$ calls on them as low as desired, in particular smaller than $2 \eps$; cf \cref{res:low:burnin}.

Recall that $\eps = \min\bra{\tfrac18(1/\alpha-1), 10^{-3}}$ and $\tilde \alpha = \tfrac12 (\alpha + 1)$.
Finally, it is not difficult to check
that
\(
	\log(1/\eps) / \log(1/\tilde \alpha)
\le
	\rbr{ - 2 \log(1 - \alpha) + 14 } / (1 - \alpha)
=
	\sigma_*(\alpha)
\)
This completes the proof.
\end{Proof}

\begin{rmkt*}
	We could formulate \cref{res:tr:low:burn-in} in the following way:
		``for
			all $\alpha < 1$,
			all $\eps \in (0,1)$
		and
			all $\sigma \ge f(\alpha,\sigma)$,
		...''
	where $f$ is some explicit function.
	We have not done this in order to save additional technicalities; instead we simply chose an $\eps \in (0,1)$ which is suitable for later.
\end{rmkt*}

The set $\mfs$ is now (in essence) the same as that in \S\ref{sec:low}, and hence we may then use exactly the same variable length coupling argument as in \S\ref{sec:low:var-coupling} to prove the following result.

\begin{prop}
	For
		$\eps \cq \min\bra{\tfrac18(1/\alpha - 1), 10^{-3}}$
	and
		all $K$ sufficiently large,
	in the scenario of \cref{res:var-length-coup},
	using the coupling from \cref{def:tr:coupling} and $\mfs$ defined above,
	there exists a stopping time $\tau$ with
	$\gamma \le \tfrac23$ and $M \cq 6$ a valid choice.
\end{prop}

\begin{Proof}[Proof References]
\renewcommand{\qedsymbol}{\ensuremath{\triangle}}
	See the results of \S\ref{sec:low:var-coupling}; specifically,
	see
		\cref{def:low:st} for the definition of the stopping time $\tau$
	and
		\cref{res:low:gamma,res:low:M} for the bounds on $\gamma$ and $M$, respectively.
\end{Proof}

The proof of \cref{res:low} given in \S\ref{sec:low:conclusion} applies here to prove \cref{res:tr:low}.

\begin{Proof}[Proof of \cref{res:tr:low}]
Noting that $\gamma^{M} = (2/3)^6 \le e^{-1/18}$, which is the same bound as we used for the original slow arrivals regime, by the same argument we find that
\[
	\prb{ X^t \neq Y^t, \, \msg[0,t] } \le 1/n
\Quad{for all}
	t \ge 40 \log n.
\]
Our burn-in phase is of length $20 \log n$, and so we have
\[
	\prb{ X^t \neq Y^t } \le 9 / n
\Quad{for all}
	t \ge 60 \log n.
\]

\smallskip

The final claim follows from the fact that we stochastically upper bounded the original system by a subcritical Erlang system and the same result holds for such a system.
\end{Proof}

\subsection{Fast Mixing with Fast Arrivals: $\alpha > 1$ and High-Blocking}
\label{sec:tr:high}

In this subsection
we consider the `fast arrivals' regime, ie $\alpha > 1$.
We always use the coupling \mbpcoup from \cref{def:tr:coupling}.
The aim of is to prove the following theorem.
Recall~that
\[
	\sigma_*(\alpha)
=
	\log 4 / \log \alpha
\Quad{for}
	\alpha > 1.
\]

\begin{thm}
\label{res:tr:high}
	For all $\alpha > 1$,
	all $\sigma \ge \sigma_*(\alpha)$
	and all $K$ and $n$ sufficiently large,
	if $X$ and $Y$ are $\DAR_n^\sigma(\alpha,K)$ systems, then
	under the coupling \mbpcoup,
	for all $t \ge 20 \log n$,
	we have
	\[
		\MAX{(x,y) \in \Omega^2}
		\prb[x,y]{ X^t \neq Y^t }
	\le
		4/n
	=
		\oh1.
	\]
	Thus
	for all $K$ sufficiently large,
	all $\eps \in (0,1)$
	and all $n$ sufficiently large,
	we have
	$\tmix(\eps) \le 20 \log n$.
\end{thm}

Again, the proof of this proposition will be very similar to that of \cref{res:high}, except that instead of using a careful sequence of stochastic dominations for the burn-in the trunk reservation will allow us to do one, even more trivial than in the last section, stochastic domination.

We now explain how to do the burn-in phase.
Consider two $\DAR_n^\sigma(\lambda,K)$ systems, $X$ and $Y$.

\begin{prop}
	For all $\alpha > 1$,
	all $\eps \in (0,1)$
	and all $\sigma \ge \log(1/\eps)/\log\alpha$,
	setting
	$t \cq 10 \log n$,
	\[
		\mfs &\cq \brb{ e \in \Omega \midb g(e) \ge 1 - 2 \eps }
	\quad	\text{and}
	\\
		\mfb &\cq \brb{ e \in \Omega \midb \prb[e]{ X^s \in \mfs \: \forall \, s \le n } \ge 1 - 1/n },
	\]
	for all $K$ and $n$ sufficiently large (depending on $\alpha$ and $\eps$),
	we have
	\[
		\prb{ (X^t,Y^t) \notin \mfb^2 }
	\le
		4/n.
	\]
\end{prop}

\begin{Proof}
We may always ignore the rerouting and lower bound $X \gtrsim E \sim \Er(\alpha,K)^n$; see \cref{res:prelim:stoch-dom:cor}.
Writing $\pi_E$ for the invariant distribution of $E$, we have $\pi_E(K) \ge 1 - 1/\alpha$, and $\pi_E(K) \to 1 - 1/\alpha$ as $\Kinf$.
Also, by comparison with a biased simple random walk, we see that the invariant distribution decays exponentially: as $\Kinf$, we have
\[
	\pi_E(K-j) \to \alpha^{-j}(1-1/\alpha)
\Quad{and}
	\pi_E\rbb{\bra{K-j, ..., K}} \to 1 - \alpha^{-(j+1)}
\Quad{for}
	\text{fixed $j$}.
\]
In particular, for $\sigma \ge \log(1/\eps) / \log \alpha$, we get
\[
	\pi_E\rbb{\bra{K-\sigma,...,K}} \ge 1 - \eps
\]
for $K$ sufficiently large.
Hence $\pi_E(\mfb) \ge 1 - 1/n$ for $n$ sufficiently large, using the concentration of the Binomial and of the Poisson distributions (as done in previous sections).

By \cref{res:prelim:product-mixing}, the $1/n$ mixing time of an $\Er(\alpha,K)^n$ system is at most $10 \log n$.
Hence taking our burn-in time to be $10 \log n$, we see that $X$ is in $\mfb$ at this time with probability at least $1 - 2/n$.
Extending this to two systems, this probability becomes at least $1 - 4/n$.
	%
\end{Proof}

When $X$ is in $\mfs$, the probability that a reroute is accepted is $(1-g)^2 \le 4 \eps^2$. If we take $\eps \cq \tfrac14$, then this is at most $\tfrac14$, which is strictly less than $\tfrac12$.
With this choice of $\eps$, for $\alpha > 1$, we have $\log(1/\eps) / \log \alpha = \log 4 / \log \alpha = \sigma_*(\alpha)$.
Given this, we then use exactly the same argument as in \S\ref{sec:high:var-coupling} to prove the following result.

\begin{prop}
	For
		$\eps \cq \tfrac14$
	and
		all $K$ sufficiently large,
	in the scenario of \cref{res:var-length-coup},
	using the coupling from \cref{def:tr:coupling} and $\mfs$ defined above,
	there exists a stopping time $\tau$ with
	$\gamma \le \tfrac45 \le e^{-1/5}$ and $M \cq 3$ a valid choice.
\end{prop}

\begin{Proof}
In proving \cref{res:high:gamma,res:high:M}, we had a parameter $\xi > 0$ and we showed that we have $\gamma \le 1 - \xi$ and may take $M \cq \log(2/\xi)$.
We now recall how $\xi$ was defined:
	it was used in the definition of $\mfs$, in \cref{res:high:burnin}.
In particular, using our current definition of $\mfs$, we can define $\xi$ by $\varphi_c + \xi = 1 - 2 \eps$, where we recall that $\varphi_c = 1 - 1/\sqrt2$, so $\xi \cq 1/\sqrt2 - 2 \eps$. Taking $\eps \cq \tfrac14$, we get $\xi \ge \tfrac15$.
Hence we have $\gamma \le \tfrac45 \le e^{-1/5}$, and $\log(2/\xi) \le \log(10) \le 3$ so we may take $M \cq 3$.
\end{Proof}

The proof of \cref{res:high} given in \S\ref{sec:high:conclusion} applies here with $\xi \cq \tfrac15$ and $M \cq 3$ to prove \cref{res:tr:high}---except now $\xi$ and $M$ do not depend on $\alpha$, so in fact we apply the proof of \cref{res:low} given in \S\ref{sec:low:conclusion} in an (almost) identical way as we did for the low-blocking case.

\begin{Proof}[Proof of \cref{res:tr:high}]
This proof is word-for-word the same as for the low-blocking regime, ie the proof of \cref{res:tr:low} (ie slow arrivals with trunk reservation), except that where we considered time $40 \log n$ and $20 \log n$, we replace these with $10 \log n$; so the total time is $20 \log n$.
\end{Proof}


\section*{Acknowledgements}
\addcontentsline{toc}{section}{Acknowledgements}

\noindent
The question of studying mixing times for this model was originally raised by Nathana\"el Berestycki.
I would like to thank Perla Sousi, my PhD supervisor, for reading this paper and giving lots of constructive feedback.
I would also like to thank Frank Kelly,
for numerous helpful discussions on this work and related \emph{stochastic networks} discussions.
He introduced me to the topic through his Cambridge Part III lecture course and his book \cite{KY:stoc-net-book} with Elena Yudovina; I have become thoroughly interested in the topic as a result.

I also thank the anonymous referee for helpful comments which improved the clarity and presentation of the paper.
They also alerted me to the analogous metastable behaviour exhibited by the Chayes--Machta dynamics in the random cluster model and to the references \cite{BS:rc-meta,GLP:rc-meta:jour}.

\renewcommand{\bibfont}{\sffamily}
\printbibliography[heading=bibintoc]

\appendix
\counterwithout{thm}{subsection}
\counterwithin{thm}{section}
\section{Mixing Time of an Erlang Link or System}
\label{app:er-mix}

The primary aim of this section is to determine the mixing time of a product Erlang system, $\Er(\beta, K)^n$, which we denote $\tmix^{\beta,K;n}(\cdot)$.
To this end, first we find the mixing time of a single Erlang link, $\Er(\beta, K)$, ie $\tmix^{\beta,K;1}(\cdot)$.
We also analyse systems of trunk-reserved Erlang links, ie of $\ET(\alpha, \beta, \sigma, K)$; we denote the corresponding mixing times $\tmix^{\alpha,\beta,\sigma,K;n}(\cdot)$.

Recall that the trunk-reserved link is a generalisation of the standard link: taking $\alpha = \beta$ and $\sigma = 0$ reduces the $\ET$-link to an $\Er$-link.
It thus suffices to consider the trunk-reserved links.

\begin{lem}
\label{res:er-mix:single-link:upper}
	For
		all $\alpha, \beta \in (0,\infty)$ with $\alpha \le \beta$,
		all $\sigma, K \in \mbn_0$ with $K \ge \max\bra{\sigma,1}$
	and
		all $\eps \in (0,1)$,
	we have
	\[
		\tmix^{\alpha,\beta,\sigma,K;1}(\eps) \le \log K + \log(1/\eps).
	\]
	Moreover, since this holds for all $\eps$ and $K$, we may let $\eps$ depend on $K$.
\end{lem}

\begin{Proof}
Consider two links $X = \Xt$ and $Y = \Yt$. Use the `natural' coupling:
	couple the arrivals;
	pair calls where possible and couple their departures and let the `extra' calls depart independently.
This is a \textit{coalescent coupling} in the sense that if $X_t = Y_t$ then $X_s = Y_s$ for all $s \ge t$.
Write $\tau_0 = \inf\bra{ t \ge 0 \mid X_t = Y_t }$.
Then
\[
	\prb[x,y]{X_t \neq Y_t}
=
	\prb[|x-y|]{ \tau_0 > t }.
\]
If $X_t,Y_t < K$, then arrivals leave $X-Y$ unchanged and the departure of an `extra' call decreases $|X-Y|$ by 1.
If $X_t = K > Y_t$, then arrivals also decrease $|X-Y|$ by 1; we upper bound $\tau_0$ by ignoring this.
Hence we can stochastically dominate $\tau_0$ as
\[ \textstyle
	\tau_0 \lesssim \sum_{i=1}^{\abs{x-y}} \mce_i
\Quad{where}
	\mce_i \ \text{is an independent exponential-$i$ for each $i$}.
\]
Determining the quantity on the right-hand side is the coupon-collector problem. The worst-case is clearly when $\abs{x-y} = K$.
Taking $t = \log K + \log(1/\eps)$ and using the union bound, we obtain
\[
	\prb[x,y]{X_t \neq Y_t}
\le
	K e^{-t}
=
	K e^{-\log K - \log(1/\eps)}
=
	\eps.
\]
Hence
\(
	\tmix^{\beta,K;1}(\eps) \le \log K + \log(1/\eps),
\)
using the coupling representation of TV distance.
\end{Proof}

As a corollary of this estimate, we can bound the relaxation time, which we denote $\trel^{\alpha,\beta,\sigma,K;1}$.

\begin{cor}
\label{res:er-mix:single-link:gap}
	For
		all $\alpha, \beta \in (0,\infty)$ with $\alpha \le \beta$,
		all $\sigma, K \in \mbn_0 $ with $K \ge \max\bra{\sigma,4}$
	and
		all $\eps \in (0,1)$,
	we have
	\[
		\trel^{\alpha,\beta,\sigma,K;1} \le 5 (\beta + 1)K.
	\]
\end{cor}

\begin{Proof}
By \cite[Theorem 12.5]{LPW:markov-mixing},
for discrete-time, reversible chains,
for all $\eps \in (0,1)$,
we have
\[
	\tmix(\eps) \ge (\trel - 1) \log\rbb{1/(2 \eps)},
\]
where $\tmix(\cdot)$ is the mixing and $\trel$ the relaxation time.
To apply this to $\Er(\beta,K)$, we must first discretise the chain.
Define the discrete-time process $(Z^m)_{m \in \mbn_0}$ via the following step distribution.
\begin{quote}
	Draw $B \sim \Bern\rbr{ \beta/(\beta + 1) }$.
	
	\begin{itemize}[noitemsep]
		\item 
		If $B = 0$, then select a slot uniformly amongst all $K$ and set this slot to be empty.
		
		\item 
		If $B = 1$, then add a call if the link if there are at least $\sigma+1$ circuits free.
		Otherwise, draw $B' \sim \Bern(\alpha/\beta)$ and add a call if $B' = 1$ and the link is not full.
	\end{itemize}
\end{quote}
Let $(S_m)_{m \in \mbn_0}$ be a rate-$(\beta + 1)K$ Poisson process with $0 \cq S_0 < S_1 < \cdots$.
Define the continuous-time process $Y \cq (Y^t)_{t\ge0}$ by setting $Y^t \cq Z^m$ for $t$ with $S_m \le t < S_{m+1}$.

Analogously to \cref{res:prelim:discrete} we have $Y \sim \ET(\alpha, \beta, \sigma, K)$.
Applying the above mixing--relaxation time inequality with $\eps \cq 1/K$ we deduce that $\trel^{\alpha,\beta,\sigma,K;1} \le 5 (\beta + 1)K$ when $K \ge 4$.
Note that there is a time-change of $(\beta + 1)K$ between the discrete- and continuous-time chains.
\end{Proof}


We use this to bound the mixing of a product Erlang system, $\ET(\alpha,\beta,\sigma,K)^n$.

\begin{cor}
\label{res:er-mix:multi-link:mixing}
	For
		all $\alpha, \beta \in (0,\infty)$ with $\alpha \le \beta$,
		all $\sigma, K \in \mbn_0$ with $K \ge \max\bra{\sigma,1}$
	and
		all $n$ sufficiently large,
	we have
	\[
		\tmix^{\alpha,\beta,\sigma,K;n}(1/n) \le 10 \log n.
	\]
\end{cor}

\begin{Proof}
A simple observation of the proof shows that the mixing result of \cite[Theorem 20.7]{LPW:markov-mixing} applies when $\eps$ depends on $n$. Applying this and the bound of \cref{res:er-mix:single-link:gap}, we obtain our claim.
(Note that there is a time-change of $(\beta + 1)K$ between the discrete- and continuous-time chains.)
\end{Proof}

\section{Deferred Proofs of Lemmas from Fast Mixing Analysis}
\label{app:fast}

In this appendix we give the deferred proofs of \cref{res:low:erlang_inv,res:low:high_inv-msre,res:high:claim:a>sqrt2:inv-msre,res:high:claim:a<sqrt2:inv-msre}.

\begin{Proof}[Proof of \cref{res:low:erlang_inv}]
It is well-known that an unbounded queue $\Er(\beta,\infty)$, ie $M/M/1$ queue with arrival rate $\beta$, has invariant distribution $\Poi(\beta)$.
The invariant distribution of $\Er(\beta,K)$ is simply this but conditioned to be at most $K$.
Let $X \sim \Poi(\beta)$ and $\xi \in (0,1)$.
Then
\[
	\pi_{\beta,K}\rbb{ [\xi K, K] }
=
	\pr{ X \ge \xi K \mid X \le K }
\le
	\pr{ X \ge \xi K } / \pr{ X \le K }.
\]
The Poisson distribution is known to concentrate.
Thus if $\beta < \xi < 1$ are independent of $K$, then
\[
	\LIM{\Kinf} \, \pr{ X \ge \xi K } = 0
\Quad{and}
	\LIM{\Kinf} \, \pr{ X \le K } = 1
\]
as $\ex{X} = \beta K$.
This shows that
\(
	\pi_{\beta,K}\rbr{ [\xi K, K] }
\to
	0
\)
as $\Kinf$ for any $\xi \in (\beta, 1)$.
\end{Proof}

\begin{rmkt*}
	The Poisson distribution concentrates with exponential large deviations. Thus
	\[
		- \log \pi_{\beta,K}\rbb{ [\xi K, K] }
	\asymp
		K
	\quad
		\text{as $\Kinf$}.
	\qedhere
	\]
\end{rmkt*}

\begin{Proof}[Proof of \cref{res:low:high_inv-msre}]
Recall that $p(f)$ is the equilibrium probability that an $\Er(\beta,K)$ link with $\beta \cq \alpha(1 + 2f(1-f))$ is full, if $\beta \ge 1$.
By definition of $k$, for all $j < k$, $\Er(\alpha_j,K)^n$ is a supercritical Erlang system, ie one with $\alpha_j \ge 1$.
Also, recall that $f_i = p(f_{i-1}) + 2\eta$ for all $i$ (with $\eta > 0$).

Since $\Er(\alpha_j,K)^n$ is a product system, we may write $\Pi_j = \pi_j^n$, where $\pi_j$ is the invariant distribution of a single $\Er(\alpha_j,K)$ link. For each $i$, in the limit as $\Kinf$, we have
\[
	\pi_{i-1}(K)
\to
	1 - 1/\alpha_{i-1}
=
	p(f_{i-1})
=
	f_i - 2 \eta.
\]
Hence $\pi_{i-1}(K) \le f_i - \eta$ for sufficiently large $K$.
Now let $Y_{i-1} \sim \Pi_{i-1}$. We then have
\[
	\tfrac1n \absb{ \bra{ j \mid Y_{i-1, j} = K } }
\sim
	\Bin\rbb{n, \pi_{i-1}(K)}
\lesssim
	\Bin(n, f_i - \eta),
\]
by the independence of the links.
Consider $i \ne k+1$ first.
Hence, by Hoeffding's inequality,
\[
	\Pi_{i-1}(\mfs_i^c)
=
	\prb{ \tfrac1n \absb{ \bra{ j \mid Y_{i-1, j} = K } } > f_i }
\le
	\prb{ \Bin(n, f_i - \eta) > f_i n }
\le
	\expb{ - 2 \eta^2 n }.
\]

Using \cref{res:prelim:discrete}, Poisson concentration and the union bound we deduce that
\[
	K n^2 \expb{ - 2 \eta^2 n }
&\gtrsim
	\prb[\Pi_{i-1}]{ E_{i-1}^s \in \mfs_i \: \forall \, s \le n }
\\&=
	\sumt{e \in \Omega}
	\Pi_{i-1}(e) \cdot \prb[e]{ E_{i-1}^s \in \mfs_i \: \forall \, s \le n }
\ge
	\Pi_{i-1}\rbb{(\mfb_i)^c} \cdot 1/n,
\]
since, by definition of $\mfb_i$, if $e \notin \mfb_i$ then the probability in the summation is at least $1/n$.
This establishes the claim for a fixed $i \le k$.

\smallskip

Consider now $i = k+1$.
By construction and \cref{res:low:erlang_inv}, we have $\sum_{\ell = \xi K}^K \pi_k(\ell) \le \eps$.
An analogous argument to that used above now applies.

\smallskip

Since $k$ is a fixed number, depending only on $\alpha$, not on $K$ or $n$, we can choose $K$ and $n$ large enough so that the results hold for all $i = 1,...,k,k+1$ simultaneously.
\end{Proof}

\begin{Proof}[Proof of \cref{res:high:claim:a>sqrt2:inv-msre}]
We have $\alpha > \sqrt2$ and $p(0) - 2 \xi \ge \varphi_c + 2 \xi$, and
\[
	\mfs = \brb{ e \in \Omega \midb \varphi(e) \ge \varphi_c + 2 \xi }
\Quad{and}
	\mfb = \brb{ e \in \Omega \midb \prb[e]{ X^s \in \mfs \: \forall \, s \le n } \ge 1 - 1/n }.
\]

Let $\pi_0$ be the invariant distribution of a single $\Er(\alpha,K)$ link; so $\Pi_0 = \pi_0^n$.
We then have
\[
	\pi_0(K) \ge p(0) \ge \varphi_c + 4 \xi.
\]
Now let $Y_0 \sim \Pi_0$. We then have
\[
	\tfrac1n \absb{ \bra{ j \mid Y_{0,j} = K } }
\sim
	\Bin\rbb{n, \pi_0(K)}
\gtrsim
	\Bin\rbb{n, \varphi_c + 4 \xi},
\]
by the independence of links.
Hence, by Hoeffding's inequality,
\[
	\Pi_0(\mfs^c)
&
=
	\prb{ \tfrac1n \absb{ \bra{ j \mid Y_{0,j} = K } } > \varphi_c + 2 \xi }
\\&
\le
	\prb{ \Bin\rbb{n,\varphi_c + 4 \xi} < (\varphi_c + 2 \xi)n }
\le
	\expb{ - 8 \xi^2 n }.
\]
Using Poisson concentration for the number of updates to the system and the union bound,
we now deduce the claim for $\Pi_0(\mfb)$ in exactly the same was as we did in \cref{res:low:high_inv-msre}.
\end{Proof}

\begin{Proof}[Proof of \cref{res:high:claim:a<sqrt2:inv-msre}]
Consider first $i \ne 0$.
Recall that, by definition of $k$, we have $f_i \le \varphi_c$ for all $i < k$, $\alpha_i \ge \alpha > 1$ for all $i$ and $f_i = \min\bra{p(f_{i-1}) - \eta, \tfrac13}$.
For each $i \ne 0$, we then have
\[
	\pi_{i-1}(K)
\ge
	1 - 1/\alpha_{i-1}
=
	p(f_{i-1})
\ge
	f_i + \eta.
\]
Now let $Y_{i-1} \sim \Pi_{i-1}$.
We then have
\[
	\tfrac1n \absb{ \bra{ j \mid Y_{i-1,j} = K } }
\sim
	\Bin\rbb{n, \pi_{i-1}(K)}
\gtrsim
	\Bin(n, f_i + \eta),
\]
by independence of links.
Hence, by Hoeffding's inequality,
\[
	\Pi_{i-1}(\mfs_i^c)
=
	\prb{ \tfrac1n \absb{ \bra{ j \mid Y_{i-1,j} = K } } < f_i }
\le
	\prb{ \Bin(n, f_i + \eta) < f_i n }
\le
	\expb{ - 2 \eta^2 n }.
\]
Using \cref{res:prelim:discrete}, Poisson concentration and the union bound we deduce the claim for a fixed $i \ge 1$ in exactly the same was as we did in \cref{res:low:high_inv-msre}.

\smallskip

Consider now $i = 0$.
Recall that $p(f) \le 1 - \tfrac23/\sqrt2 < \tfrac35 < \tfrac23$, and hence $\pi_{-1}(K) < \tfrac35$ for $K$ sufficiently large.
The same argument as used above now applies.

\smallskip

Since $k$ is a fixed number, depending only on $\alpha$, not on $K$ or $n$, we can choose $K$ and $n$ large enough so that the results hold for all $i = 0,1,...,k$ simultaneously.
\end{Proof}

\section{High-Degree Polynomial Calculations}
\label{app:poly-calcs}

In this section we prove results on the high-degree polynomials from the retries section, \S\ref{sec:retries}.
To make the proofs more natural and easier to read, we change the notation slightly, eg writing $f$ or $g$ for the functions with inputs $x$ or $y$.

%

\begin{lem}[Stationary Points; \cref{res:retry:unique-maxima,res:retry:fsp>1/2}]
	For $\rho \in \mbn$,
	define $f_\rho : \mbr \to \mbr$ by
	\[
		f_\rho(x)
	\cq
		x \rbb{ 1 - \rbr{1 - (1 - x)^2 }^\rho } / (1 - x)
	\Quad{for}
		x \in \mbr \setminus \bra{1}
	\Quad{and}
		f_\rho(1)
	\cq
		0.
	\]
	Then, for each $\rho$, the following hold:
		$f_\rho$ is a polynomial, and hence smooth;
		there exists a unique $x_\rho \in [0,1]$ with
		\(
			f'_\rho(x_\rho) = 0;
		\)
		the extremum $x_\rho$ is a maximum and lies in $[\tfrac12, 1)$. 
\end{lem}

\begin{Proof}
Abbreviate $f \cq f_\rho$.
First, observe that the numerator in $f$ (which is a polynomial) has a zero at $x = 1$; thus $f$ is a polynomial.
Expanding around $1$ shows that $1$ is a double root, but not a triple root, of the numerator; hence $f(1) = 0$ but $f'(1) \ne 0$.
Similarly, $f(0) = 0$ but $f'(0) \ne 0$.

\smallskip

We now establish the existence and uniqueness of the turning point $x_\rho \in [0,1]$, which must lie in $(0,1)$ as $f'(0) \ne 0 \ne f'(1)$.
It is then convenient to reparametrise by $y \cq 1 - x$:
\[
	g(y)
\cq
	f(1 - y)
=
	(1 - y) \rbb{ 1 - \rbr{ 1 - y^2 }^\rho } / y
\Quad{for}
	y \in \mbr \setminus \bra{0}
\Quad{and}
	g(0)
\cq
	0.
\]
We now differentiate $g$:
\[
	\tilde g(y)
\cq
	g'(y) \cdot y^2
=
	\rbr{1 - y^2}^\rho \rbr{ 1 + y + 2 \rho y } / (1 + y) - 1
=
	\rbr{1 - y^2}^\rho \rbb{ 1 + \tfrac{2 \rho y^2}{1+y} } - 1.
\]
Note that $\tilde g(0) = 0$ and $g'(1) = \tilde g(1) = -1 < 0$;
by expanding around $0$, we see that $g(y) > 0$ for $y > 0$ sufficiently small (and that $g'(0) = \rho > 0$).
We now show that $\tilde g$ has a unique turning point in $(0,1)$, ie there is a unique $y_* \in (0,1)$ with $\tilde g'(y_*) = 0$.
(Note that $\tilde g$ is a polynomial.)
This implies that $\tilde g$, and hence $g'$ and $f'$, has a unique zero in $(0,1)$.

It remains to show the claim for $\tilde g$.
To do this, we differentiate $\tilde g$:
\[
	\tilde g'(y)
=
	- 2 \rho y \rbr{1 - y^2}^\rho \rbb{ (2\rho-1) y^2 + 2y - 1 } / (1+y).
\]
Hence $\tilde g'(y) = 0$ with $y \in (0,1)$ if and only if $y \in \bra{y_-, y_+} \cap (0,1)$ where
\[
	y_\pm
\cq
	\tfrac1{2\rho-1}\rbb{ -1 \pm \sqrt{ 1 + (2\rho-1)^2 } };
\Quad{note that}
	y_- < 0 < y_+ < 1.
\]
Hence $\tilde g$ has a unique turning point in $(0,1)$.
This completes the proof of existence and uniqueness.

\smallskip

We turn to the last part of the statement: we show that $y_\rho \cq 1 - x_\rho \le \tfrac12$.
For $y \ne 0$, we have
\[
	\tilde g(y) = 0
\Quad{if and only if}
	g'(y) = 0
\Quad{if and only if}
	f'(1 - y) = 0.
\]
Observe that $\tilde g(y) = 0$ if and only if
\[
	(1 - y^2) \rbb{ 1 + \tfrac{2 \rho y^2}{1 + y} }^{1/\rho} = 1.
\]
If $\rho = 1$, then this is satisfied by $y = \tfrac12$. (We already know this from previously, as $f_1(x) = x(1-x)$.)
Numerical calculations show that $y_\rho < \tfrac12$ for $\rho \le 9$.
Suppose now that $\rho \ge 9$ and $y \le \tfrac12$.
Then
\[
	(1 - y^2) \rbb{ 1 + \tfrac{2 \rho y^2}{1 + y} }^{1/\rho} = 1
\Quad{implies that}
	1 \le \tfrac34 (1 + \rho)^{1/\rho}.
\]
It is easy to see that this last inequality is only satisfies for $\rho \in \bra{1, ..., 7}$.
This is a contradiction.
Hence $y_\rho \le \tfrac12$, ie $x_\rho \ge \tfrac12$, for all $\rho \in \mbn$.
This completes the proof.
\end{Proof}

\begin{lem}[Distinct Solutions; \cref{res:retry:num-solns}]
	For $\rho \in \mbn$ and $\alpha > 0$,
	define the polynomial $h_{\rho,\alpha}$ by
	\[
		h_{\rho,\alpha}(x)
	\cq
		x \rbb{ 1 - 2 \rbr{ 1 - (1 - x)^2 }^\rho } + 1 - 1/\alpha
	\Quad{for}
		x \in \mbr.
	\]
	For each $\rho \in \mbn$, define
	\[
		\alpha_c(\rho)
	\cq
		\inf\brb{ \alpha > 0 \midb h_{\rho,\alpha}(x) \ne 0 \ \forall \, x \in [0,1] }.
	\]
	Then,
	for $\alpha \in (\alpha_c(\rho), 1)$,
	there are precisely two (distinct) solutions $x \in (0,1)$ to $h_{\rho,\alpha}(x) = 0$,
	while,
	for $\alpha > 1$,
	there is a unique solution $x \in (0,1)$ to $h_{\rho,\alpha}(x) = 0$.
\end{lem}

\begin{Proof}
For $\rho \in \mbn$ and $x \in \mbr$,
define
\[
	f_\rho(x)
\cq
	x \rbb{ 1 - 2\rbr{ 1 - (1-x)}^\rho };
\Quad{then}
	h_{\rho,\alpha}(x) = f_\rho(x) + 1 - 1/\alpha.
\]
Abbreviate $f \cq f_\rho$.
We show that $f$ has a unique turning point in $(0,1)$, which is a maximum.
Along with the fact that $f(0) = 0$, ie $h_{\rho,\alpha}(0) = 1 - 1/\alpha$, this verifies the claims in the statement.

To establish the existence and uniqueness of the turning point in $(0,1)$ for the polynomial $f$, we differentiate.
First, it is convenient to reparametrise by $y \cq 1 - x$:
\[
	g(y)
\cq
	f(1 - y)
\cq
	(1 - y) \rbb{ 1 - 2\rbr{ 1 - y^2 }^\rho }
\Quad{for}
	y \in \mbr.
\]
Clearly it suffices to prove the claim for $g$ instead of $f$.
We differentiate:
\[
	g'(y)
=
	2 \rbr{ 1 - y^2 }^\rho \rbb{ 1 + (1 + 2 \rho) y } / (1 + y) - 1.
\]
(This is very similar, but not exactly the same, as $\tilde g$ from the previous proof. We apply the same style of analysis.)
Note that $g'(0) = 1 > 0$ and $g'(1) = -1 < 0$.
We now show that there is a unique $y_* \in (0,1)$ with $g''(y_*) = 0$.
This implies that $g'$ has a unique zero in $(0,1)$, as required.

It remains to show the claim for $g''$.
To do this, we differentiate $g'$:
\[
	g''(y)
=
	-4 \rho (1 - y^2)^{\rho-1} \rbb{ (2\rho+1) y^2 + 2y - 1 } / (1 + y).
\]
Hence $g''(y) = 0$ with $y \in (0,1)$ if and only if $y \in \bra{y_-, y_+} \cap (0,1)$ where
\[
	y_\pm
\cq
	\tfrac1{2\rho+1} \rbb{ -1 \pm \sqrt{ 1 + (2\rho+1)^2 } };
\Quad{note that}
	y_- < 0 < y_+ < 1.
\]
Hence $g'$ has a unique turning point in $(0,1)$.
This completes the proof.
\end{Proof}

\end{document}

%% file: DAR_preamble_biblio.tex
\usepackage[doi=false,isbn=false,url=false,
	hyperref=auto,
	sorting=nyt,
	maxnames=10,
	maxcitenames=3,
	backend=biber,
	texencoding=auto,
	giveninits=true,
	block=space,
]{biblatex}

\DefineBibliographyStrings{english}{%
	andothers = {et al}
} 
 
\DeclareFieldFormat{eprint:arxiv}{%
	\href{https://arxiv.org/abs/#1}{\emph{arXiv:\allowbreak #1}}}

\DeclareFieldFormat{eprint:mrnumber}{%
	\href{http://www.ams.org/mathscinet-getitem?mr=MR#1}{MR\allowbreak #1}}

\DeclareFieldFormat{eprint:jstor}{%
	\href{http://www.jstor.org/stable/#1}{JSTOR\allowbreak #1}}

\DeclareFieldFormat{eprint:customeprint}{%
	\em Available at \upshape \href{https://www.#1}{#1}}

\DeclareFieldFormat{eprint:online}{%
	Available \href{https://#1}{online}}

\DeclareFieldFormat{eprint:inprep}{%
	\emph{In preparation}}

\DeclareFieldFormat{eprint:manual}{%
	\emph{#1}}

\DeclareFieldFormat{eprint:onarxiv}{%
	\emph{#1}}

\DeclareFieldFormat{eprint:toappear}{%
	\mbox{\emph{To appear}}}

\DeclareFieldFormat{eprint:accepted}{%
	\emph{Accepted at {#1}; to appear}}

\DeclareSourcemap{
	\maps[datatype=bibtex]{
		\map{
			\step[fieldsource=mrnumber,		fieldtarget=eprint, final]
			\step[fieldset=eprinttype,		fieldvalue=mrnumber]
		}
		\map{
			\step[fieldsource=arxiv,		fieldtarget=eprint, final]
			\step[fieldset=eprinttype,		fieldvalue=arxiv]
		}
		\map{
			\step[fieldsource=jstor,		fieldtarget=eprint, final]
			\step[fieldset=eprinttype,		fieldvalue=jstor]
		}
		\map{
			\step[fieldsource=customeprint,	fieldtarget=eprint, final]
			\step[fieldset=eprinttype,		fieldvalue=customeprint]
		}
		\map{
			\step[fieldsource=online,		fieldtarget=eprint, final]
			\step[fieldset=eprinttype,		fieldvalue=online]
		}
		\map{
			\step[fieldsource=inprep,		fieldtarget=eprint, final]
			\step[fieldset=eprinttype,		fieldvalue=inprep]
		}
		\map{
			\step[fieldsource=manual,		fieldtarget=eprint, final]
			\step[fieldset=eprinttype,		fieldvalue=manual]
		}
		\map{
			\step[fieldsource=onarxiv,		fieldtarget=eprint, final]
			\step[fieldset=eprinttype,		fieldvalue=onarxiv]
		}
		\map{
			\step[fieldsource=toappear,		fieldtarget=eprint, final]
			\step[fieldset=eprinttype,		fieldvalue=toappear]
		}
		\map{
			\step[fieldsource=accepted,		fieldtarget=eprint, final]
			\step[fieldset=eprinttype,		fieldvalue=accepted]
		}
	}
}

\DeclareFieldFormat{doi}{%
	\href{https://doi.org/#1}{DOI}%
}

\DeclareFieldFormat{url}{%
	\href{https://#1}{\texttt{#1}}%
}

\DeclareFieldFormat{comments}{%
	\emph{#1}%
} 
 
\DeclareFieldFormat
	[article,inbook,incollection,inproceedings,patent,thesis,unpublished]
	{title}{{#1\isdot}}

\DeclareFieldFormat
	[article,book,inbook,incollection,inproceedings,patent,thesis,unpublished, online]
	{date}{\mkbibparens{#1}}
\DeclareFieldFormat
	[article,book,inbook,incollection,inproceedings,patent,thesis,unpublished]
	{volume}{\mkbibbold{#1}}
\DeclareFieldFormat
	{pages}{\mkbibparens{#1}}

\DeclareFieldFormat{edition}{%
	\ifinteger{#1}
	{\mkbibordedition{#1}~\bibstring{edition}\addcomma}
	{#1~\bibstring{edition}\addcomma}}

\renewbibmacro*{volume+number+eid}{%
	\printfield{volume}%
\setunit*{\adddot}%
	\printfield{number}%
\setunit{\space}%
	\printfield{eid}%
}

\DeclareBibliographyDriver{article}{%
	\usebibmacro{bibindex}%
	\usebibmacro{begentry}%
	\usebibmacro{author}%
\setunit{\addspace}%
	\usebibmacro{date}%
\newunit\newblock
	\usebibmacro{title}%
\newunit\newblock
	\printfield{journaltitle}\addperiod%
\setunit*{\addspace}%
	\usebibmacro{volume+number+eid}%
\setunit*{\addspace}\newblock
	\printfield{pages}
\setunit*{\addspace}
	\printfield{comments}%
	\usebibmacro{eprint}%
\setunit*{\addspace}%
	\printfield{doi}%
	\usebibmacro{finentry}%
}

\DeclareBibliographyDriver{online}{%
	\usebibmacro{bibindex}%
	\usebibmacro{begentry}%
	\usebibmacro{author}%
\setunit{\addspace}%
	\usebibmacro{date}%
\newunit\newblock
	\usebibmacro{title}%
\newunit\newblock
	\usebibmacro{eprint}%
	\usebibmacro{finentry}%
}

\DeclareBibliographyDriver{book}{%
	\usebibmacro{bibindex}%
	\usebibmacro{begentry}%
	\usebibmacro{author}%
\setunit{\addspace}\newblock
	\usebibmacro{date}%
\newunit\newblock
	\usebibmacro{title}%
\setunit*{\addspace}\newblock
	\printfield{edition}%
\newunit\newblock
	\printlist{publisher}
\setunit*{\addspace}%
	\printfield{volume}%
\setunit*{\addspace}\newblock%
	\printfield{comments}%
	\usebibmacro{eprint}%
\setunit*{\addspace}
	\printfield{doi}
	\usebibmacro{finentry}%
}

\DeclareBibliographyDriver{thesis}{%
	\usebibmacro{bibindex}%
	\usebibmacro{begentry}%
	\usebibmacro{author}%
\setunit{\addspace}\newblock
	\usebibmacro{date}%
\newunit\newblock
	\usebibmacro{title}.%
\setunit*{\addspace}\newblock
	\space Thesis, \printlist{institution}
	\usebibmacro{eprint}%
	\printfield{comments}%
\setunit*{\addspace}
	\printfield{doi}%
	\usebibmacro{finentry}%
}

\DeclareBibliographyDriver{inproceedings}{%
	\usebibmacro{bibindex}%
	\usebibmacro{begentry}%
	\usebibmacro{author}%
\setunit{\addspace}%
	\usebibmacro{date}%
\newunit\newblock
	\usebibmacro{title}%
\newunit\newblock
	\printfield{booktitle}\addcomma%
\setunit*{\addspace}%
	\printfield{series}\addcomma%
\setunit*{\addspace}\newblock
	\usebibmacro{editor}\addcomma
\setunit*{\addspace}\newblock
	\printlist{publisher}
\setunit*{\addspace}%
	\printfield{volume}%
\setunit*{\addspace}%
	\printfield{pages}
\setunit*{\addspace}
	\usebibmacro{eprint}%
	\printfield{comments}%
\setunit*{\addnnspace}
	\printfield{doi}
	\usebibmacro{finentry}%
}

\DeclareBibliographyDriver{inbook}{%
	\usebibmacro{bibindex}%
	\usebibmacro{begentry}%
	\usebibmacro{author}%
\setunit{\addspace}%
	\usebibmacro{date}%
\newunit\newblock
	\usebibmacro{title}%
\newunit\newblock
	\printfield{booktitle}\addcomma%
\setunit*{\addspace}%
	\printfield{series}\addcomma%
\setunit*{\addspace}\newblock
	\usebibmacro{editor}\addcomma
\setunit*{\addspace}\newblock
	\printlist{publisher}
\setunit*{\addspace}%
	\printfield{volume}%
\setunit*{\addspace}%
	\printfield{pages}
\setunit*{\addspace}
	\usebibmacro{eprint}%
	\printfield{comments}%
\setunit*{\addspace}
	\printfield{doi}
	\usebibmacro{finentry}%
}

\DeclareBibliographyDriver{incollection}{%
	\usebibmacro{bibindex}%
	\usebibmacro{begentry}%
	\usebibmacro{author}%
\setunit{\addspace}%
	\usebibmacro{date}%
\newunit\newblock
	\usebibmacro{title}%
\newunit\newblock
	\printfield{booktitle}\addcomma%
\setunit*{\addspace}%
	\printfield{series}\addcomma%
\setunit*{\addspace}\newblock
	\printlist{publisher}
\setunit*{\addspace}%
	\printfield{volume}%
\setunit*{\addspace}%
	\printfield{pages}
\setunit*{\addspace}
	\usebibmacro{eprint}%
	\printfield{comments}%
\setunit*{\addspace}
	\printfield{doi}
	\usebibmacro{finentry}%
}

\AtEveryBibitem{%
	\clearfield{day}%
	\clearfield{month}%
} 
 
\DeclareNameAlias{sortname}{given-family}
\DeclareNameAlias{default}{given-family}



%% file: DAR_preamble_main.tex
\usepackage[utf8]{inputenc}

\usepackage{geometry}
\usepackage{etoc}

\usepackage{titlesec}
\titleformat{\section}
{\sffamily \Large \bfseries \boldmath}{\thesection}{1em}{}
\titleformat{\subsection}
{\sffamily \large \bfseries \boldmath}{\thesubsection}{1em}{}
\titleformat{\subsubsection}
{\sffamily \normalsize \bfseries \boldmath}{\thesubsubsection}{1em}{}
\titleformat{\paragraph}[runin]
{\sffamily \normalsize \bfseries \boldmath}{\theparagraph}{1em}{}
\titleformat{\subparagraph}[runin]
{\sffamily \normalsize \bfseries \boldmath}{\thesubparagraph}{1em}{}

\usepackage{chngcntr}
\counterwithin{figure}{section}

\usepackage[UKenglish]{babel}
\usepackage{amsmath}
\usepackage{amsthm}
\usepackage{amssymb}
\usepackage{mathrsfs}
\usepackage{mathtools}
\newcommand{\cq}{\coloneqq}
\usepackage{empheq}
\usepackage{wasysym}

\usepackage{cases}

\def\IfAmpersandUseAlign#1#2&#3\EndIfAmpersandUseAlign
{%
	\if\relax\detokenize{#3}\relax
	\begin{equation*}%
	#1%
	\end{equation*}%
	\else
	\begin{align*}%
	#1%
	\end{align*}%
	\fi
}
\def\[#1\]%
{%
	\IfAmpersandUseAlign{#1}#1&\EndIfAmpersandUseAlign
}

\newcommand{\beq}{\begin{equation*}}
\newcommand{\eeq}{\end{equation*}}

\usepackage[nottoc,numbib]{tocbibind}

\numberwithin{equation}{section}

\usepackage{xifthen}

\usepackage[normalem]{ulem}

\usepackage{stackrel}
\usepackage{bm}

\newcommand{\romannumbering}{%
	\renewcommand{\labelenumi}{\upshape(\roman{enumi})}%
	\renewcommand{\theenumi}{\upshape(\roman{enumi})}%
	\renewcommand{\labelenumii}{\upshape(\textit{\alph{enumii}})}%
	\renewcommand{\theenumii}{\upshape(\roman{enumi}.\textit{\alph{enumii}})}
}

\usepackage{enumitem}

\newcommand{\bcdot}{\ensuremath{\bm\cdot}}
{\begin{itemize}[itemsep=0pt,topsep=\smallskipamount,label={\ensuremath{\bm\cdot}}]}
	{\end{itemize}}

\setlist[description]{%
	topsep=0pt,		
	noitemsep,		
	font={\mdseries\itshape},	
}

\newlist{ralist}{enumerate}{2}
\setlist[ralist,1]{
	label	= \upshape(\roman*)}
\setlist[ralist,2]{
	label	= \upshape(\textit{\alph*}),
	ref		= \upshape(\roman{ralisti}.\textit{\alph*})}

\newlist{ranlist}{enumerate}{3}
\setlist[ranlist,1]{
	label=(\roman*) }
\setlist[ranlist,2]{
	label=(\textit{\alph*}),
	ref=(\roman{ranlisti}.\textit{\alph*})	}
\setlist[ranlist,3]{
	label=\arabic*.,
	ref=(\roman{ranlisti}.\textit{\alph{ranlistii}}.\arabic*)	}

\newlist{ranRlist}{enumerate}{4}
\setlist[ranRlist,1]{
	label=(\roman*) }
\setlist[ranRlist,2]{
	label=(\textit{\alph*}),
	ref=(\roman{ranRlisti}.\textit{\alph*})	}
\setlist[ranRlist,3]{
	label=\arabic*.,
	ref=(\roman{ranRlisti}.\textit{\alph{ranRlistii}}.\arabic*)	}
\setlist[ranRlist,4]{
	label=\Roman*.,
	ref=(\roman{ranRlisti}.\textit{\alph{ranRlistii}}.\arabic{ranRlistiii}.\Roman*)	}

\newlist{rnlist}{enumerate}{2}
\setlist[rnlist,1]{
	label=(\roman*) }
\setlist[rnlist,2]{
	label=\arabic*.,
	ref=(\roman{rnlisti}.\arabic*)	}

\newlist{anlist}{enumerate}{3}
\setlist[anlist,1]{
	label=(\textit{\alph{anlisti}}),
	ref=(\textit{\alph{anlisti}})		}
\setlist[anlist,2]{
	label=\arabic{anlistii}.,
	ref=(\textit{\alph{anlisti}}.\arabic{anlistii})	}
\setlist[anlist,3]{
	label=\Roman{anlistiii}.,
	ref=(\textit{\alph{anlisti}}.\arabic{anlistii}.\Roman{anlistiii})	}

\newlist{Rnlist}{enumerate}{2}
\setlist[Rnlist,1]{
	label=\Roman*. }
\setlist[Rnlist,2]{
	label=\arabic*.,
	ref=(\Roman{Rnlisti}.\arabic*)	}

\newlist{Ranlist}{enumerate}{3}
\setlist[Ranlist,1]{
	label=\Roman*.,
	ref=\Roman* 	}
\setlist[Ranlist,2]{
	label=(\textit{\alph*}),
	ref=(\Roman{Ranlisti}.\textit{\alph*})	}
\setlist[Ranlist,3]{
	label=\arabic*.,
	ref=(\Roman{Ranlisti}.\textit{\alph{Ranlistii}}.\arabic*)		}

\usepackage[dvipsnames,hyperref]{xcolor}

\usepackage{manyfoot}
\SetFootnoteHook{\hspace*{-1.8em}}
\DeclareNewFootnote{bl}[gobble]
\newcommand{\blfootnote}[1]{\footnotebl{\sffamily#1}}

\newcommand{\nt}{\addtocounter{equation}{1}\tag{\theequation}}

\newcommand{\Quad}[1]{\quad \text{#1} \quad}
\newcommand{\QUAD}[1]{
	\mathchoice
	{\quad\text{#1}\quad}
	{\text{ #1 }}
	{\text{ #1 }}
	{\text{ #1 }}
}

\renewcommand{\epsilon}{\varepsilon}

\newcommand{\eps}{\epsilon}

\newcommand{\tref}[1]{\text{\ref{#1}}}

\newcommand{\prref} [1]{\mathbb{P}( \tref{#1} )}
\newcommand{\prrefb}[1]{\mathbb{P}\bigl( \tref{#1} \bigr)}

\newcommand{\maxt}[1]{ \textstyle \max_{#1} \displaystyle }

\newcommand{\binomt}[2]{ \textstyle \binom{#1}{#2} \displaystyle }

\usepackage{calc}
\newlength{\halfplusheight}
\newcommand{\MAX}[1]{\mathop{\raisebox{\halfplusheight}{\(\displaystyle\max_{#1}\)}}}
\newcommand{\MIN}[1]{\mathop{\raisebox{\halfplusheight}{\(\displaystyle\min_{#1}\)}}}
\newcommand{\LIM}[1]{\mathop{\raisebox{\halfplusheight}{\(\displaystyle\lim_{#1}\)}}}

\newcommand{\sumt}[2][]{
	\ifthenelse{\isempty{#1}}
	{\textstyle \sum_{#2}      \displaystyle}
	{\textstyle \sum_{#2}^{#1} \displaystyle}
}
\newcommand{\sumd}[2][]{
	\ifthenelse{\isempty{#1}}
	{\sum_{#2}}
	{\sum_{#2}^{#1}}
}

\newcommand{\intt}[2][]{
	\ifthenelse{\isempty{#1}}
	{\textstyle \int_{#2}      \displaystyle}
	{\textstyle \int_{#2}^{#1} \displaystyle}
}

\newcommand{\prodt}[2][]{
	\ifthenelse{\isempty{#1}}
	{\textstyle \prod_{#2}      \displaystyle}
	{\textstyle \prod_{#2}^{#1} \displaystyle}
}

\usepackage{xspace}

\let\originalexp\exp
\let\exp\relax

\DeclareRobustCommand{\exp} [1]{\originalexp(#1)}
\newcommand{\expb} [1]{\originalexp\bigl( #1 \bigr)}

\newcommand{\abs}  [1]{ | #1 | }
\newcommand{\absb} [1]{\bigl| #1 \bigr|}

\newcommand{\rbr} [1]{ ( #1 ) }
\newcommand{\rbb} [1]{\bigl( #1 \bigr)}
\newcommand{\rbB} [1]{\Bigl( #1 \Bigr)}
\newcommand{\rbbb}[1]{\biggl( #1 \biggr)}

\newcommand{\bra} [1]{ \{ #1 \} }
\newcommand{\brb} [1]{\bigl\{ #1 \bigr\}}

\DeclareMathOperator{\diam}{diam}
\newcommand{\mix}{\textnormal{mix}}

\newcommand{\rel}{\textnormal{rel}}

\newcommand{\ST}{{\textnormal{ st }}}

\newcommand{\SP}{\textnormal{SP}}

\newcommand{\DAR}{\textup{DAR}}

\newcommand{\Er}{\textnormal{Er}}
\newcommand{\ET}{\textnormal{ET}}

\newcommand{\Stop}		{{\upshape \bfseries \sffamily Stop}\xspace}
\newcommand{\Repeat}	{{\upshape \bfseries \sffamily Repeat}\xspace}
\newcommand{\Return}	{{\upshape \bfseries \sffamily Return}\xspace}
\newcommand{\Proceed}	{{\upshape \bfseries \sffamily Proceed}\xspace}
\newcommand{\FOR}		{{\upshape \bfseries \sffamily for}\xspace}
\newcommand{\IF}		{{\upshape \bfseries \sffamily if}\xspace}
\newcommand{\END}		{{\upshape \bfseries \sffamily end}\xspace}

\newcommand{\Bin}{\textnormal{Bin}}

\newcommand{\Poi}{\textnormal{Poi}}
\newcommand{\Bern}{\textnormal{Bern}}

\newcommand{\tmix}{t_\mix}
\newcommand{\trel}{t_\rel}

\newcommand{\mbpcoup}{\ensuremath{(\mathbb{P}_{x,y})_{x,y\in\Omega}}\xspace}
\newcommand{\varphie}[1][]{
	\ifthenelse{\equal{}{#1}}
	{\tfrac1n \abs{ \bra{ j \mid e_j = K } } }
	{\tfrac1n \abs{ \bra{ j \mid #1_j = K } } }
}

\newcommand{\Xt}{(X_t)_{t\ge0}}

\newcommand{\Yt}{(Y_t)_{t\ge0}}

\newcommand{\cups}{\cup \cdots \cup}

\newcommand{\floor}[1]{\lfloor #1 \rfloor}

\newcommand{\ceil}[1]{\lceil #1 \rceil}

\newcommand{\midb}{ \, \bigm| \, }

\newcommand{\one}  [1]{\textnormal{\textbf{1}}\{ #1 \}}

\newcommand{\tvb}[1]{\bigl\lVert #1 \bigr\rVert_{\textnormal{TV}}}

\newcommand{\logn}[1][]{
	\ifthenelse{\equal{}{#1}}
	{\log n}
	{(\log n)^{#1}}
}

\newcommand{\loglogn}[1][]{
	\ifthenelse{\equal{}{#1}}
	{\log\log n}
	{(\log\log n)^{#1}}
}

\newcommand{\pr}[2][]{
	\ifthenelse{\equal{}{#1}}
	{\mathbb{P}(#2)}
	{\mathbb{P}_{#1}(#2)}
}
\newcommand{\prb}[2][]{
	\ifthenelse{\equal{}{#1}}
	{\mathbb{P}\bigl( #2 \bigr)}
	{\mathbb{P}_{#1}\bigl( #2 \bigr)}
}
\newcommand{\prB}[2][]{
	\ifthenelse{\equal{}{#1}}
	{\mathbb{P}\Bigl( #2 \Bigr)}
	{\mathbb{P}_{#1}\Bigl( #2 \Bigr)}
}
\newcommand{\prbb}[2][]{
	\ifthenelse{\equal{}{#1}}
	{\mathbb{P}\biggl( #2 \biggr)}
	{\mathbb{P}_{#1}\biggl( #2 \biggr)}
}
\newcommand{\prBB}[2][]{
	\ifthenelse{\equal{}{#1}}
	{\mathbb{P}\Biggl( #2 \Biggr)}
	{\mathbb{P}_{#1}\Biggl( #2 \Biggr)}
}
\newcommand{\prs}[2][]{
	\ifthenelse{\equal{}{#1}}
	{\mathbb{P}\left( #2 \right)}
	{\mathbb{P}_{#1}\left( #2 \right)}
}

\newcommand{\ex}[2][]{
	\ifthenelse{\equal{}{#1}}
	{\mathbb{E}(#2)}
	{\mathbb{E}_{#1}(#2)}
}
\newcommand{\exb}[2][]{
	\ifthenelse{\equal{}{#1}}
	{\mathbb{E}\bigl( #2 \bigr)}
	{\mathbb{E}_{#1}\bigr( #2 \bigr)}
}
\newcommand{\exB}[2][]{
	\ifthenelse{\equal{}{#1}}
	{\mathbb{E}\Bigl( #2 \Bigr)}
	{\mathbb{E}_{#1}\Bigl( #2 \Bigr)}
}
\newcommand{\exbb}[2][]{
	\ifthenelse{\equal{}{#1}}
	{\mathbb{E}\biggl( #2 \biggr)}
	{\mathbb{E}_{#1}\biggl( #2 \biggr)}
}
\newcommand{\exBB}[2][]{
	\ifthenelse{\equal{}{#1}}
	{\mathbb{E}\Biggl( #2 \Biggr)}
	{\mathbb{E}_{#1}\Biggl( #2 \Biggr)}
}

\newcommand{\Varb}[2][]{
	\ifthenelse{\equal{}{#1}}
	{\mathbb{V}\textnormal{ar} \bigl(#2\bigr)}
	{\mathbb{V}\textnormal{ar}_{#1} \bigl(#2\bigr)}
}

\newcommand{\Oh}  [1]{\mathcal{O}( #1 )}

\newcommand{\oh}  [1]{o( #1 )}

\newcommand{\mbn}{\mathbb{N}}

\newcommand{\mbp}{\mathbb{P}}

\newcommand{\mbr}{\mathbb{R}}

\newcommand{\mbz}{\mathbb{Z}}

\newcommand{\mce}{\mathcal{E}}
\newcommand{\mcf}{\mathcal{F}}
\newcommand{\mcg}{\mathcal{G}}

\newcommand{\mcl}{\mathcal{L}}

\newcommand{\msg}{\mathscr{G}}

\newcommand{\mfb}{\mathfrak{B}}

\newcommand{\mfl}{\mathfrak{L}}

\newcommand{\mfs}{\mathfrak{S}}

\newcommand{\mfu}{\mathfrak{U}}

\newcommand{\ninf}{\ensuremath{n\to\infty}\xspace}

\newcommand{\Kinf}{\ensuremath{K\to\infty}\xspace}

\newcommand{\asKinf}{\text{as $\Kinf$}}

\def\IfAmpersandUseAlign#1#2&#3\EndIfAmpersandUseAlign%
	{%
		\if\relax\detokenize{#3}\relax
		\begin{equation*}%
		#1%
		\end{equation*}%
		\else
		\begin{align*}%
		#1%
		\end{align*}%
		\fi
	}
	\def\[#1\]%
	{%
		\IfAmpersandUseAlign{#1}#1&\EndIfAmpersandUseAlign
	}

\usepackage[
	colorlinks,
	pdfauthor = {Sam\ Olesker-Taylor},
	pdfusetitle
]{hyperref}
\hypersetup{
	colorlinks,
	linkcolor	= {blue},
	urlcolor	= {blue},
	citecolor	= {ForestGreen}
}


\newcounter{parentnumber}

\usepackage[capitalise]{cleveref}

\newenvironment{Proof}[1][\proofname]{%
	\proof[\sffamily\bfseries{\boldmath#1}]
}{\endproof}
\newcommand{\qedtriangle}{\renewcommand{\qedsymbol}{\ensuremath{\triangle}}}

\newtheoremstyle{sfsl}
{1\baselineskip}		
{1\baselineskip}		
{\slshape}				
{}						
{\bfseries\sffamily}	
{.}						
{0.5em}					
{\thmname{#1}\thmnumber{ #2}\thmnote{ \textnormal{\sffamily(#3)}}}

\newtheoremstyle{sfup}
{1\baselineskip}		
{1\baselineskip}		
{\upshape}				
{}						
{\bfseries\sffamily}	
{.}						
{0.5em}					
{\thmname{#1}\thmnumber{ #2}\thmnote{ \textnormal{\sffamily(#3)}}}

\theoremstyle{sfsl}

\newtheorem{mainthm}{Theorem}

\newtheorem*{mainthm*}{Theorem}

\newtheorem{thm}  {Theorem}[section]
\newtheorem{lem} [thm]{Lemma}
\newtheorem{prop} [thm]{Proposition}
\newtheorem{cor} [thm]{Corollary}

\newtheorem{defn}[thm]{Definition}
\newtheorem*{defn*}{Definition}

\newtheorem*{rmk*}{Remark}

\newtheorem{alg}[thm]{Algorithm}
\newtheorem*{alg*}{Algorithm}

\newtheorem{clm}[thm]{Claim}
\newtheorem*{clm*}{Claim}

\theoremstyle{sfup}

\newtheorem*{rmkTT}{Remark}
\newenvironment{rmkt*}
	{\pushQED{\qed}\renewcommand{\qedsymbol}{$\triangle$}\rmkTT}
	{\popQED\endrmkTT}

\newtheorem*{rmks*} {Remarks}
\newtheorem*{rmksTT}{Remarks}
\newenvironment{rmkst*}
	{\pushQED{\qed}\renewcommand{\qedsymbol}{$\triangle$}\rmksTT}
	{\popQED\endrmksTT}

\newtheorem*{notation*} {Notation}
\newtheorem*{notationTT}{Notation}
\newenvironment{notationt*}
	{\pushQED{\qed}\renewcommand{\qedsymbol}{$\triangle$}\notationTT}
	{\popQED\endnotationTT}

\theoremstyle{sfsl}

\newtheorem{innercustomgeneric}{\customgenericname}
\providecommand{\customgenericname}{}
\newcommand{\newcustomtheorem}[2]{%
	\newenvironment{#1}[1]
	{%
		\renewcommand\customgenericname{#2}%
		\renewcommand\theinnercustomgeneric{##1}%
		\innercustomgeneric
	}
	{\endinnercustomgeneric}
}

\newcustomtheorem{customthm} {Theorem}
\newcustomtheorem{customprop}{Proposition}
\newcustomtheorem{customlem} {Lemma}
\newcustomtheorem{customcor} {Corollary}
\newcustomtheorem{customres}{\!\!}

\crefname{thm}{Theorem}{Theorems}
\crefname{thmT}{Theorem}{Theorems}
\crefname{mainthm}{Theorem}{Theorems}
\crefname{lem}{Lemma}{Lemmas}
\crefname{lemT}{Lemma}{Lemmas}
\crefname{cor}{Corollary}{Corollaries}
\crefname{corT}{Corollary}{Corollaries}
\crefname{prop}{Proposition}{Propositions}
\crefname{propT}{Proposition}{Propositions}
\crefname{alg}{Algorithm}{Algorithms}
\crefname{algT}{Algorithm}{Algorithms}
\crefname{clm}{Claim}{Claims}

\usepackage{xparse}
\ExplSyntaxOn
\NewDocumentCommand{\mref}{m}{\quinn_mref:n {#1}}
\seq_new:N \l_quinn_mref_seq
\cs_new:Npn \quinn_mref:n #1
{
	\seq_set_split:Nnn \l_quinn_mref_seq { , } { #1 }
	\seq_pop_right:NN \l_quinn_mref_seq \l_tmpa_tl
	( 
	\seq_map_inline:Nn \l_quinn_mref_seq
	{ \ref{##1},\nobreakspace } 
	\exp_args:NV \ref \l_tmpa_tl 
	) 
}
\ExplSyntaxOff

\newtheorem*{note*} {Note}
\newtheorem*{noteTT}{Note}
	\newenvironment{notet*}
	{\pushQED{\qed}\renewcommand{\qedsymbol}{$\triangle$}\noteTT}
	{\popQED\endnoteTT}